\documentclass{amsart}

\usepackage{amsfonts, amssymb, amsmath, amsthm}                               
\usepackage{mathrsfs}
\usepackage{geometry}                                                         
\usepackage{setspace}                                                         
\usepackage{url}                                                              
\usepackage{color}
\usepackage{slashed}
\usepackage{comment}

\newtheorem{theorem}{Theorem} 	      	      	                              
\newtheorem{corollary}[theorem]{Corollary}     	      	      	      	      
\newtheorem{lemma}[theorem]{Lemma}     	       	      	      	      	      
\newtheorem{proposition}[theorem]{Proposition} 	      	      	      	      
\newtheorem{definition}[theorem]{Definition} 	      	      	      	        
\newtheorem{assumption}[theorem]{Assumption}     	      	      	      	    
\newtheorem*{remark}{Remark}                                                  

\numberwithin{equation}{section}                                              
\numberwithin{theorem}{section}                                               
\numberwithin{figure}{section}                                                

\newcommand{\mf}[1]{\mathfrak{#1}}                                            
\newcommand{\mc}[1]{\mathcal{#1}}                                             

\newcommand{\N}{\mathbb{N}}                                                   
\newcommand{\R}{\mathbb{R}}                                                   
\newcommand{\C}{\mathbb{C}}                                                   
\newcommand{\Sph}{\mathbb{S}}                                                 

\newcommand{\zfac}{ \langle\langle \xi \rangle\rangle_\ell}                   
\newcommand{\pd}[2]{\smash{#1^{\scriptscriptstyle (#2)}}}                     

\allowdisplaybreaks
\begin{document}

\title[Critically Weakly Hyperbolic and Singular Systems]{Asymptotics and Scattering for Critically Weakly Hyperbolic and Singular Systems}

\author{Bolys Sabitbek}
\address{School of Mathematical Sciences, Queen Mary University of London, London E1 4NS, United Kingdom\\ and Institute of Mathematics and Mathematical Modeling, Almaty, Kazakhstan}
\email{sabitbek.bolys@ugent.be}

\author{Arick Shao}
\address{School of Mathematical Sciences, Queen Mary University of London, London E1 4NS, United Kingdom}
\email{a.shao@qmul.ac.uk}

\begin{abstract}
We study a very general class of first-order linear hyperbolic systems that both become weakly hyperbolic and contain lower-order coefficients that blow up at a single time $t = 0$.
In ``critical" weakly hyperbolic settings, it is well-known that solutions lose a finite amount of regularity at the degenerate time $t = 0$.
In this paper, we both improve upon the results in the weakly hyperbolic setting, and we extend this analysis to systems containing critically singular coefficients, which may also exhibit significantly modified asymptotics at $t = 0$.

In particular, we give precise quantifications for (1) the asymptotics of solutions as $t$ approaches $0$; (2) the scattering problem of solving the system with asymptotic data at $t = 0$; and (3) the loss of regularity due to the degeneracies at $t = 0$.
Finally, we discuss a variety of applications for these results, including to weakly hyperbolic and singular wave equations, equations of higher order, and equations arising from relativity and cosmology, e.g.\ at big bang singularities.
\end{abstract}

\maketitle

\section{Introduction} \label{sec.intro}

In this article, we study a general class of first-order linear hyperbolic systems,
\begin{equation}
\label{eq.intro_system_pre} \partial_t \mc{U} = \mc{A} ( t, \nabla_x ) \mc{U} + \mc{F} \text{,}
\end{equation}
that both become \emph{weakly hyperbolic} (or \emph{degenerate}) and contain \emph{singular} lower-order coefficients at a time $t = 0$ and in a ``critical" manner.
One purpose in studying \eqref{eq.intro_system_pre} is twofold:
\begin{itemize}
\item We aim to determine \emph{precise asymptotics} at the degenerate and singular time $t = 0$---given data at a time $T > 0$, we wish to derive quantities associated to the corresponding solution $\mc{U}$ of \eqref{eq.intro_system_pre} that remain bounded and attain an asymptotic limit at $t = 0$.

\item We also address the \emph{scattering} problem from $t = 0$---given data for the above asymptotic quantities, we find a unique solution $\mc{U}$ of \eqref{eq.intro_system_pre} attaining this asymptotic data at $t = 0$.
\end{itemize}
Furthermore, it is well-known that solutions of such weakly hyperbolic systems exhibit a loss of regularity at $t = 0$.
Therefore, another objective is to obtain precise quantifications of this regularity loss in our aforementioned asymptotic and scattering theories.

Yet another goal is to obtain such results in an extensively general class of systems.
In particular, the systems \eqref{eq.intro_system_pre} will include, as special cases, not only the degenerate wave equations previously studied in the weakly hyperbolic literature, but also those additionally containing critically singular lower-order coefficients.
In fact, such a mix of weakly hyperbolic and singular behaviors naturally arise in a variety of physical models.
Therefore, as further applications of our main results, we will examine both the scalar wave and linearized Einstein equations on Kasner spacetimes, which serve as prototypical models for big bang singularities in relativity and cosmology.

In contrast to \eqref{eq.intro_system_pre}, equations arising from physics often tend to be nonlinear and hence more complex.
As a result, the last objective of the present paper to serve as a starting point for studying similar phenomena in both semilinear and quasilinear partial differential equations, with the goal of building toward a more detailed understanding of their asymptotic behaviors.

\subsection{Background and Motivations}

To motivate our results, as well as the more precise form of the system \eqref{eq.intro_system_pre} that we will study, we initiate our discussions with some concrete special cases.

\subsubsection{Model Equations}

Consider first, in one spatial dimension, the wave equation
\begin{equation}
\label{eq.intro_wave_1} \partial_t^2 \phi - t^{ 2 \ell } \partial_x^2 \phi + c t^{ \ell - 1 } \partial_x \phi = 0 \text{,} \qquad \ell \in \N \text{,}
\end{equation}
on times $t \geq 0$.
Observe that solutions propagate with characteristic velocities $\pm t^\ell$, which coincide at $t = 0$.
This is often described as \eqref{eq.intro_wave_1} being \emph{weakly hyperbolic} at $t = 0$, signifying that, when written appropriately as a first-order system, the matrix of principal coefficients, while still always possessing real eigenvalues, is no longer diagonalizable at $t = 0$.

It is well-known that this weak hyperbolicity can cause a \emph{loss of regularity at $t = 0$}, \emph{depending on the behavior of lower-order coefficients} in the equation.
Various quantitative relations between lower-order coefficients (e.g.\ before the ``$\partial_x$") and regularity loss have been referred in the literature as \emph{Levi conditions} \cite{levi_1908, levi_1909}; see, for instance, \cite{ivrii_petkov_1974, oleinik_1970}.
Roughly speaking, in the present context:
\begin{itemize}
\item If the lower-order coefficients vanish as $t \searrow 0$ quickly compared to $t^{ 2 \ell }$, then there is no loss of derivatives, and the lower-order terms can be seen as purely perturbative.

\item If the coefficients vanish as $t \searrow 0$ slowly compared to $t^{ 2 \ell }$, then there is infinite derivative loss, and one may only have well-posedness from $t = 0$ in Gevrey spaces \cite{shinkai_1991}.
\end{itemize}
Between the above two regimes, there is a ``critical" rate---precisely the ``$t^{ \ell - 1 }$" in \eqref{eq.intro_wave_1}---for which one loses a finite number of derivatives.
Furthermore, the precise number of derivatives lost depends nontrivially on the coefficients.
Thus, it is natural to ask \emph{what the mechanism is for the regularity loss at $t = 0$}, as well as \emph{provide precise formulas quantifying this loss}.

An early historical example, due to \cite{qi_1958}, concerns the case $\ell = 1$ and $c := 1 + 4 k$, with $k \in \N$:
\begin{equation}
\label{eq.intro_wave_1a} \partial_t^2 \phi - t^2 \partial_x^2 \phi + ( 1 + 4 k ) \partial_x \phi = 0 \text{,} \qquad ( \phi, \partial_t \phi ) |_{ t = 0 } := ( \phi_0, 0 ) \text{.}
\end{equation}
When the initial data at $t = 0$ has vanishing velocity, \eqref{eq.intro_wave_1a} has the explicit solution
\[
\phi ( t, x ) = \sum_{ j = 0 }^k c_{ jk } t^{ 2j } \, \partial_x^j \phi_0 \big( x + \tfrac{ t^2 }{2} \big) \text{,}
\]
for some non-zero constants $c_{jk}$.
In other words, there is an unavoidable loss of $k$ derivatives when one goes from the initial data at $t = 0$ to the solution on $t > 0$.

Another family of simple solutions comes from taking general $\ell \in \N$ and $c := \ell$,
\begin{equation}
\label{eq.intro_wave_1b} \partial_t^2 \phi - t^{ 2 \ell } \partial_x^2 \phi + \ell t^{ \ell - 1 } \partial_x \phi = 0 \text{,} \qquad ( \phi, \partial_t \phi ) |_{ t = 0 } := ( \phi_0, \phi_1 ) \text{,}
\end{equation}
along with general initial data.
Taking spatial Fourier transforms ($x \mapsto \xi$) of \eqref{eq.intro_wave_1b} and solving the resulting initial-value problem, one can derive the following explicit solutions:
\begin{equation}
\label{eq.intro_wave_1bs} \hat{\phi} ( t, \xi ) = \hat{\phi}_0 ( \xi ) \, e^{ - \frac{ \imath t^{ \ell + 1 } \xi }{ \ell + 1 } } + \hat{\phi}_1 ( \xi ) \, e^{ - \frac{ \imath t^{ \ell + 1 } \xi }{ \ell + 1 } } \int_0^t e^\frac{ 2 \imath \tau^{ \ell + 1 } \xi }{ \ell + 1 } \, d \tau \text{.}
\end{equation}
Note the analogue of \eqref{eq.intro_wave_1bs} for $\ell := 0$ is the standard d'Alembert formula for the free wave equation.
However, for $\ell \geq 1$, the stationary phase at $t = 0$ in the integral in \eqref{eq.intro_wave_1bs} introduces a degeneracy.
For \eqref{eq.intro_wave_1bs}, this is a Fresnel integral, which is known to satisfy
\[
\int_0^t e^\frac{ 2 \imath \tau^{ \ell + 1 } \xi }{ \ell + 1 } \, d \tau \lesssim_t \langle \xi \rangle^{ - \frac{1}{ \ell + 1 } } \text{.}
\]
This results in a \emph{fractional} loss of $\frac{ \ell }{ \ell + 1 }$-derivatives when $\phi_1 \not\equiv 0$:
\[
\langle \xi \rangle | \hat{\phi} ( t, \xi ) | \lesssim \langle \xi \rangle | \hat{\phi}_0 ( \xi ) | + \langle \xi \rangle^\frac{ \ell }{ \ell + 1 } | \hat{\phi}_1 ( \xi ) | \text{.}
\]

In fact, solutions of \eqref{eq.intro_wave_1}, for any $\ell$ and $c$, can be expressed using the \emph{confluent hypergeometric function} ${}_1 F_1$.
Indeed, from \cite{taniguchi_tozaki_1980}, the solution of \eqref{eq.intro_wave_1}, with initial data $( \phi_0, \phi_1 )$ at $t = 0$, is
\begin{align}
\label{eq.intro_wave_1s} \hat{\phi} ( t, \xi ) &= e^{ - \frac{ \imath t^{ \ell + 1 } \xi }{ \ell + 1 } } \, {}_1 F_1 \big( \tfrac{ \ell - c }{ 2 ( \ell + 1 ) }, \tfrac{ \ell }{ \ell + 1 }, \tfrac{ 2 \imath t^{ \ell + 1 } \xi }{ \ell + 1 } \big) \hat{\phi}_0 ( \xi ) \\
\notag &\qquad + t e^{ - \frac{ \imath t^{ \ell + 1 } \xi }{ \ell + 1 } } \big( \tfrac{ 2 \imath t^{ \ell + 1 } \xi }{ \ell + 1 } \big)^\frac{1}{ \ell + 1 } \, {}_1 F_1 \big( \tfrac{ \ell - c + 2 }{ 2 ( \ell + 1 ) }, \tfrac{ \ell + 2 }{ \ell + 1 }, \tfrac{ 2 \imath t^{ \ell + 1 } \xi }{ \ell + 1 } \big) \hat{\phi}_1 ( \xi ) \text{.}
\end{align}
From this point of view, the regularity loss can be attributed to singular points of the confluent hypergeometric equation satisfied by the map $z \mapsto {}_1 F_1 ( a, b, z )$ for any parameters $a, b$---namely, the regular point at $z = 0$ and the irregular point at $z = \infty$.
These imply that at leading order, ${}_1 F_1 ( a, b, z )$ behaves like certain powers $z^{ p_0 ( a, b ) }$ near $z = 0$ and $z^{ p_\infty ( a, b ) }$ near $z = \infty$.

Applying this insight to \eqref{eq.intro_wave_1s}, one obtains corresponding powers of $| \xi |$ in the formula for $\smash{ \hat{\phi} ( t, \xi ) }$, from which a careful computation leads to the following estimate for $t > 0$:
\begin{equation}
\label{eq.intro_wave_1scat} \| \phi (t) \|_{ H^s } + \| \partial_t \phi (t) \|_{ H^{ s - 1 } } \lesssim \| \phi (0) \|_{ H^{ s + \frac{ |c| - \ell }{ 2 ( \ell + 1 ) } } } + \| \partial_t \phi (0) \|_{ H^{ s - 1 + \frac{ |c| + \ell }{ 2 ( \ell + 1 ) } } } \text{.}
\end{equation}
Moreover, a similar computation produces the reverse estimate, with a different regularity shift:
\begin{equation}
\label{eq.intro_wave_1asymp} \| \phi (0) \|_{ H^{ s - \frac{ \ell }{ \ell + 1 } } } + \| \partial_t \phi (0) \|_{ H^{ s - 1 } } \lesssim \| \phi (t) \|_{ H^{ s + \frac{ |c| - \ell }{ 2 ( \ell + 1 ) } } } + \| \partial_t \phi (t) \|_{ H^{ s - 1 + \frac{ |c| - \ell }{ 2 ( \ell + 1 ) } } } \text{.}
\end{equation}

\begin{remark}
A similar loss of regularity holds for \emph{inhomogeneous} wave equations.
An explicit model example is given in \cite{mandai_1982}, for which the solution $\phi$ (with trivial initial data) depends on additional derivatives of the given forcing term $f$.
While we will, for the sake of brevity, refrain from discussing inhomogeneous equations further in the introduction, the main results of this paper, and elsewhere in the literature, also apply to equations containing such forcing terms.
\end{remark}

\subsubsection{Weakly Hyperbolic Equations}

With the above in hand, one could next ask whether various generalizations of \eqref{eq.intro_wave_1} exhibit similar behaviors.
Consider, for instance, the wave equation
\begin{equation}
\label{eq.intro_wave_2} \partial_t^2 \phi - t^{ 2 \ell } a (t) \, \partial_x^2 \phi + 2 t^\ell b (t) \, \partial^2_{tx} \phi + t^{ \ell - 1 } c (t) \, \partial_x \phi + \tilde{g} (t) \, \partial_t \phi + \tilde{h} (t) \, \phi = 0 \text{,}
\end{equation}
where $a, b, c, g, h$ are sufficiently smooth functions of (only) $t$, with $a$ also being uniformly positive.
In essence, \eqref{eq.intro_wave_2} has the same leading-order behavior as \eqref{eq.intro_wave_1}, but now with more general coefficients.
One could also study the natural analogues of \eqref{eq.intro_wave_2} in higher dimensional space $\R^d$,
\begin{align}
\label{eq.intro_wave_3} \partial_t^2 \phi - t^{ 2 \ell } \sum_{ i, j = 1 }^d a_{ij} (t) \, \partial_{ x_i x_j }^2 \phi + 2 t^\ell \sum_{ i = 1 }^d b_i (t) \, \partial^2_{ t x_i } \phi \qquad & \\
\notag + t^{ \ell - 1 } \sum_{ i = 1 }^d c_i (t) \, \partial_{ x_i } \phi + \tilde{g} (t) \, \partial_t \phi + \tilde{h} (t) \, \phi &= 0
\end{align}
(where the $a_{ij}$'s have the usual uniform ellipticity property), and again ask whether solutions satisfy estimates with derivative loss similar to \eqref{eq.intro_wave_1scat}--\eqref{eq.intro_wave_1asymp}.

Classical results for general classes of wave equations include \cite{nersesian_1966, oleinik_1970}.
In one spatial dimension, \cite{nersesian_1966} obtained, in the context of \eqref{eq.intro_wave_2}, finite loss of regularity depending on $\ell$, $a$, and $c$.
Similarly, in higher dimensions, \cite{oleinik_1970} derived, in the setting of \eqref{eq.intro_wave_3}, finite regularity loss depending on $\ell$, the $a_{ij}$'s, and the $c_i$'s.
Both \cite{nersesian_1966, oleinik_1970} utilized physical space methods, hence their results only quantified integer derivative losses, which were non-optimal.
(On the other hand, these physical space methods could also treat equations with at least some coefficients also depending on $x$.)

More recent treatments of \eqref{eq.intro_wave_2} and \eqref{eq.intro_wave_3} that more carefully quantified the finite loss of regularity include \cite{dreher_reissig_2000, dreher_witt_2002, ebert_kapp_filho_2009, reissig_1997} (see also \cite{dreher_2002, dreher_witt_2005, ebert_2005} for more general weakly hyperbolic systems).
For the subsequent discussion, we narrow our focus to the results of \cite{dreher_reissig_2000, dreher_witt_2002, reissig_1997}, which applied more refined Fourier methods in order to attain the optimal (fractional) derivative loss \eqref{eq.intro_wave_1scat} when applied to the model equation \eqref{eq.intro_wave_1}, as well as in several other cases.

A basic difficulty in treating wave equations \eqref{eq.intro_wave_2}--\eqref{eq.intro_wave_3} is that one no longer has explicit solutions \eqref{eq.intro_wave_1s}.
Furthermore, solutions of \eqref{eq.intro_wave_2} generally fail to approach a model solution \eqref{eq.intro_wave_1s} in the limit $t \searrow 0$, as it turns out that one must modify the expansions for the model solution---in particular, for ${}_1 F_1$---at (possibly arbitrarily) higher orders.
Despite the above, a number of qualitative properties of the model solutions \eqref{eq.intro_wave_1s} do carry over to this more general setting.

The first key idea is that the nature of the singularities of ${}_1 F_1$ do in fact extend to solutions of \eqref{eq.intro_wave_2}--\eqref{eq.intro_wave_3} and are fundamental to capturing the asymptotic behavior of these solutions.
Note that in \eqref{eq.intro_wave_1s}, the dynamical parts of the solution are written entirely in terms of $t^{ \ell + 1 } \xi$.
With this in mind, it is often natural to analyze \eqref{eq.intro_wave_2}--\eqref{eq.intro_wave_3} using a frequency-rescaled time,
\begin{equation}
\label{eq.intro_z} z \simeq t \langle \xi \rangle^\frac{1}{ \ell + 1 } \text{.}
\end{equation}
Indeed, taking Fourier transforms of \eqref{eq.intro_wave_2}--\eqref{eq.intro_wave_3}, one sees that the leading parts of the resulting differential equations can be expressed purely in terms of $z$ (and not of $\xi$).
This leads to uniform estimates over all frequencies, resulting in energy estimates similar to \eqref{eq.intro_wave_1scat}--\eqref{eq.intro_wave_1asymp}.

Another fundamental feature of the analysis in \cite{dreher_reissig_2000, dreher_witt_2002} is the microlocal partitioning of the time-frequency domain $( 0, T ]_t \times \R^d_\xi$ into the so-called \emph{pseudodifferential} and \emph{hyperbolic zones}:
\begin{equation}
\label{eq.intro_zones} Z_P := \{ z \lesssim 1 \} \text{,} \qquad Z_H := \{ z \gtrsim 1 \} \text{.}
\end{equation}
Roughly speaking, this allows one to isolate the singular behaviors at $z = 0$ and $z = \infty$, each of which provides a separate mechanism for derivative loss.
In particular, the solution $\phi$ has vastly different properties on $Z_P$ and $Z_H$, and one must analyze the regions separately.

In solving \eqref{eq.intro_wave_2}--\eqref{eq.intro_wave_3}, it is standard to rewrite the wave equation appropriately as a first-order system.
For example, for one spatial dimension, \eqref{eq.intro_wave_2} can be expressed on $Z_H$ as
\begin{equation}
\label{eq.intro_system_2} \partial_t \left[ \begin{matrix} \imath t^\ell \xi \, \hat{\phi} \\ \partial_t \hat{\phi} \end{matrix} \right] = \left\{ \imath t^\ell \xi \left[ \begin{matrix} 0 & 1 \\ a (t) & - 2 b (t) \end{matrix} \right] + t^{-1} \left[ \begin{matrix} \ell & 0 \\ - c (t) & 0 \end{matrix} \right] + \dots \right\} \left[ \begin{matrix} \imath t^\ell \xi \, \hat{\phi} \\ \partial_t \hat{\phi} \end{matrix} \right] \text{.}
\end{equation}
More precise first-order systems will be provided later on, but for now, one should notice that the first coefficient $\imath t^\ell \xi [ \cdot ]$ on the right-hand side of \eqref{eq.intro_system_2} captures the hyperbolicity, which degenerates as $t \searrow 0$.
Moreover, the lower-order coefficient $\ell t^{ \ell - 1 } c (t)$ of $\partial_x \phi$ in \eqref{eq.intro_wave_2} is captured in a \emph{Fuchsian} term $t^{-1} [ \cdot ]$ in \eqref{eq.intro_system_2}, which is ``critically singular" in that it scales in the same manner as ``$\partial_t$", and it barely fails to be $t$-integrable.
Finally, the remaining coefficients ``$\dots$" are remainders, in that they behave better than the Fuchsian term and are $t$-integrable.

In fact, the Fuchsian term $t^{-1} [ \dots ]$ provides the main mechanism for the loss of regularity as $t \searrow 0$.
To see this more closely, we rewrite \eqref{eq.intro_system_2} in terms of $z$ as
\begin{equation}
\label{eq.intro_system_2z} \partial_z \left[ \begin{matrix} \imath t^\ell \xi \, \hat{\phi} \\ \partial_t \hat{\phi} \end{matrix} \right] = \left\{ \imath z^\ell \left[ O (1) \right] + z^{-1} \left[ O (1) \right] + \dots \right\} \left[ \begin{matrix} \imath t^\ell \xi \, \hat{\phi} \\ \partial_t \hat{\phi} \end{matrix} \right] \text{.}
\end{equation}
This term $z^{-1} [ \dots ]$ is then handled by including integrating factors of the form $z^p$, which has the effect of adding extra powers of $\xi$---that is, extra derivatives---to the unknowns.

However, these integrating factors further complicate the analysis by introducing the same powers $z^p$ to the off-diagonal remainder coefficients ``$\dots$" in \eqref{eq.intro_system_2z}, which can have the effect of making the ``remainder" terms no longer integrable.
This obstacle was overcome in \cite{dreher_reissig_2000} by adopting a higher-order diagonalization scheme on $Z_H$.
(This was based on a scheme in \cite{reissig_yagdjian_2000} for a different problem of obtaining decay estimates as $t \nearrow \infty$.)
Here, one renormalizes \eqref{eq.intro_system_2z} into a form
\[
\partial_z \pd{U}{m} = \left\{ \sum_{ k = 1 }^m \pd{D}{k} + \mc{R} \right\} \pd{U}{m} \text{,}
\]
where $\pd{D}{1}, \dots, \pd{D}{m}$ yield an expansion (in negative powers of $z$) of diagonal matrices, and where the remainder $\mc{R}$ now contains a sufficiently negative power of $z$ that can absorb powers $z^p$ arising from the integrating factor while still remaining integrable.
From the above renormalized system, one can then obtain uniform estimates for the renormalized $\pd{U}{m}$.

The analysis on $Z_P$ is easier in comparison, and this was done in \cite{dreher_reissig_2000, dreher_witt_2002} in a more ad-hoc manner.
From all this, one can generate an overall \emph{renormalized unknown}---defined separately on $Z_P$ and $Z_H$---that is \emph{uniformly bounded} from $t = 0$ to positive times and vice versa, \emph{without any loss of derivatives}.
In particular, this allows one not only to produce \emph{asymptotics} as $t \searrow 0$ for solutions starting at some positive time, but also to solve the converse ``\emph{scattering}" problem from data at $t = 0$ for the renormalized unknown.
(Note that from this perspective, the term ``loss of regularity" is misleading, as the derivative loss exhibited in \eqref{eq.intro_wave_1scat}--\eqref{eq.intro_wave_1asymp} arises entirely from converting the renormalized quantity back to the original unknown $\hat{\phi}$.)

The aforementioned ideas extend to more general first-order systems having a similar form as \eqref{eq.intro_system_2} (or its higher-dimensional analogues)---roughly, containing a degenerating hyperbolicity and a singular Fuchsian coefficient as $t \searrow 0$; see, for instance, \cite{dreher_witt_2005}.

Finally, there are also numerous results in the literature that apply to weakly hyperbolic equations and systems with coefficients that depend on both $t$ and $x$.
Examples include the classical \cite{nersesian_1966, oleinik_1970}, as well as \cite{ebert_2005, ebert_kapp_filho_2009}, all of which employed physical space methods.
In addition, from the microlocal side, the aforementioned \cite{dreher_2002, dreher_witt_2005} treated, via pseudodifferential calculus, classes of first-order weakly hyperbolic systems with general $( t, x )$-dependent coefficients.

\begin{remark}
We mention that \cite{dreher_witt_2002, dreher_witt_2005} treated the systems \eqref{eq.intro_wave_2}--\eqref{eq.intro_wave_3} in a different manner using \emph{edge Sobolev spaces}.
These are \emph{spacetime} Sobolev spaces over $( 0, T ]_t \times \R^d_x$ containing additional pseudodifferential weights adapted to the scalings arising from the critical weak hyperbolicity.
\end{remark}

\begin{remark}
For the wave equations \eqref{eq.intro_wave_2}--\eqref{eq.intro_wave_3}, one can formulate the corresponding first-order system so that nontrivial Fuchsian term and higher-order renormalization only arises in $Z_H$.
(See \cite{dreher_witt_2005}, where this was also used as a crucial assumption for their general result on first-order systems.)
However, this will no longer be the case for the singular systems treated in this paper.
\end{remark}

\subsubsection{Critically Singular Systems}

Having identified the singular Fuchsian coefficients as the key mechanism for loss of derivatives, it then becomes natural to now widen our scope to wave equations with similar critically singular coefficients, and then to a larger class first-order hyperbolic systems with general Fuchsian coefficients.
This is one key contribution of the present article.

For singular wave equations, we can consider, first in one spatial dimension,
\begin{equation}
\label{eq.intro_wave_4} \partial_t^2 \phi - t^{ 2 \ell } a (t) \, \partial_x^2 \phi + 2 t^\ell b (t) \, \partial^2_{tx} \phi + t^{ \ell - 1 } c (t) \, \partial_x \phi + t^{-1} g (t) \, \partial_t \phi + t^{-2} h (t) \, \phi = 0 \text{,}
\end{equation}
where the functions $a, b, c, g, h$ are as before, and where the exponent satisfies $\ell > -1$.
(The crucial requirement for $\ell$ is simply that $t^\ell$ is integrable near $t = 0$.)
Then:
\begin{itemize}
\item On $Z_H$, one can write \eqref{eq.intro_wave_4} as the following first-order system:
\begin{equation}
\label{eq.intro_wave_4h} \partial_t \left[ \begin{matrix} \imath t^\ell \xi \hat{\phi} \\ \partial_t \hat{\phi} \end{matrix} \right] = \left\{ \imath t^\ell \xi \left[ \begin{matrix} 0 & 1 \\ a (t) & - 2 b(t) \end{matrix} \right] + t^{-1} \left[ \begin{matrix} \ell & 0 \\ - c (t) & - g (t) \end{matrix} \right] + \dots \right\} \left[ \begin{matrix} \imath t^\ell \xi \hat{\phi} \\ \partial_t \hat{\phi} \end{matrix} \right] \text{.}
\end{equation}

\item On $Z_P$, one can express \eqref{eq.intro_wave_4} as follows:
\begin{equation}
\label{eq.intro_wave_4p} \partial_t \left[ \begin{matrix} t^{-1} \hat{\phi} \\ \partial_t \hat{\phi} \end{matrix} \right] = \left\{ t^{-1} \left[ \begin{matrix} -1 & 1 \\ - h (t) & - g (t) \end{matrix} \right] + \dots \right\} \left[ \begin{matrix} t^{-1} \hat{\phi} \\ \partial_t \hat{\phi} \end{matrix} \right] \text{.}
\end{equation}
\end{itemize}
As before, ``$\dots$" denotes error terms that do not alter qualitative behaviors of solutions.
(On $Z_H$, these contain additional negative powers of $z$, while on $Z_P$, these contain positive powers of $z$.)

\begin{remark}
Note the weakly hyperbolic term $\imath t^\ell \xi [ \dots ]$ is only relevant on $Z_H$, since it serves merely as a remainder term on $Z_P$.
Therefore, this term is not explicitly listed in \eqref{eq.intro_wave_4p}.
\end{remark}

Observe that \eqref{eq.intro_wave_4h}--\eqref{eq.intro_wave_4p} exhibit \emph{nontrivial Fuchsian contributions on both $Z_P$ and $Z_H$}.
The upshot of this is that one may need to perform separate \emph{higher-order renormalizations} (to arbitrarily high order) \emph{on both $Z_P$ and $Z_H$}.
In particular, on $Z_P$, the need for such a renormalization implies that only a carefully constructed higher-order expansion involving $\hat{\phi}$ and its derivatives may converge to an asymptotic limit as $t \searrow 0$.
As we shall see, in some cases, the asymptotic expansion may be \emph{polyhomogeneous}, i.e.\ it may also contain powers of $\log z$.
Thus, a key challenge for this paper is to generate this asymptotic quantity, even when elaborate logarithmic corrections are required.

Analogous first-order formulations hold for singular wave equations in higher-dimensions,
\begin{align}
\label{eq.intro_wave_5} \partial_t^2 \phi - t^{ 2 \ell } \sum_{ i, j = 1 }^d a_{ij} (t) \, \partial_{ x_i x_j }^2 \phi + 2 t^\ell \sum_{ i = 1 }^d b_i (t) \, \partial^2_{ t x_i } \phi \qquad & \\
\notag + t^{ \ell - 1 } \sum_{ i = 1 }^d c_i (t) \, \partial_{ x_i } \phi + t^{-1} g (t) \, \partial_t \phi + t^{-2} h (t) \, \phi &= 0 \text{,}
\end{align}
though we defer the precise expressions until Section \ref{sec.wave}.
These higher-dimensional systems can be treated in a similar manner as the one-dimensional case above.

For a more physically motivated example, as well as one that lies beyond the scopes of \eqref{eq.intro_wave_4} and \eqref{eq.intro_wave_5}, we consider scalar wave equations on general Kasner spacetimes:
\begin{equation}
\label{eq.intro_wave_kasner} \partial_t^2 \phi - \sum_{ i = 1 }^d t^{ 2 \ell_i } \partial_{ x_i }^2 \phi + t^{-1} \partial_t \phi = 0 \text{,} \qquad \ell_1, \dots, \ell_d > -1 \text{.}
\end{equation}
Observe that \eqref{eq.intro_wave_kasner} has a similar weakly hyperbolic and singular structure as \eqref{eq.intro_wave_5}.
However, the main novelty in \eqref{eq.intro_wave_kasner}, relative to \eqref{eq.intro_wave_5}, is the \emph{anisotropy}---the hyperbolicity degenerates as $t \searrow 0$ at different rates (i.e.\ different powers of $t$) in the various spatial directions.

The motivation for \eqref{eq.intro_wave_kasner} is that it is the wave operator associated to the \emph{Kasner metrics}
\begin{equation}
\label{eq.intro_kasner_metric} g := - dt^2 + \sum_{ i = 1 }^d t^{ -2 \ell_i } \, d x_i^2 \text{.}
\end{equation}
(The parameters $-\ell_1, \dots, -\ell_d$ are commonly known as the \emph{Kasner exponents}.)
In relativity and cosmology, these geometries serve as common models of \emph{big bang singularities}.
From this point of view, \eqref{eq.intro_wave_kasner} is useful as a linearized singularity model that is far easier to study than the full Einstein equations.
Kasner dynamics also serve as models for more general spacelike singularities---for instance, Schwarzschild black hole interiors \cite{fournodavlos_sbierski_2020} or asymptotically de Sitter infinity \cite{cicortas_2023}.

It is well-known within mathematical relativity that solutions of \eqref{eq.intro_wave_kasner} have modified asymptotics with a logarithmic correction and lose derivatives as $t \searrow 0$; see, for instance, \cite{alho_fournodavlos_franzen_2019, petersen_2016, ringstrom_2020}.
However, these results (barely) fail to obtain the optimal derivative loss.
Of particular interest is the more recent \cite{li_2024}, which achieved the precise derivative loss by modifying the asymptotic quantity by an extra logarithmic derivative (in our current context, arising from the discrepancy between $t$ and $z$).

Roughly speaking, \cite{li_2024} showed the following quantities had finite limits as $t \searrow 0$,
\begin{equation}
\label{eq.intro_kasner_limit} \lim_{ t \searrow 0 } ( \hat{\phi} - t \log z \, \partial_t \hat{\phi} ) ( t, \xi ) = \hat{\varphi}_{ +, 0 } ( \xi ) \text{,} \qquad \lim_{ t \searrow 0 } t \partial_t \hat{\phi} ( t, \xi ) = \hat{\varphi}_{ -, 0 } ( \xi ) \text{,}
\end{equation}
with $z$ an appropriately defined rescaling of $t$; see \eqref{eq.intro_zz}.
Moreover, \cite{li_2024} showed, for $t > 0$, that
\begin{equation}
\label{eq.intro_kasner_est} \| \phi (t) \|_{ H^s } + \| \partial_t \phi (t) \|_{ H^{ s - 1 } } \simeq \| \varphi_{ +, 0 } \|_{ H^{ s - \frac{1}{2} } } + \| \varphi_{ -, 0 } \|_{ H^{ s - \frac{1}{2} } }
\end{equation}
(with the constants depending on $t$), demonstrating a precise loss of $\frac{1}{2}$-derivatives at $t = 0$.
Also, \cite{li_2024} obtained scattering solutions starting from $t = 0$, with asymptotic data imposed for $\varphi_{ \pm, 0 }$.

A key element of the analysis in \cite{li_2024} is the (independent) rediscovery of some ideas mentioned above, in particular the partition into and separate treatment of the zones $Z_P$ and $Z_H$.
Moreover, a particularly important novelty in \cite{li_2024} is that it demonstrated the aforementioned ideas were, with some minor modifications, sufficient for treating the anisotropic degeneracies in \eqref{eq.intro_wave_kasner}.

Finally, we mention the connection between the present problem and that of decay and asymptotics for wave equations and weakly hyperbolic systems as $t \nearrow \infty$.
Indeed, rewriting the equation in terms of $\tau := t^{-1}$ converts the above to an asymptotic problem at $\tau \searrow 0$.
Since $\partial_t$ and $t^{-1}$ have the same scaling, this transformation preserves the system's Fuchsian structure.
In fact, many results involving decay at $t \nearrow \infty$ employ a similar microlocal decomposition into zones as described above (though with different natural scalings); see, e.g., \cite{reissig_yagdjian_2000, ruzhansky_wirth_2015, wirth_2017}.

Also, returning to relativity, we mention the articles \cite{anderson_2005, cicortas_2023}, which established precise asymptotics and scattering results at infinity for wave equations on de Sitter spacetime.
In contrast to \eqref{eq.intro_wave_kasner}, the degeneracies here are not anisotropic, so the results lie in the framework of \eqref{eq.intro_wave_5}.
As before, solutions exhibit modified asymptotics and shifted regularity at infinity.
Furthermore, the results were extended to a class of asymptotically de Sitter vacuum spacetimes in \cite{cicortas_2024b}.

\begin{remark}
We note the setting of \cite{li_2024} was on $( 0, T ] \times \mathbb{T}^d_x$, in which each level set of $t$ is a torus.
Similarly, in \cite{cicortas_2023}, each level set of $t$ is a sphere $\Sph^d$.
Nonetheless, the difference in setting only yields minor, cosmetic differences in the analysis---on $\mathbb{T}^d$, the continuous Fourier transform is replaced by discrete Fourier series, while on $\Sph^d$, one uses spectral decompositions instead.
\end{remark}

\subsection{The Main Results}

The main results of this article obtain sharp asymptotics, scattering, and loss of regularity for a very general class of weakly hyperbolic systems with critically singular (Fuchsian) coefficients.
We note here some key features of our result:
\begin{itemize}
\item This general class of systems includes the wave equations \eqref{eq.intro_wave_4} and \eqref{eq.intro_wave_5}, the wave equation \eqref{eq.intro_wave_kasner} on Kasner (allowing also for a larger class of lower-order coefficients), and higher-order weakly hyperbolic and singular differential equations.

\item We are able to treat general Fuchsian terms, and we provide a precise accounting of how these terms affect the asymptotics and the loss of regularity as $t \searrow 0$.

\item We also allow for general perturbative terms, i.e.\ those that ``behave better than the Fuchsian term".
While these terms may affect the asymptotics and scattering theories beyond the leading term, we show that they do not affect the loss of regularity as $t \searrow 0$.
\end{itemize}
In particular, our result combines ideas from both the weakly hyperbolic (e.g.\ \cite{dreher_reissig_2000, dreher_witt_2002, dreher_witt_2005, nunes_wirth_2015, reissig_yagdjian_2000, wirth_2017}) and mathematical relativity (e.g.\ \cite{li_2024, ringstrom_2020}) literature, along with some novel developments.

On the other hand, the key restriction of our result is that the \emph{coefficients of our system depend only on $t$}, and not on the spatial variable $x$.
The general setting of $(t, x)$-dependent coefficients will be treated in forthcoming works.
The main reason for considering only $t$-dependent coefficients is that one can perform the analysis directly in Fourier space.
This lets us more clearly highlight the mechanisms behind our renormalizations and regularity loss, without an additional layer of technicality, e.g.\ pseudodifferential operators.
Moreover, even this restricted class already encompasses some settings of physical interest, e.g.\ linearized systems in Kasner backgrounds.

\subsubsection{Setting and Assumptions}

Unfortunately, some technical setup is required in order to state our main result.
We give an informal summary here; see Section \ref{sec.system_prelim} for the precise development.

At the most basic level, we are considering on $( 0, T ]_t \times \R^d_x$ a first-order system of the form
\begin{equation}
\label{eq.intro_system} \partial_t \check{U} = \mc{A} ( t, \nabla_x ) \check{U} \text{,} \qquad \check{U}: ( 0, T ]_t \times \R^d_x \rightarrow \C^n \text{.}
\end{equation}
Since the operator $\mc{A}$ has no $x$-dependence, we can consider the Fourier transform of \eqref{eq.intro_system},
\begin{equation}
\label{eq.intro_systemf} \partial_t U = \mc{A} ( t, \xi ) U \text{,} \qquad U: ( 0, T ]_t \times \R^d_\xi \rightarrow \C^n \text{,}
\end{equation}
that is, we can view \eqref{eq.intro_systemf} as a system of ordinary differential equations for every fixed frequency $\xi \in \R^d$.
Thus, for our main result, we must impose appropriate assumptions on the coefficients $\mc{A}$.

The first step is to quantify the degenerating hyperbolicity of our system.
For this, we set
\begin{equation}
\label{eq.intro_hyp} \mf{H}^2 ( t, \xi ) := \sum_{ i, j = 1 }^d t^{ \ell_i + \ell_j } \lambda_{ij} (t) \, \xi_i \xi_j \text{,}
\end{equation}
where $\ell_1, \dots, \ell_d > -1$, and where $\lambda_{ij} (t)$ are coefficients such that $\mf{H}^2$ satisfies an ellipticity conditions; see \eqref{eq.system_ellip}.
Note that by choosing different values for the $\ell_i$'s, we capture settings where the hyperbolicity degenerates in an anisotropic manner, such as in Kasner backgrounds.

With the hyperbolicity in place, we can now define our analogue of ``$z$" from before,
\begin{equation}
\label{eq.intro_zz} z := \langle \langle \xi \rangle \rangle t \text{,} \qquad \langle \langle \xi \rangle \rangle := \max_{ 1 \leq i \leq d } \langle \xi_i \rangle^\frac{1}{ \ell_i + 1 } \text{,}
\end{equation}
representing the ($\xi$-dependent) rescaled time that is adapted to the hyperbolicity $\mf{H}$.
Using this $z$, we then partition our region $( 0, T ]_t \times \R^d_\xi$ into \emph{three} microlocal zones:
\begin{itemize}
\item \emph{Pseudodifferential zone}: $Z_P := \{ z \ll 1 \}$.

\item \emph{Intermediate zone}: $Z_I := \{ z \simeq 1 \}$.

\item \emph{Hyperbolic zone}: $Z_H := \{ z \gg 1 \}$.
\end{itemize}
The reason for now requiring three zones is that in general, the renormalizations on $Z_P$ and $Z_H$ will only be well-defined for sufficiently small and large $z$, respectively.
However, as the new intermediate zone $Z_I$ lies away from all the singular behavior, the analysis there will in practice be trivial.

We can now state the assumptions that we impose on our system.
The key complication, which makes the main result difficult to state, is we must impose different assumptions in each microlocal zone.
This is because the aspects our system that are considered ``leading-order", as well as which powers of $z$ are considered ``good", depend on the zone in question.

First, on $Z_I$ (see Assumption \ref{ass.system_i}), we impose only a minimal condition:
\begin{equation}
\label{eq.intro_ass_i} \mc{A} |_{ Z_I } = \langle\langle \xi \rangle\rangle \, O (1) \text{.}
\end{equation}
In practice, \eqref{eq.intro_ass_i} is trivially satisfied, since $z$ is bounded from above and below on $Z_I$, so that the presence of any derivatives will not cause trouble for our analysis.
In contrast, for the remaining zones, the main assumption will entail diagonalizing the leading-order part of $\mc{A}$.

On $Z_P$ (see Assumption \ref{ass.system_p}), we assume there exists a matrix-valued function $M_P$ such that
\begin{equation}
\label{eq.intro_ass_p} M_P ( \mc{A} |_{ Z_P } ) M_P^{-1} + ( \partial_t M_P ) M_P^{-1} = t^{-1} B_P + O ( t^{-1} z^\varepsilon ) \text{,}
\end{equation}
for some $\varepsilon > 0$.
Here, $M_P$ serves to diagonalize (as much as possible) the leading part of $\mc{A}$, which on $Z_P$ is the Fuchsian term.
Thus, we will assume in \eqref{eq.intro_ass_p} that \emph{$B_P$ is everywhere in Jordan normal form}.
The remaining term ``$O ( t^{-1} z^\varepsilon )$" behaves strictly better than the Fuchsian term $t^{-1} B_P$ (note positive powers of $z$ are small on $Z_P$) and can be viewed as ``remainder".

Finally, on $Z_H$ (see Assumption \ref{ass.system_h}), we assume there is a matrix-valued $M_H$ such that
\begin{equation}
\label{eq.intro_ass_h} M_H ( \mc{A} |_{ Z_H } ) M_H^{-1} + ( \partial_t M_H ) M_H^{-1} = \imath \mf{H} D_H + t^{-1} B_H + O ( t^{-1} z^{ -\varepsilon } ) \text{,}
\end{equation}
for some $\varepsilon > 0$.
Again, $M_H$ serves as a leading-order diagonalizer of $\mc{A}$, though the leading-order part is now the hyperbolic part.
Thus, we assume \emph{$D_H$ in \eqref{eq.intro_ass_h} is everywhere diagonal and real-valued}, i.e.\ the system is strictly hyperbolic at positive times.
On top of the above, we also assume \emph{the eigenvalues of $D_H$ are uniformly separated from each other}, that is, the system has uniformly separated characteristic speeds, modulo the degenerating hyperbolicity given by $\mf{H}$.
Next, ``$t^{-1} B_H$" is the Fuchsian term in this setting, though we can no longer assume $B_H$ has any diagonal structure.
Lastly, ``$O ( t^{-1} z^{ -\varepsilon } )$" captures the remainder terms that behave strictly better than the Fuchsian part (note, in contrast to before, that negative powers of $z$ are well-behaved on $Z_H$).

\begin{remark}
Observe that on $Z_P$, the hyperbolic part of $\mc{A}$ is contained within the ``remainder" terms, since it contains additional positive powers of $z$ compared to the Fuchsian term.
\end{remark}

\subsubsection{The Theorem Statements}

We can now give informal statements of our main results:

\begin{theorem}[Asymptotics] \label{thm.intro_asymp}
Assuming the above development and data $u_T: \R^d_\xi \rightarrow \C^n$, there exists a unique solution $U$ of \eqref{eq.intro_systemf} with $U (T) = u_T$.
Furthermore, there exists a (systematically computed) renormalization $U_A$ of $U$ such that $U_A (t)$ has finite limit as $t \searrow 0$:
\begin{equation}
\label{eq.intro_asymp_limit} \lim_{ \tau \searrow 0 } U_A ( \tau, \xi ) = u_A ( \xi ) \in \C^n \text{,} \qquad \xi \in \R^d \text{.}
\end{equation}
In addition, $U_A$ satisfies the following uniform bounds:
\begin{equation}
\label{eq.intro_asymp} | u_A ( \xi ) | \lesssim | U_A ( T, \xi ) | \text{,} \qquad \xi \in \R^d \text{.}
\end{equation}
\end{theorem}

\begin{theorem}[Scattering] \label{thm.intro_scattering}
Assuming the above development and asymptotic data $u_A: \R^d_\xi \rightarrow \C^n$, there exists a unique solution $U$ of \eqref{eq.intro_systemf} that achieves this asymptotic data as $t \searrow 0$,
\begin{equation}
\label{eq.intro_scattering_limit} \lim_{ \tau \searrow 0 } U_A ( \tau, \xi ) = u_A ( \xi ) \text{,} \qquad \xi \in \R^d \text{,}
\end{equation}
with $U_A$ being the same renormalization as in Theorem \ref{thm.intro_asymp}.
Also, $U_A$ satisfies the bounds
\begin{equation}
\label{eq.intro_scattering} | U_A ( T, \xi ) | \lesssim | u_A ( \xi ) | \text{.}
\end{equation}
\end{theorem}

For the full, precise versions of Theorems \ref{thm.intro_asymp} and \ref{thm.intro_scattering}, the reader is referred to Theorem \ref{thm.energy_main}.
In the meantime, a number of remarks regarding Theorems \ref{thm.intro_asymp} and \ref{thm.intro_scattering} are in order.

\begin{remark}
Although Theorems \ref{thm.intro_asymp} and \ref{thm.intro_scattering} were stated for homogeneous systems \eqref{eq.intro_system}--\eqref{eq.intro_systemf}, both results extend directly to \emph{inhomogeneous} systems with a \emph{forcing term} $\mc{F}$, i.e.
\[
\partial_t U = \mc{A} ( t, \xi ) U + \mc{F} ( t, \xi ) \text{,}
\]
provided $\mc{F}$ (or rather, its corresponding renormalization) is sufficiently time-integrable.
For details, see the precise result, Theorem \ref{thm.energy_main}, which includes such forcing terms in its treatment.
\end{remark}

\begin{remark}
The renormalized quantity $U_A$ in Theorems \ref{thm.intro_asymp} and \ref{thm.intro_scattering} is constructed by applying separate renormalizations on $Z_P$ and $Z_H$; see \eqref{eq.system_UA}, along with Propositions \ref{thm.energy_p} and \ref{thm.energy_h}, for the precise formulas.
In addition, a closer inspection of the proofs of Theorems \ref{thm.intro_asymp}--\ref{thm.intro_scattering} yields \emph{algorithms for systematically computing both renormalizations to any arbitrary finite order}.
\end{remark}

\begin{remark}
Note \eqref{eq.intro_asymp} and \eqref{eq.intro_scattering} imply $U_A$ satisfies uniform estimates \emph{without any loss of derivatives}.
In other words, Theorems \ref{thm.intro_asymp} and \ref{thm.intro_scattering} are converses of each other, and \emph{the asymptotics and scattering theories are fully reversible}.
From this perspective, the aforementioned loss of regularity arises when one transforms from the renormalized $U_A$ to the original unknown $U$.
\end{remark}

\begin{remark}
From the proof of Theorems \ref{thm.intro_asymp}--\ref{thm.intro_scattering}, we can pinpoint the precise source of the loss of regularity and modified asymptotics as $t \searrow 0$.
In particular, these are consequences of:
\begin{itemize}
\item The \emph{zero-time limit} of $B_P$, that is, $B_{ P, 0 } := \lim_{ t \searrow 0 } B_P$.

\item The \emph{infinite-frequency limit} of $B_H$, that is, $B_{ H, \infty } := \lim_{ | \xi | \nearrow \infty } B_H$.
\end{itemize}
See the discussions near Assumptions \ref{ass.system_p} and \ref{ass.system_h}, as well as Propositions \ref{thm.energy_p} and \ref{thm.energy_h}.
\end{remark}

\begin{remark}
Note the asymptotic limits stated in \eqref{eq.intro_asymp_limit} and \eqref{eq.intro_scattering_limit} are pointwise in frequency.
However, from the point of view of \eqref{eq.intro_system}, it is far more natural to consider convergences in Sobolev norms.
We show that, with respect to $H^s$-norms, there is a loss of $\varepsilon$ (for any $\varepsilon > 0$) derivatives in the asymptotic limit; see Proposition \ref{thm.energy_conv} for a precise statement.
This discrepancy arises from the fact that the estimates one proves for \eqref{eq.intro_systemf} are uniform in the rescaled time $z$ rather than in $t$.
\end{remark}

\begin{remark}
Recall Theorems \ref{thm.intro_asymp}--\ref{thm.intro_scattering} assumed that the eigenvalues of $D_H$ are uniformly separated.
However, \emph{one can still obtain weaker results when this assumption on $D_H$ is removed}; see Theorem \ref{thm.energy_ex} for the precise statement.
In short, one has analogues of Theorems \ref{thm.intro_asymp}--\ref{thm.intro_scattering}, but the asymptotics and scattering theories may no longer be converses of each other, as one may now have \emph{different renormalized unknowns $U_A$ in the asymptotics and scattering theories}.
These more general results will be crucial when we study the linearized Einstein-scalar system in Section \ref{sec.einstein}.
\end{remark}

\subsubsection{Ideas of the Proof}

The proofs of Theorems \ref{thm.intro_asymp}--\ref{thm.intro_scattering} (more accurately, the precise Theorem \ref{thm.energy_main}) can be found throughout Section \ref{sec.system}.
Below, we informally discuss a few key ideas.

As mentioned before, the key step is to construct, via higher-order diagonalization processes, separate renormalized quantities on $Z_P$ and $Z_H$ that remain controlled down to $t \searrow 0$.
Roughly speaking, on $Z_P$, the strategy is to obtain a linear transformation $\pd{\mc{Q}_P}{m}$, for any $m \geq 0$, such that
\begin{equation}
\label{eq.intro_proof_p1} \partial_z ( \pd{\mc{Q}_P}{m} U ) = \left\{ z^{-1} B_P + \sum_{ k = 1 }^m \pd{D_P}{k} + \pd{\mc{R}_P}{m} \right\} ( \pd{\mc{Q}_P}{m} U ) \text{,}
\end{equation}
where the $\pd{D_P}{k}$'s are diagonal, and where $\pd{\mc{R}_P}{m}$ is higher-order remainder satisfying
\[
\pd{D_P}{k} = O ( z^{ -1 + k \varepsilon } ) \text{,} \qquad \pd{\mc{R}_P}{m} = O ( z^{ -1 + m \varepsilon } ) \text{,} \qquad 1 \leq k < m \text{,} \quad \varepsilon > 0 \text{.}
\]
Then, as long as $m$ is sufficiently large, the remainder $\pd{\mc{R}_P}{m}$ will be nice enough to absorb any amplifying powers of $z$ arising from the integrating factors due to the Fuchsian term $z^{-1} B_P$, and one will be able apply standard estimates to the diagonalized system \eqref{eq.intro_proof_p1}.

The diagonalization here is similar to the scheme utilized in \cite{wirth_2017} (in the context of wave decay as $t \nearrow \infty$; see also \cite{wirth_2017a}); details of the procedure were provided in \cite{nunes_wirth_2015}.
One applies an induction on $m$; given $\pd{\mc{Q}_P}{m}$, one then obtains the next-order transformation, $\pd{\mc{Q}_P}{m+1}$, by solving for the difference $\pd{\mc{Q}_P}{m+1} - \pd{\mc{Q}_P}{m} = O ( z^{ ( m + 1 ) \varepsilon } )$.
In particular, when this difference solves an appropriate system of differential equations (see \eqref{eql.perf_diag_p_2}), the resulting remainder $\pd{\mc{R}_P}{m+1}$ gains an extra power $z^\varepsilon$ over the previous $\pd{\mc{R}_P}{m}$; see Lemma \ref{thm.perf_diag_p} for details.
With the above in hand, the final renormalized quantity on $Z_P$ is obtained by applying an integrating factor to do away with the Fuchsian term:
\[
U_A |_{ Z_P } = \exp ( - z^{-1} B_P ) \pd{\mc{Q}_P}{m} U \text{.}
\]
Observe in particular that this renormalized quantity takes the form
\begin{equation}
\label{eq.intro_proof_p2} U_A |_{ Z_P } \sim \operatorname{diag} ( z^{ p_1 }, \dots, z^{ p_n } ) \bigg[ I_n + \sum_{ k = 1 }^m O ( z^{ k \varepsilon } ) \bigg] \, U \text{.}
\end{equation}

Thus far, all the preceding development assumed \emph{$B_P$ is diagonal}.
One novelty of our result is that we can also extend the above scheme to treat $B_P$'s that are merely in Jordan normal form.
In this case, \emph{the expansion for $U_A |_{ Z_P }$ will be polyhomogeneous}, in that it will also contain powers of $\log z$.
Furthermore, there exist special values of $B_P$ (both diagonal or otherwise) for which one obtains powers of $\log z$ in higher-order terms of the expansion.
While these cases were excluded in \cite{wirth_2017}, our renormalization scheme also treats these in a unified and systematic manner.

The objective on $Z_H$ is similar---to obtain a linear transformation $\pd{\mc{Q}_H}{m}$, for $m \geq 1$, with
\begin{equation}
\label{eq.intro_proof_h1} \partial_z ( \pd{\mc{Q}_H}{m} U ) = \left\{ \imath \langle\langle \xi \rangle\rangle^{-1} \mf{H} D_H + \sum_{ k = 1 }^m \pd{D_H}{k} + \pd{\mc{R}_H}{m} \right\} ( \pd{\mc{Q}_H}{m} U ) \text{,}
\end{equation}
again with the $\pd{D_H}{k}$'s diagonal, and with $\pd{\mc{R}_P}{m}$ being a higher-order remainder, satisfying
\[
\pd{D_H}{k} = O ( z^{ -1 - (k - 1) \varepsilon } ) \text{,} \qquad \pd{\mc{R}_P}{m} = O ( z^{ -1 - (m - 1) \varepsilon } ) \text{,} \qquad 1 \leq k < m \text{,} \quad \varepsilon > 0 \text{.}
\]
(In particular, observe that the Fuchsian part of the diagonal expansion is encoded in $\pd{D_H}{1}$.)
However, the key difference with $Z_P$ is that the leading coefficient in $Z_H$ is the hyperbolic part $\imath \mf{H} D_H$, so the renormalization process will have to proceed differently.

Our renormalization procedure on $Z_H$ is directly inspired by that used in the decay result of \cite{reissig_yagdjian_2000}.
The rough idea is that, by using the hyperbolic term in \eqref{eq.intro_proof_h1} as integrating factor, then the uniform separation of eigenvalues of $D_H$ yields nontrivial phases in the non-diagonal components of the system, from which one can extract negative (good) powers of $z$ via stationary phase.

In practice, the diagonalization process is again an induction on $m$.
The key technical difference is that the equation for $\pd{\mc{Q}_H}{m+1} - \pd{\mc{Q}_H}{m} = O ( z^{ -m \varepsilon } )$ that yields an improvement in the subsequent $\pd{\mc{R}_H}{m+1}$ over $\pd{\mc{R}_H}{m}$ is now algebraic, rather than differential, in nature; see \eqref{eql.perf_diag_h_1}.
(In particular, this equation involves dividing by the differences of wave speeds.)
The full details of the diagonalization can be found in Lemma \ref{thm.perf_diag_h}, though we note here that the Fuchsian part of the diagonal expansion is given precisely by the diagonal elements of $B_H$.
Therefore, the contribution of $Z_H$ to the modified asymptotics and the loss of regularity at $t \searrow 0$ is driven by these diagonal elements of $B_H$.

The final renormalized quantity on $Z_H$ is then roughly of the form
\[
U_A |_{ Z_H } = \exp \bigg( \mathord{-} \int \pd{D_H}{1} dz \bigg) \pd{\mc{Q}_H}{m} U \text{.}
\]
Less precisely, the renormalized quantity will have the schematic form
\begin{equation}
\label{eq.intro_proof_h2} U_A |_{ Z_H } \sim \mc{E}_H \bigg[ I_n + \sum_{ k = 1 }^m O ( z^{ -k \varepsilon } ) \bigg] \, U \text{.}
\end{equation}
(In practice, the integrating factor $\mc{E}_H$ is a diagonal matrix usually involving powers of $z$ or $\mf{H}$.)

Finally, a novel aspect of our result is that we also systematically obtain asymptotics and scattering \emph{without assuming uniform separation of the eigenvalues of $D_H$}; see Theorem \ref{thm.energy_ex} for precise statements.
This is attained by adapting some ideas briefly mentioned in \cite{dreher_witt_2005}.
Here, the key difficulty is that we can only partially diagonalize up to higher orders; more specifically, we can only do away with non-diagonal components for which the two corresponding wave speeds are uniformly separated.
In other words, \emph{the $\pd{D_H}{k}$'s in \eqref{eq.intro_proof_h1} are now only block diagonal matrices}.

To get around this, we can bound each block of the Fuchsian block diagonal matrix $\pd{D_H}{1}$ from above or from below (in the sense of quadratic forms) by a diagonal matrix, depending on whether one is solving forward or backwards in time.
This yields diagonal Fuchsian parts in our renormalized system, but the price to be paid is that one in general obtains \emph{different diagonal Fuchsian terms}, and hence \emph{different renormalized quantities}, for the asymptotics and scattering settings.

In spite of the above, we are still able to solve our system \eqref{eq.intro_systemf} uniquely, both backwards to $t = 0$ (asymptotics), and forwards from $t = 0$ (scattering).
However, the above-mentioned discrepancy in the asymptotics and scattering renormalizations means that \emph{we cannot generally prove a (normed space) isomorphism between appropriate data at $t = T$ and $t = 0$}.

\subsection{Beyond the Main Results}

After Section \ref{sec.system}, the remainder of the paper is concerned with how our general results can be applied to a multitude of special cases.

\subsubsection{Applications}

In Section \ref{sec.wave_basic}, we study \emph{weakly hyperbolic wave equations with singular coefficients}, \eqref{eq.intro_wave_5}, and we demonstrate in detail how such equations fit into the framework of Section \ref{sec.system}.
In particular, \eqref{eq.intro_wave_5} satisfies the requisite assumptions if one takes, as unknowns,
\begin{equation}
\label{eq.intro_wave_U} U |_{ Z_P } := \left[ \begin{matrix} t^{-1} \hat{\phi} \\ \partial_t \hat{\phi} \end{matrix} \right] \text{,} \qquad U |_{ Z_H } := \left[ \begin{matrix} \imath t^\ell | \xi | \hat{\phi} \\ \partial_t \hat{\phi} \end{matrix} \right] \text{.}
\end{equation}
The precise estimates capturing the loss of regularity are found in Theorem \ref{thm.wave_main}.

Furthermore, by setting $g (0) = h (0) = 0$, we can consider the special case of weakly hyperbolic \emph{non-singular} wave equations covered in the literature, e.g.\ \cite{dreher_reissig_2000, dreher_witt_2002, dreher_witt_2005}.
The main estimates in this case are stated in Corollary \ref{thm.wave_nonsing} and reduce to the estimates \eqref{eq.intro_wave_1scat}--\eqref{eq.intro_wave_1asymp} for the model equation \eqref{eq.intro_wave_1}.

\begin{remark}
In some cases, Corollary \ref{thm.wave_nonsing} yields a slight improvement over the estimates of \cite{dreher_reissig_2000, dreher_witt_2002, dreher_witt_2005}.
More specifically, while the regularity loss was previously expressed as a supremum over $\smash{\R^d_\xi}$ of (parts of) $B_H$, here we improve this to a supremum over all directions $\omega \in \Sph^{ d - 1 }$ of the infinite-frequency limit of $B_H$.
(An analogous infinite-frequency limit was also used in the decay results of \cite{wirth_2017}.)
\end{remark}

Next, in Section \ref{sec.wave_kasner}, we turn our attention to \emph{wave equations with anisotropic degeneracies}, once again showing how these fit into the framework of Section \ref{sec.system}.
This class includes, as special cases, wave equations on Kasner spacetimes.
For Kasner settings, we demonstrate that our main result recovers the optimal asymptotics and derivative loss obtained in \cite{li_2024}.

In Section \ref{sec.wave_higher}, we outline how the analysis of Section \ref{sec.wave_basic} can be extended from wave equations to \emph{higher-order critically weakly hyperbolic equations with singular coefficients}.
(Higher-order weakly hyperbolic equations \emph{without} singular coefficients were treated in \cite{dreher_witt_2005}.)
While general formulas for the modified asymptotics and loss of regularity become too complicated to state explicitly, here we at least demonstrate how these can be systematically obtained.

\begin{remark}
For weakly hyperbolic wave and higher-order equations from Sections \ref{sec.wave_basic}--\ref{sec.wave_higher} with regular, non-singular coefficients, $\smash{\hat{\phi}}$ itself has a finite asymptotic limit at $t \searrow 0$.
Observe that in this setting, our results immediately imply the standard $C^\infty$-well-posedness for these equations from $t = 0$.
\end{remark}

Finally, in Section \ref{sec.einstein}, we revisit the linearized Einstein-scalar system about Kasner backgrounds studied in \cite{li_2024}, and we apply the methods of Section \ref{sec.system} to recover the energy estimates---and hence the main asymptotics and scattering results---obtained in \cite{li_2024} for this system.
Further, while \cite{li_2024} only considered Kasner exponents satisfying an additional subcriticality condition (see \eqref{eq.kasner_subcritical}), here we \emph{extend the energy estimates to all non-degenerate Kasner exponents}.

The novel ingredient leading to the above-mentioned improvement comes from a refined analysis on $Z_P$.
In particular, when the subcriticality condition holds, one can obtain the requisite estimates on $Z_P$ without any higher-order renormalization.
\emph{In the absence of subcriticality}, however, one must \emph{diagonalize the system to higher-order (i.e.\ renormalize the unknowns to higher order)} in order to construct quantities that have a finite asymptotic limit as $t \searrow 0$.

On $Z_H$, the desired estimates can be already obtained using the high-frequency energy estimates of \cite{li_2024}, however here we revisit the analysis in the framework of Section \ref{sec.system}.
From this perspective, one difficulty is that many components of the system propagate at common speeds, so we are in the setting of the weaker Theorem \ref{thm.energy_ex}.
Thus, one must find additional structure in the Fuchsian part $B_H$ in order to obtain time-reversible (i.e.\ converse) asymptotics and scattering theories.

\subsubsection{Outlook}

As mentioned before, the key restriction of our results here is that they only apply to systems with coefficients depending only on $t$.
The natural next step---which we address in upcoming work---is to handle more general coefficients depending on both $t$ and $x$.
In this setting, one expects that the loss of regularity will also depend on $x$ and can vary between points in space.
In fact, \cite{dreher_2002, dreher_witt_2005} treated coefficients depending on both $t$ and $x$ for smaller classes of weakly hyperbolic systems (without singular coefficients and anisotropic degeneracies); while they quantified loss of regularity in this setting, the systems they studied do not exhibit modified asymptotics at $t = 0$.

While the analysis in Section \ref{sec.einstein} highlights the applicability of our results to the Einstein equations, it would be of interest to extend this to a full treatment of the linearized Einstein-scalar system about Kasner backgrounds, i.e.\ a full extension of \cite{li_2024} to all non-degenerate Kasner exponents.
This would require treating additional issues pertaining to higher-order renormalizations but outside the immediate scope of Section \ref{sec.system}, e.g.\ asymptotic gauge choices and behaviors of the constraints.

Next, this article serves as a first step in a larger effort to study \emph{nonlinear} degenerate hyperbolic systems with singular coefficients.
Of particular interest are \emph{quasilinear} systems, where the principal part itself depends on the unknown.
One prerequisite for analyzing quasilinear equations is a robust linear theory that can handle large classes of principal parts, and this paper furnishes a key step in this direction by treating an extensive catalog of degeneracies and Fuchsian structures.

The structures studied in this article are found in numerous equations of physical interest, with one key example being the Einstein equations of general relativity near big bang (and more general spacelike) singularities.
One potential application of a nonlinear extension of our results would be to provide precise asymptotics of solutions of the Einstein equations at such singularities, provided one is in a setting where the nonlinear theory is known a priori to reach the singularity.
Examples of this include recent stability results for Kasner spacetimes (see \cite{fournodavlos_rodnianski_speck_2023} and references within) and for cosmological regions of Kerr-de Sitter spacetimes \cite{fournodavlos_schlue_2024}.

Another potential direction is toward scattering problems for the full Einstein equations, with asymptotic data imposed at the singular time.
Although this question is largely open, we do mention the recent scattering results of \cite{cicortas_2024a, cicortas_2024b} for the Einstein vacuum equations in asymptotically de Sitter settings and the constructions \cite{fournodavlos_luk_2023} of spacetimes with prescribed Kasner-type singularities.

Finally, recall that the systems studied in this article are, in a loose sense, based around model wave equations with solutions expressed in terms of confluent hypergeometric functions.
There do exist many other equations of physical interest with model solutions based on more general Heun functions \cite{hortacsu_2011}; a well-known example from relativity involves wave-type systems that model black hole perturbations about Kerr spacetimes \cite{blandin_pons_marcilhacy_1983, minucci_macedo_2025}.
It would also be of interest to investigate whether methods from this paper can be adapted to analyze these settings.

\subsection{Notations}

We adopt the following standard notations throughout the paper:
\begin{itemize}
\item We let $\N$, $\R$, $\C$ denote the sets of natural numbers (excluding $0$), real numbers, and complex numbers, respectively.
In addition, we set $\N_0 := \N \cup \{ 0 \}$.

\item We let $\C^n$ ($\R^n$) denote the set of all $n$-vectors with complex (real) components.
We also let $\C^n \otimes \C^m$ ($\R^n \otimes \R^m$) denote the set of all $n \times m$ matrices with complex (real) entries.

\item We let $C ( X, Y )$ and $C^k ( X; Y )$ denote the space of all continuous and $C^k$ (respectively) functions mapping from some subset $X \subseteq \R^q$ to some vector space $Y$.

\item We write $a \lesssim b$ and $a \gtrsim b$ to denote $a \leq C b$ and $a \geq C b$ (respectively) for some universal constant $C$; any dependence of $C$ on other parameters will be mentioned explicitly.
We will also write $a \simeq b$ to mean that both $a \lesssim b$ and $a \gtrsim b$ hold.

\item We let ``$\imath$" denote the imaginary number, to distinguish it from our use of ``$i$" as an index.
\end{itemize}

\subsection{Acknowledgments}

The authors thank Spyros Alexakis, Grigorios Fournodavlos, Warren Li, Todd Oliynik, Michael Ruzhansky, and Zoe Wyatt for discussions, as well as the anonymous referee for their careful reading and helpful comments.
B.S.\ is supported by EPSRC grant EP/V005529/2.
A.S.\ was supported by EPSRC grant EP/Y021487/1 for part of this work.

\section{General Hyperbolic Systems} \label{sec.system}

In this section, we study general hyperbolic systems that become degenerate and critically singular at a single time $t = 0$.
We state and prove the main result of the paper, Theorem \ref{thm.energy_main}, which quantifies the precise asymptotic and scattering theories at the critical time.

\subsection{Preliminaries} \label{sec.system_prelim}

We begin by collecting various assumptions that we will impose on our system, as well as some basic definitions that will be useful later.

\subsubsection{The General Setting}

We begin by defining an admissible class of coefficients:

\begin{definition} \label{def.coeff}
Let $T > 0$, let $\alpha > 0$, and let $X$ be a normed vector space.
We define $\mc{T}^\alpha ( [ 0, T ]; X )$ to be the space of all functions $f: [ 0, T ] \rightarrow X$ of the form
\begin{equation}
\label{eq.coeff} f (t) = f_0 + f_\ast (t) \text{,} \qquad t \in [ 0, T ] \text{,}
\end{equation}
where $f_0 \in X$, and where $f_\ast: [ 0, T ] \rightarrow X$ is smooth on $( 0, T ]$ and satisfies, for any $k \in \N$,
\begin{equation}
\label{eq.coeff_est} | ( t \partial_t )^k f_\ast |_X \leq C_k t^\alpha \text{,} \qquad C_k > 0 \text{.}
\end{equation}
\end{definition}

\begin{remark}
In practice, the spaces $\mc{T}^\alpha ( [ 0, T ]; X )$ will contain the various coefficients for the systems that we study.
Note a prototypical example of a function $f \in \mc{T}^\alpha ( [ 0, T ]; \R )$ is
\[
f (t) = f_0 + f_\alpha t^\alpha \text{,} \qquad f_0, f_\alpha \in \R \text{.}
\]
\end{remark}

A preliminary description of our system can now be given as follows:

\begin{assumption}[General Setting] \label{ass.system}
We fix the following quantities pertaining to our setting:
\begin{itemize}
\item Fix the spatial dimension $d \in \N$ and the timespan $T > 0$.

\item Fix $n \in \N$, representing the number of components in our system.

\item Fix $\ell \in ( -1, \infty )^d$, $\alpha > 0$, and $\lambda \in \mc{T}^\alpha ( [ 0, T ]; \R^d \otimes \R^d )$.
We then define the function
\begin{equation}
\label{eq.system_degen} \mf{H}: ( 0, T ]_t \times \R^d_\xi \rightarrow ( 0, \infty ) \text{,} \qquad \mf{H} ( t, \xi ) := \bigg[ \sum_{ i, j = 1 }^d t^{ \ell_i + \ell_j } \, \lambda_{ij} (t) \, \xi_i \xi_j \bigg]^\frac{1}{2} \text{,}
\end{equation}
representing the \emph{degenerate hyperbolicity} of our system.
To ensure that \eqref{eq.system_degen} is well-defined, we additionally assume the following uniform ellipticity condition:
\begin{equation}
\label{eq.system_ellip} \mf{H}^2 ( t, \xi ) \geq \lambda_0 \sum_{ i = 1 }^d t^{ 2 \ell_i } \, \xi_i^2 \text{,} \qquad \lambda_0 > 0 \text{.}
\end{equation}
\end{itemize}
In general, we consider the following system in the time-frequency domain $( 0, T ]_t \times \R^d_\xi$,
\begin{equation}
\label{eq.system_gen} \partial_t U = \mc{A} \, U + F \text{,}
\end{equation}
where $\smash{ U: ( 0, T ]_t \times \R^d_\xi \rightarrow \C^n }$ and $\smash{ F \in C ( ( 0, T ]_t \times \R^d_\xi; \C^n ) }$ represent spatial Fourier transforms of the unknown and the forcing term of our system, respectively, and $\mc{A}: ( 0, T ]_t \times \R^d_\xi \rightarrow \C^n \otimes \C^n$ denotes the symbol of the operator defining our weakly hyperbolic system.
\end{assumption}

\begin{remark}
In other words, the weakly hyperbolic system under consideration is given by
\begin{equation}
\label{eq.system_actual} \partial_t \check{U} = \mc{A} ( t, - \imath \nabla_x ) \, \check{U} + \check{F} \text{.}
\end{equation}
\end{remark}

We defer our detailed assumptions on $\mc{A}$---in particular on its relation to the degeneracy $\mf{H}$---to further below, as these will vary depending on the particular region in the time-frequency domain.

\begin{remark}
Note in particular that $\mc{A}$ depends on $t$, but not on the spatial coordinates $x$.
\end{remark}

Next, we define the frequency-rescaled time $z$ that captures the scaling of our system:

\begin{definition} \label{def.z}
We define the following \emph{rescaled frequency measure}:
\begin{equation}
\label{eq.z_factor} \zfac := \max_{ 1 \leq i \leq d } \langle \xi_i \rangle^\frac{1}{ \ell_i + 1 } \text{.}
\end{equation}
We then define the \emph{frequency-rescaled time} $z$ by
\begin{equation}
\label{eq.z} z: [ 0, T ]_t \times \R^d_\xi \rightarrow \R \text{,} \qquad z ( t, \xi ) := \zfac \, t \text{,}
\end{equation}
as well as the similarly \emph{rescaled hyperbolic degeneracy}:
\begin{equation}
\label{eq.z_degen} \mf{Z}: [ 0, T ]_t \times \R^d_\xi \rightarrow \R \text{,} \qquad \mf{Z} ( t, \xi ) := \zfac^{-1} \, \mf{H} ( t, \xi ) \text{.}
\end{equation}
Furthermore, for future convenience, we define the parameter
\begin{equation}
\label{eq.ell} \ell_\ast := \min ( \ell_1, \dots, \ell_d ) \text{.}
\end{equation}
\end{definition}

In our technical analysis, we will often work with the rescaled coordinates $( z, \xi )$, instead of with $( t, \xi )$.
Observe that in these transformed coordinates, we have
\begin{equation}
\label{eq.z_coord} \partial_z = \zfac^{-1} \, \partial_t \text{,} \qquad d z = \zfac \, dt \text{.}
\end{equation}

Next, we define the appropriate partition of $( 0, T ]_t \times \R_\xi$ into microlocal zones:

\begin{definition} \label{def.zones}
Fix $0 < \rho_0 \ll_T 1$, and define the following regions within $( 0, T ]_t \times \R^d_\xi$,
\begin{equation}
\label{eq.zones} Z_P := \{ z \leq \rho_0 \} \text{,} \qquad Z_I := \{ \rho_0 \leq z \leq \rho_0^{-1} \} \text{,} \qquad Z_H := \{ z \geq \rho_0^{-1} \} \text{,}
\end{equation}
called the \emph{pseudodifferential zone}, \emph{intermediate zone}, and \emph{hyperbolic zone}, respectively.
\end{definition}

\begin{proposition} \label{thm.zones_ass}
The following hold for sufficiently small $\rho_0$ (depending on $T$):
\begin{equation}
\label{eq.zones_ass} Z_P \subseteq \{ t \leq \rho_0 \} \text{,} \qquad Z_H \subseteq \big\{ | \xi | \geq \tfrac{1}{2} \big( \tfrac{1}{ T \rho_0 } \big)^{ \ell_\ast + 1 } \big\} \text{.}
\end{equation}
\end{proposition}

\begin{proof}
First, if $( t, \xi ) \in Z_P$, then by \eqref{eq.z} and \eqref{eq.zones},
\[
t \leq \zfac^{-1} z \leq \rho_0 \text{,}
\]
proving the first part of \eqref{eq.zones_ass}.
Similarly, if $( t, \xi ) \in Z_H$, then by \eqref{eq.z_factor}, \eqref{eq.z} and \eqref{eq.zones},
\[
\langle \xi \rangle \geq \zfac^{ \ell_\ast + 1 } \geq \big( \tfrac{1}{ T \rho_0 } \big)^{ \ell_\ast + 1 } \text{.}
\]
The second part of \eqref{eq.zones_ass} follows, since $\langle \xi \rangle \leq 2 | \xi |$ for sufficiently large $\xi$.
\end{proof}

\begin{remark}
As a result of Proposition \ref{thm.zones_ass}, we will, in practice, always choose $\rho_0 \ll_T 1$ so that
\[
Z_P \subseteq \{ t \ll 1 \} \text{,} \qquad Z_H \subseteq \{ | \xi | \gg 1 \} \text{.}
\]
\end{remark}

Similar to Definition \ref{def.coeff}, it will be useful to have spaces for various ``remainder" quantities:

\begin{definition} \label{def.symbol}
Let $\mc{V} \subseteq ( 0, T ]_t \times \R^d_\xi$, let $X$ be a normed vector space, and fix $a \in \R$.
\begin{itemize}
\item Let $\mc{S}^a ( \mc{V}; X )$ be the space of all functions $f \in C( \mc{V}; X )$ such that $f ( \cdot, \xi )$ is smooth (in $t$) for any $\xi \in \R^d$, and $f$ satisfies the following uniform bounds on $\mc{V}$ for all $k \in \N_0$:
\begin{equation}
\label{eq.symbol} | ( z \partial_z )^k f |_X \leq C_k z^a \text{,} \qquad C_k > 0 \text{.}
\end{equation}

\item Let $\mc{S}^a_\ast ( \mc{V}; X )$ be the space of all functions $f \in C ( \mc{V}; X )$ such that $f ( \cdot, \xi )$ is smooth (in $t$) for any $\xi \in \R^d$, and there exists $p \in \N_0$ such that the following hold on $\mc{V}$ for all $k \in \N_0$:
\begin{equation}
\label{eq.symbol_poly} | ( z \partial_z )^k f |_X \leq C_k z^a ( 1 + | \log z | )^p \text{,} \qquad C_k > 0 \text{.}
\end{equation}
\end{itemize}
\end{definition}

\begin{remark}
By \eqref{eq.z} and \eqref{eq.z_coord}, the left-hand sides of \eqref{eq.symbol} and \eqref{eq.symbol_poly} can also be written
\[
| ( z \partial_z )^k f |_X = | ( t \partial_t )^k f |_X \text{.}
\]
\end{remark}

\begin{remark}
Note the following inclusions hold for any $\delta > 0$ (with $a$, $X$ as in Definition \ref{def.symbol}):
\begin{align*}
\mc{S}^a ( Z_P; X ) \subseteq \mc{S}^a_\ast &( Z_P; X ) \subseteq \mc{S}^{ a - \delta } ( Z_P; X ) \text{,} \\
\mc{S}^{ a  } ( Z_H; X ) \subseteq \mc{S}^a_\ast &( Z_H; X ) \subseteq \mc{S}^{ a + \delta } ( Z_H; X ) \text{.}
\end{align*}
\end{remark}

\begin{proposition} \label{thm.degen_est}
For sufficiently small $\rho_0$ (depending on $T$), the following hold:
\begin{align}
\label{eq.degen_est} t \mf{H} \in \mc{S}^{ \ell_\ast + 1 } ( Z_P \cup Z_I; \R ) \text{,} &\qquad ( t \mf{H} )^{-1} \in \mc{S}^{ - ( \ell_\ast + 1 ) } ( Z_H; \R ) \text{,} \\
\notag \mf{Z} \in \mc{S}^{ \ell_\ast } ( Z_P \cup Z_I; \R ) \text{,} &\qquad \mf{Z}^{-1} \in \mc{S}^{ -\ell_\ast } ( Z_H; \R ) \text{.}
\end{align}
\end{proposition}

\begin{proof}
First, by \eqref{eq.system_degen}--\eqref{eq.system_ellip}, we have (everywhere on $( 0, T ]_t \times \R^n_\xi$) that
\begin{equation}
\label{eql.degen_est_0} ( t \mf{H} )^2 \simeq \sum_{ i = 1 }^d t^{ 2 ( \ell_i + 1 ) } \xi_i^2
\end{equation}
Moreover, differentiating \eqref{eq.system_degen} and recalling Definition \ref{def.coeff}, \eqref{eq.z}, and \eqref{eq.z_coord}, we obtain
\begin{align}
\label{eql.degen_est_1} \big| ( z \partial_z )^k [ ( t \mf{H} )^2 ] \big| &= \big| ( t \partial_t )^k [ ( t \mf{H} )^2 ] \big| \\
\notag &\lesssim \sum_{ i = 1 }^d t^{ 2 ( \ell_i + 1 ) } \xi_i^2
\end{align}
for any $k \in \N$, where the constant in \eqref{eql.degen_est_1} also depends on $k$.

Let us restrict now to $Z_P \cup Z_I$.
Here, by Definitions \ref{def.z}--\ref{def.zones}, we have
\begin{align}
\label{eql.degen_est_10} \sum_{ i = 1 }^d t^{ 2 ( \ell_i + 1 ) } \xi_i^2 &\lesssim \sum_{ i = 1 }^d t^{ 2 ( \ell_i + 1 ) } \zfac^{ 2 ( \ell_i + 1 ) } \\
\notag &\lesssim z^{ 2 ( \ell_\ast + 1 ) } \text{.}
\end{align}
Combining \eqref{eql.degen_est_0}--\eqref{eql.degen_est_10}, we then obtain, for any $k \in \N_0$,
\begin{align*}
\big| ( z \partial_z )^k ( t \mf{H} ) \big| &\lesssim \bigg[ \sum_{ i = 1 }^d t^{ 2 ( \ell_i + 1 ) } \xi_i^2 \bigg]^\frac{1}{2} \\
&\lesssim z^{ \ell_\ast + 1 } \text{,}
\end{align*}
in particular proving the first part of \eqref{eq.degen_est}.

Next, we restrict to $Z_H$.
For $\xi \in \R^d$, we let $1 \leq j ( \xi ) \leq n$ such that
\[
\zfac = \langle \xi \cdot e_{ j ( \xi ) } \rangle^\frac{1}{ \ell_j + 1 } \text{.}
\]
Recalling \eqref{eq.system_ellip}, Definition \ref{def.zones}, and that $| \xi | \gg 1$ on $Z_H$, we obtain
\begin{align}
\label{eql.degen_est_11} \sum_{ i = 1 }^d t^{ 2 ( \ell_i + 1 ) } \xi_i^2 &\gtrsim t^{ 2 ( \ell_{ j ( \xi ) } + 1 ) } \zfac^{ 2 ( \ell_{ j ( \xi ) } + 1 ) } \\
\notag &\gtrsim z^{ 2 ( \ell_\ast + 1 ) } \text{.}
\end{align}
Combining \eqref{eql.degen_est_0}--\eqref{eql.degen_est_1} and \eqref{eql.degen_est_11} yields, for any $k \in \N_0$,
\begin{align*}
\big| ( z \partial_z )^k [ ( t \mf{H} )^{-1} ] \big| &\lesssim \bigg[ \sum_{ i = 1 }^d t^{ 2 ( \ell_i + 1 ) } \xi_i^2 \bigg]^{ -\frac{1}{2} } \\
&\lesssim z^{ - ( \ell_\ast + 1 ) } \text{,}
\end{align*}
which completes the proof of the second part of \eqref{eq.degen_est}.

For the remaining properties for $\mf{Z}$ on \eqref{eq.degen_est}, we simply use the properties for $t \mf{H}$ in \eqref{eq.degen_est} that we just proved, and we note from \eqref{eq.z_degen} that $\mf{Z} = z^{-1} ( t \mf{H} )$.
\end{proof}

The following integration lemma will be useful later in our analysis:

\begin{proposition} \label{thm.time_int}
Let $\mc{V} \subseteq ( 0, T ]_t \times \R^d_\xi$ and $a \in \R$.
Given $g \in \mc{S}^a_\ast ( \mc{V}; \C )$, there exists $\mc{I} g \in \mc{S}^{ a + 1 }_\ast ( \mc{V}; \C )$ such that $\partial_z ( \mc{I} g ) = g$.
Furthermore, if $g \in \mc{S}^a ( \mc{V}; \C )$ and $a \neq -1$, then $\mc{I} g \in \mc{S}^{ a + 1 } ( \mc{V}; \C )$.
\end{proposition}

\begin{proof}
First, we extend $g$ to $( 0, \infty )_t \times \R^d_\xi$ by defining its values to be zero outside $\mc{V}$.
Working now in $( z, \xi )$-coordinates, we can then construct $\mc{I} g$ as follows:
\[
( \mc{I} g ) ( z, \xi ) := \begin{cases} \int_0^z g ( \zeta, \xi ) \, d \zeta & a > -1 \text{,} \\
\int_\infty^z g ( \zeta, \xi ) \, d \zeta & a < -1 \text{,} \\ \int_1^z g ( \zeta, \xi ) \, d \zeta & a = -1 \text{.} \end{cases}
\]
Note the above can pick up an additional power of $\log z$ only when $a = -1$.
\end{proof}

\subsubsection{Detailed Assumptions}

We now describe the precise assumptions we will impose on our system \eqref{eq.system_gen}.
In particular, we will impose different assumptions on each of the regions $Z_I$, $Z_P$, $Z_H$.
These reflect the fact that we must diagonalize our system differently in each of the regions.

\begin{assumption}[Assumptions on $Z_I$] \label{ass.system_i}
$\mc{A}$ satisfies the following on $Z_I$:
\begin{equation}
\label{eq.A_i} \mc{A} \in \zfac \mc{S}^0 ( Z_I; \C^n \otimes \C^n ) \text{.}
\end{equation}
\end{assumption}

\begin{remark}
In practice, Assumption \ref{ass.system_i} will trivially hold, since $z \simeq 1$ on $Z_I$ by \eqref{eq.zones}.
\end{remark}

Note that under Assumption \ref{ass.system_i}, our system on $Z_I$ can be written as
\begin{equation}
\label{eq.system_i} \partial_z U = R_I U + \zfac^{-1} F \text{,} \qquad R_I := \zfac^{-1} \mc{A} \in \mc{S}^0 ( Z_I; \C^n \otimes \C^n ) \text{.}
\end{equation}

\begin{assumption}[Assumptions on $Z_P$] \label{ass.system_p}
$\mc{A}$ can be expressed on $Z_P$ as
\begin{equation}
\label{eq.A_p} M_P \mc{A} M_P^{-1} + ( \partial_t M_P ) M_P^{-1} = t^{-1} B_P + \zfac R_P \text{,} 
\end{equation}
where the quantities in \eqref{eq.A_p} satisfy the following:
\begin{itemize}
\item $M_P \in \mc{S}^0 ( Z_P; \C^n \otimes \C^n )$ is invertible, and $M_P^{-1} \in \mc{S}^0 ( Z_P; \C^n \otimes \C^n )$.

\item $B_P \in \mc{S}^0 ( Z_P; \C^n \otimes \C^n )$ is independent of $t$, that is,
\begin{equation}
\label{eq.B_p} B_P ( t, \xi ) = B_{ P, 0 } ( \xi ) \text{,} \qquad B_{ P, 0 }: \R_\xi^d \rightarrow \C^n \otimes \C^n \text{.}
\end{equation}

\item $B_{ P, 0 }$ is everywhere in Jordan normal form.

\item $R_P \in \mc{S}^{ a_P }_\ast ( Z_P; \C^n \otimes \C^n )$ for some $a_P > -1$.
\end{itemize}
\end{assumption}

\begin{remark}
The assumption \eqref{eq.B_p} is imposed only for convenience, as one can also treat $B_P$'s that are $t$-dependent.
In this more general setting, one can in practice replace $B_P$ by the $t$-independent $B_P ( 0, \cdot )$ and treat the difference $t^{-1} [ B_P ( t, \cdot ) - B_P ( 0, \cdot ) ]$ as part of the remainder term $\zfac R_P$.

For instance, if the $B_{ P; ij } ( \cdot, \xi )$ are in the coefficient class $\mc{T}^\alpha ( [ 0, T ]; \C^n \otimes \C^n )$ for any $\xi \in \R^d$ (see Definition \ref{def.coeff}), then $t^{-1} [ B_P ( t, \cdot ) - B_P ( 0, \cdot ) ]$ lies in $\zfac \mc{S}^{ -1 + \alpha } ( Z_P; \C^n \otimes \C^n )$. 
\end{remark}

\begin{remark}
In our upcoming applications, $B_P$ will be a constant matrix.
\end{remark}

Under Assumption \ref{ass.system_p} (and recalling \eqref{eq.z}, \eqref{eq.z_coord}), our system on $Z_P$ can then be written
\begin{align}
\label{eq.system_p} \partial_z U_P &= ( z^{-1} B_P + R_P ) U_P + \zfac^{-1} F_P \text{,} \\
\notag ( U_P, F_P ) :\!\!&= ( M_P U, M_P F ) \text{.}
\end{align}
Here, $z^{-1} B_P$ captures the critically singular part of the system, while $R_P$ has remainder terms.

\begin{assumption}[Assumptions on $Z_H$] \label{ass.system_h}
$\mc{A}$ can be expressed on $Z_H$ as
\begin{equation}
\label{eq.A_h} M_H \mc{A} M_H^{-1} + ( \partial_t M_H ) M_H^{-1} = \imath \mf{H} \, D_H + t^{-1} B_H + \zfac R_H \text{,}
\end{equation}
where the terms on the right-hand side of \eqref{eq.A_h} satisfy the following:
\begin{itemize}
\item $M_H \in \mc{S}^0 ( Z_H; \C^n \otimes \C^n )$ is invertible, and $M_H^{-1} \in \mc{S}^0 ( Z_H; \C^n \otimes \C^n )$.

\item $B_H \in \mc{S}^0 ( Z_H; \C^n \otimes \C^n )$ is independent of $| \xi |$, i.e.
\begin{equation}
\label{eq.B_h} B_H ( t, \xi ) = B_{ H, \infty } \big( t, \tfrac{ \xi }{ | \xi | } \big) \text{,} \qquad \xi \neq 0 \text{,} \quad B_{ H, \infty }: ( 0, T ] \times \Sph^{ d - 1 } \rightarrow \C^n \otimes \C^n \text{.}
\end{equation}

\item $D_H \in \mc{S}^0 ( Z_H; \C^n \otimes \C^n )$ is everywhere diagonal and real-valued.

\item $R_H \in \mc{S}^{ a_H } ( Z_H; \C^n \otimes \C^n )$ for some $a_H < -1$.
\end{itemize}
\end{assumption}

\begin{remark}
Again, the assumption \eqref{eq.B_h} is for convenience, as one can also treat many $| \xi |$-dependent $B_H$'s.
In practice, one replaces $B_H$ by its infinite-frequency limit $B_{ H, \infty } := \lim_{ | \xi | \nearrow \infty } B_H$ and treats the difference $t^{-1} [ B_H ( t, \cdot ) - B_{ H, \infty } ( t, \cdot ) ]$ as part of the remainder term $\zfac R_H$.
\end{remark}

Under Assumption \ref{ass.system_h} (and using \eqref{eq.z}, \eqref{eq.z_coord}), our system on $Z_H$ becomes
\begin{align}
\label{eq.system_h} \partial_z U_H &= ( \imath \mf{Z} \, D_H + z^{-1} B_H + R_H ) U_H + \zfac^{-1} F_H \text{,} \\
\notag ( U_H, F_H ) :\!\!&= ( M_H U, M_H F ) \text{.}
\end{align}
Note that $\imath \mf{Z} \, D_H$ captures the hyperbolicity of \eqref{eq.system_h} at $t > 0$ and its degeneration at $t = 0$, while $z^{-1} B_H$ and $R_H$ capture the critically singular part and the remainder, respectively.

For the main result, Theorem \ref{thm.energy_main}, as well as for most of our key examples and applications in later sections, we will impose one additional condition on $Z_H$ on top of Assumption \ref{ass.system_h}:

\begin{definition} \label{def.system_h_strict}
We say \eqref{eq.system_gen} is \emph{semi-strictly hyperbolic} iff there exists $d_0 > 0$ with
\begin{equation}
\label{eq.system_h_strict} | D_{ H, ii } - D_{ H, jj } | \geq d_0 \text{,} \qquad 1 \leq i, j \leq n \text{,} \quad i \neq j \text{.}
\end{equation}
\end{definition}

\subsubsection{The Intermediate Zone}

We first restrict our attention to $Z_I$, for which the analysis is the most trivial, since $z$ is bounded both from above and from below on $Z_I$.

\begin{proposition} \label{thm.energy_i}
For any $\xi \in \R^d$ and $0 \leq t_0, t_1 \leq T$ with $( t_0, \xi ), ( t_1, \xi ) \in Z_I$, we have
\begin{equation}
\label{eq.energy_i} | U ( t_1, \xi ) | \lesssim | U ( t_0, \xi ) | + \bigg| \int_{ t_0 }^{ t_1 } | F ( \tau, \xi ) | \, d \tau \bigg| \text{.}
\end{equation}
\end{proposition}

\begin{proof}
For convenience, we set $z_0 := \zfac \, t_0$ and $z_1 := \zfac \, t_1$, and we adopt $( z, \xi )$-coordinates throughout the proof.
Since \eqref{eq.zones}, Assumption \ref{ass.system_i}, and \eqref{eq.system_i} imply
\[
\int_{ \rho_0 }^{ \rho_0^{-1} } | R_I ( \zeta, \xi ) | \, d \zeta \lesssim 1 \text{,}
\]
then applying the Gronwall inequality and the above to \eqref{eq.system_i} yields
\[
| U ( z_1, \xi ) | \lesssim | U ( z_0, \xi ) | + \zfac^{-1} \bigg| \int_{ z_0 }^{ z_1 } | F ( \zeta, \xi ) | \, d \zeta \bigg| \text{.}
\]
Rewriting the above in $( t, \xi )$-coordinates and recalling \eqref{eq.z_coord} results in \eqref{eq.energy_i}.
\end{proof}

\subsection{The Pseudodifferential Zone} \label{sec.system_zp}

Next, we turn our attention to $Z_P$.
In this case, we generally cannot treat \eqref{eq.system_p} directly, as the integrating factor from $z^{-1} B_P$ leads to a power of $z$ that can make $R_P$ non-integrable.
Thus, we will need additional renormalizations to improve the error.

\subsubsection{Higher-Order Renormalization}

Before discussing the renormalization itself, we first present the following technical lemma that will play a key role in this process.

\begin{lemma} \label{thm.perf_diag_pN}
For $a \in \R$ and $Y \in \mc{S}^a_\ast ( Z_P; \C^n \otimes \C^n )$, there exists $N \in \mc{S}^{ a + 1 }_\ast ( Z_P; \C^n \otimes \C^n )$ so that
\begin{equation}
\label{eq.perf_diag_pN} \partial_z N + [ N, z^{-1} B_P ] = Y \text{.}
\end{equation}
\end{lemma}

\begin{proof}
It suffices to construct $N_{ij} \in \mc{S}^{ a + 1 }_\ast ( Z_P; \C )$ for every fixed $1 \leq i, j \leq n$.
Note that for such $i$ and $j$, the corresponding component of \eqref{eq.perf_diag_pN} can be written as
\begin{equation}
\label{eql.perf_diag_pN_0} \partial_z N_{ij} + z^{-1} ( B_{ P, jj } - B_{ P, ii } ) N_{ij} = Y_{ ij } + z^{-1} ( B_{ P, i (i+1) } N_{ (i+1) j } - N_{ i (j-1) } B_{ P, (j-1) j } ) \text{,}
\end{equation}
where we have, for brevity, also adopted the conventions
\begin{equation}
\label{eql.perf_diag_pN_1} B_{ P, i (n+1) } = B_{ P, 0j } = N_{ (n+1) j } = N_{i0} = 0 \text{.}
\end{equation}

Observe (note $B_P$ is $z$-independent) that \eqref{eql.perf_diag_pN_0} can be rewritten as
\begin{align}
\label{eql.perf_diag_pN_2} \partial_z ( z^{ q_{ ij } } N_{ ij } ) &= z^{ q_{ ij } } Y_{ ij } + z^{ q_{ij} - 1 } ( B_{ P, i (i+1) } N_{ (i+1) j } - N_{ i (j-1) } B_{ P, (j-1) j } ) \text{,} \\
\notag q_{ ij } :\!\!&= B_{ P, jj } ( \xi ) - B_{ P, ii } ( \xi ) \text{.}
\end{align}
As a result, we can construct $N_{ij}$ recursively via the relation
\begin{equation}
\label{eql.perf_diag_pN_3} N_{ ij } := z^{ -q_{ij} } \mc{I} ( z^{ q_{ij} - 1 } ( B_{ P, i (i + 1) } N_{ (i+1) j } - B_{ P, (j-1) j } N_{ i (j-1) } ) ) + z^{ -q_{ij} } \mc{I} ( z^{ q_{ij} } Y_{ ij } ) \text{.}
\end{equation}
	
The construction proceeds via a nested induction, starting with the bottom row.
First, note that $N_{n1}$ (the bottom-left entry) is well-defined by Proposition \ref{thm.time_int}, \eqref{eql.perf_diag_pN_1}, and \eqref{eql.perf_diag_pN_3}:
\begin{equation}
\label{eql.perf_diag_pN_10} N_{n1} = z^{ -q_{n1} } \mc{I} ( z^{ q_{n1} } Y_{n1} ) \in \mc{S}^{ a + 1 }_\ast ( Z_P; \C ) \text{.}
\end{equation}
Next, for the remaining elements in the bottom row, we have, from \eqref{eql.perf_diag_pN_1} and \eqref{eql.perf_diag_pN_3},
\[
N_{nj} = z^{ -q_{nj} } \mc{I} ( z^{ q_{nj} } Y_{nj} ) - z^{ -q_{nj} } \mc{I} ( z^{ q_{nj} - 1 } B_{ P, (j-1) j } N_{ n (j-1) } ) \text{,} \qquad 1 < j \leq n \text{.}
\]
Note the right-hand side of the above only contains components of $N$ on the bottom row to the left of $N_{nj}$.
Consequently, iterating through the bottom row from left to right, we see that each $N_{nj}$ is well-defined.
Moreover, since we have already obtained $N_{ n (j-1) } \in \mc{S}^{ a + 1 }_\ast ( Z_P; \C )$ at an earlier step of the iteration, then applying Lemma \ref{thm.time_int} to the above yields $N_{nj} \in \mc{S}^{ a + 1 }_\ast ( Z_P; \C )$.

Once the bottom row of $N$ is obtained, we can then similarly iterate through each remaining row, from bottom to top; for each row, we then iterate through the entries from left to right.
In particular, for each $1 \leq i < n$ and $1 \leq j \leq n$, we see that the right-hand side of \eqref{eql.perf_diag_pN_2} only contains entries of $N$ that are below and/or to the left of $N_{ij}$---entries that have already been determined earlier in the process to lie in $\mc{S}^{ a + 1 }_\ast ( Z_P; \C )$.
Thus, it follows that $N_{ij}$ is well-defined, and applying Proposition \ref{thm.time_int} to \eqref{eql.perf_diag_pN_2} yields $N_{ij} \in \mc{S}^{ a + 1 }_\ast ( Z_P; \C )$, as desired.
\end{proof}

The following lemma describes precisely our renormalization of \eqref{eq.system_p}, up to any finite order:

\begin{lemma} \label{thm.perf_diag_p}
There exist sequences of matrix-valued functions
\begin{align}
\label{eq.perf_diag_p_DNR} \pd{D_P}{m} &\in \mc{S}^{ m ( a_P + 1 ) - 1 }_\ast ( Z_P; \C^n \otimes \C^n ) \text{,} \\
\notag \pd{N_P}{m} &\in \mc{S}^{ m ( a_P + 1 ) }_\ast ( Z_P; \C^n \otimes \C^n ) \text{,} \\
\notag \pd{R_P}{m} &\in \mc{S}^{ m ( a_P + 1 ) + a_P }_\ast ( Z_P; \C^n \otimes \C^n ) \text{,}
\end{align}
for all $m \in \N$, such that the following properties hold:
\begin{itemize}
\item Each $\pd{D_P}{m}$, $m \in \N$, is everywhere diagonal.

\item Each $\pd{R_P}{m}$, $m \in \N$, is given by the following:
\begin{align}
\label{eq.perf_diag_p_RQ} \pd{R_P}{m} &= \sum_{ k = 1 }^m \partial_z \pd{N_P}{k} + \pd{Q_P}{m} ( z^{-1} B_P + R_P ) - \bigg( z^{-1} B_P + \sum_{ k = 1 }^m \pd{D_P}{k} \bigg) \pd{Q_P}{m} \text{,} \\
\notag \pd{Q_P}{m} :\!\!&= I_n + \sum_{ k = 1 }^m \pd{N_P}{k} \in \mc{S}^0 ( Z_P; \C^n \otimes \C^n ) \text{.}
\end{align}
\end{itemize}
Furthermore, given $m \in \N$, if $\rho_0$ is sufficiently small (with respect to $m$), then:
\begin{itemize}
\item $\pd{Q_P}{m}$ is invertible, and both $\pd{Q_P}{m}$, $( \pd{Q_P}{m} )^{-1}$ are uniformly bounded.

\item The following system of differential equations holds on $Z_P$:
\begin{align}
\label{eq.perf_diag_p} \partial_z \pd{U_P}{m} &= \bigg( z^{-1} B_P + \sum_{ k = 1 }^m \pd{D_P}{k} \bigg) \pd{U_P}{m} + \pd{R_P}{m} ( \pd{Q_P}{m} )^{-1} \pd{U_P}{m} + \zfac^{-1} \pd{F_P}{m} \text{,} \\
\notag ( \pd{U_P}{m}, \pd{F_P}{m} ) :\!\!&= ( \pd{Q_P}{m} U_P, \pd{Q_P}{m} F_P ) \text{.}
\end{align}
\end{itemize}
\end{lemma}

\begin{proof}
For convenience, we begin by setting (see Assumption \ref{ass.system_p})
\begin{equation}
\label{eql.perf_diag_p_0} \pd{R_P}{0} := R_P \in \mc{S}^{ a_P }_\ast ( Z_P; \C^n \otimes \C^n ) \text{.}
\end{equation}
We now construct the $( \pd{D_P}{m}, \pd{N_P}{m}, \pd{R_P}{m} )$'s inductively over $m \in \N$.
Fix $m \in \N$, and suppose we have constructed $( \pd{D_P}{k}, \pd{N_P}{k} )$ for all $1 \leq k < m$ and $\pd{R_P}{k}$ for all $0 \leq k < m$, satisfying \eqref{eq.perf_diag_p_DNR}--\eqref{eq.perf_diag_p_RQ} and \eqref{eql.perf_diag_p_0}, and with each $\pd{D_P}{k}$ everywhere diagonal on $Z_P$.

Recalling \eqref{eq.perf_diag_p_DNR} and \eqref{eql.perf_diag_p_0} (for $m = 1$), we can define $\pd{D_P}{m} \in \mc{S}^{ m ( a_P + 1 ) - 1 }_\ast ( Z_P; \C^n \otimes \C^n )$ by
\begin{equation}
\label{eql.perf_diag_p_1} \pd{D_P}{m} := \operatorname{diag} ( \pd{R_{ P, 11 }}{m-1}, \dots, \pd{R_{ P, nn }}{m - 1} ) \text{,}
\end{equation}
which is diagonal.
We then apply Lemma \ref{thm.perf_diag_pN}, with $Y := \pd{D_P}{m} - \pd{R_P}{m - 1}$, to conclude---with the aid of \eqref{eq.perf_diag_p_DNR} and \eqref{eql.perf_diag_p_1}---that there exists $\pd{N_P}{m} \in \mc{S}^{ m ( a_P + 1 ) }_\ast ( Z_P; \C^n \otimes \C^n )$ satisfying
\begin{equation}
\label{eql.perf_diag_p_2} \partial_z \pd{N_P}{m} + [ \pd{N_P}{m}, z^{-1} B_P ] = \pd{D_P}{m} - \pd{R_P}{m-1} \text{.}
\end{equation}

Define now $( \pd{Q_P}{m}, \pd{R_P}{m} )$ by the formulas \eqref{eq.perf_diag_p_RQ}.
Recalling our inductive hypotheses, along with \eqref{eq.perf_diag_p_RQ}, \eqref{eql.perf_diag_p_0} (when $m = 1$), and \eqref{eql.perf_diag_p_2}, we obtain, from a direct computation,
\begin{align*}
\pd{R_P}{m} &= ( \pd{R_P}{m-1} - \pd{D_P}{m} + [ \pd{N_P}{m}, z^{-1} B_P ] + \partial_z \pd{N_P}{m} ) \\
&\qquad + \pd{N_P}{m} R_P - \sum_{ k = 1 }^{ m - 1 } \pd{D_P}{k} \pd{N_P}{m} - \pd{D_P}{m} \sum_{ k = 1 }^m \pd{N_P}{k} \\
&= \pd{N_P}{m} R_P - \sum_{ k = 1 }^{ m - 1 } \pd{D_P}{k} \pd{N_P}{m} - \pd{D_P}{m} \sum_{ k = 1 }^m \pd{N_P}{k} \text{.}
\end{align*}
Combining \eqref{eq.perf_diag_p_DNR} and \eqref{eql.perf_diag_p_1}--\eqref{eql.perf_diag_p_2} with the above yields $\pd{R_P}{m} \in \mc{S}^{ m ( a_P + 1 ) + a_P }_\ast ( Z_P; \C^n \otimes \C^n )$.
This completes the induction, in particular establishing \eqref{eq.perf_diag_p_DNR}--\eqref{eq.perf_diag_p_RQ} for all $m \in \N$.

Finally, let us fix $m \in \N$.
Then, setting $\rho_0$ to be sufficiently small, each $| \pd{N_P}{k} |$, $1 \leq k \leq m$, can be made arbitrarily uniformly small on $Z_P$ (since $k ( a_P + 1 ) > 0$).
Therefore, from \eqref{eq.perf_diag_p_RQ}, we see that $\pd{Q_P}{m}$ is invertible, and its inverse is also uniformly bounded.
By \eqref{eq.system_p} and \eqref{eq.perf_diag_p_RQ},
\begin{align*}
\partial_z \pd{U_P}{m} &= \sum_{ k = 1 }^m \partial_z \pd{N_P}{k} \, U_P + \pd{Q_P}{m} ( z^{-1} B_P + R_P ) U_P + \zfac^{-1} \pd{Q_P}{m} F_P \\
&= \bigg( z^{-1} B_P + \sum_{ k = 1 }^m \pd{D_P}{k} \bigg) \pd{U_P}{m} + \pd{R_P}{m} ( \pd{Q_P}{m} )^{-1} \pd{U_P}{m} + \zfac^{-1} \pd{F_P}{m} \text{,}
\end{align*}
which yields \eqref{eq.perf_diag_p} and completes the proof of the lemma.
\end{proof}

\begin{definition}
For convenience, in addition to Lemma \ref{thm.perf_diag_p}, we also define
\begin{equation}
\label{eq.triv_diag_p} \pd{U_P}{0} := U_P \text{,} \qquad \pd{Q_P}{0} := I_n \text{,} \qquad \pd{R_P}{0} := R_P \text{.}
\end{equation}
\end{definition}

\begin{remark}
Note that \eqref{eq.perf_diag_p} in the case $m = 0$ reduces to the unrenormalized system \eqref{eq.system_p}.
In the following development, we will leave open the possibility of setting $m = 0$, corresponding to settings where no higher-order renormalization on $Z_P$ is necessary.
\end{remark}

\subsubsection{Estimates for Solutions}

We now use the system \eqref{eq.perf_diag_p} to prove our main estimate on $Z_P$.
In the following, we recall the quantity $B_{ P, 0 }$ that was defined in Assumption \ref{ass.system_p}.

\begin{definition} \label{thm.system_Ep}
For notational convenience, we define the quantity
\begin{equation}
\label{eq.system_Ep} \mc{E}_P := \exp ( -\log z \, B_P ) \text{.}
\end{equation}
\end{definition}

\begin{lemma} \label{thm.explicit_Ep}
Given any $\xi \in \R^d$, $k \in \N$, and Jordan block $\mf{B} \in \C^k \otimes \C^k$ of $B_{ P, 0 } ( \xi )$, that is,
\begin{equation}
\label{eq.explicit_Ep_jordan} \mf{B}_{ij} = \begin{cases} b \in \C & i = j \text{,} \\ 1 & j = i + 1 \text{,} \\ 0 & \text{otherwise,} \end{cases} \qquad n_0 < i, j \leq n_0 + k \text{,} \quad 0 \leq n_0 \leq n - k \text{,}
\end{equation}
the corresponding entries of $\mc{E}_P$ and $\mc{E}_P^{-1}$ are given, for any $n_0 < i, j \leq n_0 + k$, by
\begin{equation}
\label{eq.explicit_Ep} \mc{E}_{ P, ij } ( z, \xi ) = \begin{cases} \frac{ ( - \log z )^{ j - i } z^{-b} }{ ( j - i )! } & i \leq j \text{,} \\ 0 & \text{otherwise,} \end{cases} \qquad \mc{E}^{-1}_{ P, ij } ( z, \xi ) = \begin{cases} \frac{ ( \log z )^{ j - i } z^b }{ ( j - i )! } & i \leq j \text{,} \\ 0 & \text{otherwise.} \end{cases}
\end{equation}
All other entries of $\mc{E}_P ( z, \xi )$ and $\mc{E}_P^{-1} ( z, \xi )$ not along a Jordan block of $B_{ P, 0 } ( \xi )$ vanish identically.
\end{lemma}

\begin{proof}
This is a direct computation using the definition of matrix exponentials.
\end{proof}

\begin{remark}
When $B_{ P, 0 }$ is diagonal, $\mc{E}_P$ and $\mc{E}_P^{-1}$ have the explicit forms
\begin{equation}
\label{eq.explicit_Ep_diag} \mc{E}_P = \operatorname{diag} \big( z^{ - B_{ P, 11 } }, \dots, z^{ - B_{ P, nn } } \big) \text{,} \qquad \mc{E}_P^{-1} = \operatorname{diag} \big( z^{ B_{ P, 11 } }, \dots, z^{ B_{ P, nn } } \big) \text{.}
\end{equation}
\end{remark}

\begin{lemma} \label{thm.energy_Ep}
For any $m \in \N$ and $1 \leq i, j \leq n$, there exists $p \in \N_0$ such that:
\begin{align}
\label{eq.energy_Ep} \big| ( \mc{E}_P \pd{R_P}{m} ( \pd{Q_P}{m} )^{-1} \mc{E}_P^{-1} )_{ ij } ( z, \xi ) \big| &\lesssim z^{ m ( a_P + 1 ) + a_P + ( B_{ P, 0, jj } - B_{ P, 0, ii } ) ( \xi ) } ( 1 + | \log z | )^p \text{,} \\
\notag \big| ( \mc{E}_P \pd{D_P}{m} \mc{E}_P^{-1} )_{ ij } ( z, \xi ) \big| &\lesssim z^{ m ( a_P + 1 ) - 1 } ( 1 + | \log z | )^p \text{.}
\end{align}
\end{lemma}

\begin{proof}
The first property of \eqref{eq.energy_Ep} is a consequence of Lemma \ref{thm.perf_diag_p} and the formulas \eqref{eq.explicit_Ep} for $\mc{E}_P$ and $\mc{E}_P^{-1}$; note in particular that $\pd{R_P}{m} ( \pd{Q_P}{m} )^{-1} \in \mc{S}^{ m ( a_P + 1 ) + a_P }_\ast ( Z_P; \C^n \otimes \C^n )$.

Next, since $\pd{D_P}{m}$ is diagonal and $\mc{E}_P$ is upper triangular (by Lemma \ref{thm.explicit_Ep}), we have
\begin{equation}
\label{eql.energy_Ep_0} ( \mc{E}_P \pd{D_P}{m} \mc{E}_P^{-1} )_{ij} = \sum_{ i \leq k \leq j } \mc{E}_{ P, ik } \pd{D_{ P, kk }}{m} \mc{E}^{-1}_{ P, kj } \text{.}
\end{equation}
From here, we fix a particular $\xi \in \R^d$, and we split into cases:
\begin{itemize}
\item If $i > j$, then \eqref{eql.energy_Ep_0} trivially implies $( \mc{E}_P \pd{D_P}{m} \mc{E}_P^{-1} )_{ij} ( z, \xi ) \equiv 0$.

\item If $i \leq j$, and if $B_{ P, 0, ii } ( \xi ), B_{ P, 0, jj } ( \xi )$ lie in different Jordan blocks of $B_{ P, 0 } ( \xi )$, then for any $i \leq k \leq j$, either $\mc{E}^{-1}_{ P, ik } ( z, \xi )$ or $\mc{E}^{-1}_{ P, kj } ( z, \xi )$ must not lie along a Jordan block of $B_{ P, 0 } ( \xi )$ and hence vanishes by Lemma \ref{thm.explicit_Ep}.
As a result, \eqref{eql.energy_Ep_0} again yields $( \mc{E}_P \pd{D_P}{m} \mc{E}_P^{-1} )_{ij} ( z, \xi ) \equiv 0$.

\item If $i \leq j$ and $B_{ P, 0, ii } ( \xi ), B_{ P, 0, jj } ( \xi )$ lie in the same Jordan block of $B_{ P, 0 } ( \xi )$, then by \eqref{eq.explicit_Ep},
\begin{align*}
| \mc{E}_{ P, ik } ( z, \xi ) | &\lesssim z^{ -B_{ P, 0, ii } ( \xi ) } ( 1 + | \log z | )^{ n - 1 } \text{,} \\
| \mc{E}_{ P, kj }^{-1} ( z, \xi ) | &\lesssim z^{ B_{ P, 0, ii } ( \xi ) } ( 1 + | \log z | )^{ n - 1 } \text{,}
\end{align*}
for any $i \leq k \leq j$, hence the above, along with \eqref{eq.perf_diag_p_DNR} and \eqref{eql.energy_Ep_0}, yields
\[
| ( \mc{E}_P \pd{D_P}{m} \mc{E}_P^{-1} )_{ij} ( z, \xi ) | \lesssim z^{ m ( a_P + 1 ) - 1 } ( 1 + | \log z | )^p \text{,} \qquad p \in \N_0 \text{.}
\]
\end{itemize}
Combining the preceding three cases results in the second part of \eqref{eq.energy_Ep}.
\end{proof}

We can now establish our main estimate and asymptotic limits on $Z_P$:

\begin{proposition} \label{thm.energy_p}
Let $m \in \N_0$ be sufficiently large (depending on $a_P$ and $B_P$), and suppose $\rho_0$ is sufficiently small (with respect to $m$).
In addition, define the quantities
\begin{equation}
\label{eq.energy_p_U} \pd{U_{Pz}}{m} := \mc{E}_P \pd{Q_P}{m} U_P \text{,} \qquad \pd{F_{Pz}}{m} := \mc{E}_P \pd{Q_P}{m} F_P \text{,}
\end{equation}
and assume that $\pd{F_{Pz}}{m}$ is $t$-integrable on $Z_P$:
\begin{equation}
\label{eq.energy_p_ass} \int_0^{ \zfac^{-1} \rho_0 } | \pd{F_{Pz}}{m} ( \tau, \xi ) | \, d \tau < \infty \text{,} \qquad \xi \in \R^d \text{.}
\end{equation}
Then, the following asymptotic limits are both well-defined and finite:
\begin{equation}
\label{eq.energy_p_asymp} \pd{U_{Pz}}{m} ( 0, \xi ) := \lim_{ \tau \searrow 0 } \pd{U_{Pz}}{m} ( \tau, \xi ) \text{,} \qquad \xi \in \R^d \text{.}
\end{equation}
Moreover, for any $\xi \in \R^d$ and $0 \leq t_0, t_1 \leq T$ such that $( t_0, \xi ), ( t_1, \xi ) \in Z_P$, we have
\begin{equation}
\label{eq.energy_p} | \pd{U_{Pz}}{m} ( t_1, \xi ) | \lesssim | \pd{U_{Pz}}{m} ( t_0, \xi ) | + \bigg| \int_{ t_0 }^{ t_1 } | \pd{F_{Pz}}{m} ( \tau, \xi ) | \, d \tau \bigg| \text{.}
\end{equation}
\end{proposition}

\begin{proof}
We once again work in $( z, \xi )$-coordinates, and we define $z_0 := \zfac \, t_0$ and $z_1 := \zfac \, t_1$.
A direct computation using \eqref{eq.perf_diag_p} and \eqref{eq.energy_p_U} then yields that $\pd{U_{Pz}}{m}$ satisfies
\begin{align}
\label{eql.energy_p_0} \partial_z \pd{U_{Pz}}{m} &= \mc{E}_P \bigg[ \sum_{ k = 1 }^m \pd{D_P}{k} + \pd{R_P}{m} ( \pd{Q_P}{m} )^{-1} \bigg] \mc{E}_P^{-1} \pd{U_{Pz}}{m} + \zfac^{-1} \pd{F_{Pz}}{m} \text{,} \\
\notag :\!\!&= \pd{S}{m} \pd{U_{Pz}}{m} + \zfac^{-1} \pd{F_{Pz}}{m} \text{.}
\end{align}

Now, as long as $m$ is large enough, Lemma \ref{thm.energy_Ep} and \eqref{eql.energy_p_0} imply
\begin{equation}
\label{eql.energy_p_1} \int_0^{ \rho_0 } | \pd{S}{m} ( \zeta, \xi ) | \, d \zeta \lesssim 1 \text{.}
\end{equation}
Combining \eqref{eql.energy_p_0}--\eqref{eql.energy_p_1} and recalling the Gronwall inequality, we then obtain
\begin{equation}
\label{eql.energy_p_2} | \pd{U_{Pz}}{m} ( z_1, \xi ) | \lesssim | \pd{U_{Pz}}{m} ( z_0, \xi ) | + \zfac^{-1} \bigg| \int_{ z_0 }^{ z_1 } | \pd{F_{Pz}}{m} ( \zeta, \xi ) | \, d \zeta \bigg| \text{,} \qquad z_0, z_1 > 0 \text{.}
\end{equation}
Furthermore, since both $\pd{S}{m}$ and $\pd{F_{Pz}}{m}$ are $z$-integrable on $Z_P$ (the former due to \eqref{eql.energy_p_1}, and the latter due to \eqref{eq.energy_p_ass}), it follows that the asymptotic limits
\[
\pd{U_{Pz}}{m} ( 0, \xi ) := \lim_{ \zeta \searrow 0 } \pd{U_{Pz}}{m} ( \zeta, \xi ) \text{,} \qquad \xi \in \R^d
\]
exist and are finite; reverting to $( t, \xi )$-coordinates then yields the existence of \eqref{eq.energy_p_asymp}.

Similarly, again by rewriting with respect to $( t, \xi )$-coordinates, the desired estimate \eqref{eq.energy_p} now becomes an immediate consequence of \eqref{eq.energy_p_asymp} and \eqref{eql.energy_p_2}.
\end{proof}

\begin{remark}
By closer inspection of its proof, we see Proposition \ref{thm.energy_p} is applicable with any $m \in \N_0$ such that $( \mc{E}_P \pd{R_P}{m} ( \pd{Q_P}{m} )^{-1} \mc{E}_P^{-1} ) ( \cdot, \xi )$ is integrable in $z$ for all $\xi \in \R^d$.
\end{remark}

\begin{remark}
Note \eqref{eq.energy_p_ass} and \eqref{eql.energy_p_1} also imply a scattering statement---in the setting of Proposition \ref{thm.energy_p}, given $\pd{u_{Pz}}{m}: \R^d_\xi \rightarrow \C^n$, one can find a unique solution $\pd{U_{Pz}}{m}$ to \eqref{eql.energy_p_0} with $\pd{U_{Pz}}{m} ( 0, \cdot ) = \pd{u_{Pz}}{m}$.
\end{remark}

\subsection{The Hyperbolic Zone} \label{sec.system_zh}

We now consider the region $Z_H$, on which we again have to renormalize \eqref{eq.system_h} to higher orders.
However, since $B_H$ lies behind the leading-order hyperbolic part $\imath \mf{Z} D_H$, the higher-order diagonalization process must proceed differently than on $Z_P$.

\subsubsection{Higher-Order Renormalization}

The following lemma gives a precise description of our desired renormalization of \eqref{eq.system_h}, again to any arbitrary finite order:

\begin{lemma} \label{thm.perf_diag_h}
If \eqref{eq.system_gen} is semi-strictly hyperbolic, then there exist
\begin{align}
\label{eq.perf_diag_h_DNR} \pd{D_H}{m} &\in \mc{S}^{ - ( m - 1 ) ( \ell_\ast + 1 ) - 1 } ( Z_H; \C^n \otimes \C^n ) \text{,} \\
\notag \pd{N_H}{m} &\in \mc{S}^{ -m ( \ell_\ast + 1 ) } ( Z_H; \C^n \otimes \C^n ) \text{,} \\
\notag \pd{R_H}{m} &\in \mc{S}^{ -m ( \ell_\ast + 1 ) - 1 } ( Z_H; \C^n \otimes \C^n ) \text{,}
\end{align}
for all $m \in \N$, such that the following properties hold:
\begin{itemize}
\item $\pd{D_{ H, ij }}{m}$ is diagonal for any $m \in \N$.
In particular,
\begin{equation}
\label{eq.perf_diag_h_DD} \pd{D_H}{1} = \operatorname{diag} ( z^{-1} B_{ H, 11 } + R_{ H, 11 }, \dots, z^{-1} B_{ H, nn } + R_{ H, nn } ) \text{.}
\end{equation}

\item Each $\pd{R_H}{m}$, $m \in \N$, is given by the following:
\begin{align}
\label{eq.perf_diag_h_RQ} \pd{R_H}{m} &= \sum_{ k = 1 }^m \partial_z \pd{N_H}{k} + \pd{Q_H}{m} ( \imath \mf{Z} \, D_H + z^{-1} B_H + R_H ) - \bigg( \imath \mf{Z} \, D_H + \sum_{ k = 1 }^m \pd{D_H}{k} \bigg) \pd{Q_H}{m} \text{,} \\
\notag \pd{Q_H}{m} :\!\!&= I_n + \sum_{ k = 1 }^m \pd{N_H}{k} \in \mc{S}^0 ( Z_H; \C^n \otimes \C^n ) \text{.}
\end{align}
\end{itemize}
Furthermore, given $m \in \N$, if $\rho_0$ is sufficiently small (with respect to $m$), then:
\begin{itemize}
\item $\pd{Q_H}{m}$ is invertible, and both $\pd{Q_H}{m}$, $( \pd{Q_H}{m} )^{-1}$ are uniformly bounded on $Z_H$.

\item The following system of differential equations holds on $Z_H$:
\begin{align}
\label{eq.perf_diag_h} \partial_z \pd{U_H}{m} &= \bigg( \imath \mf{Z} \, D_H + \sum_{ k = 1 }^m \pd{D_H}{k} \bigg) \pd{U_H}{m} + \pd{R_H}{m} ( \pd{Q_H}{m} )^{-1} \pd{U_H}{m} + \zfac^{-1} \pd{F_H}{m} \text{,} \\
\notag ( \pd{U_H}{m}, \pd{F_H}{m} ) :\!\!&= ( \pd{Q_H}{m} U_H, \pd{Q_H}{m} F_H ) \text{.}
\end{align}
\end{itemize}
\end{lemma}

\begin{proof}
For convenience, we first set (see \eqref{eq.zones} and Assumption \ref{ass.system_h})
\begin{equation}
\label{eql.perf_diag_h_0} \pd{R_H}{0} := z^{-1} B_H + R_H \in \mc{S}^{-1} ( Z_H; \C^n \otimes \C^n ) \text{.}
\end{equation}
Fix $m \in \N$, and suppose we have $( \pd{D_H}{k}, \pd{N_H}{k} )$ for all $1 \leq k < m$ and $\pd{R_H}{k}$ for all $0 \leq k < m$, satisfying \eqref{eq.perf_diag_h_DNR} and \eqref{eq.perf_diag_h_RQ}, and with $\pd{D_{ H, ij }}{k}$ diagonal for $1 \leq k < m$.
We now proceed inductively.

From \eqref{eq.degen_est}, \eqref{eq.perf_diag_h_DNR}, and \eqref{eql.perf_diag_h_0} (for $m = 1$), we can define $\pd{D_H}{m} \in \mc{S}^{ - ( m - 1 ) ( \ell_\ast + 1 ) - 1 } ( Z_H; \C^n \otimes \C^n )$ (which is in particular diagonal) and $\pd{N_H}{m} \in \mc{S}^{ -m ( \ell_\ast + 1 ) } ( Z_H; \C^n \otimes \C^n )$ (note here we crucially use that our system is semi-strictly hyperbolic) by the formulas
\begin{equation}
\label{eql.perf_diag_h_1} \pd{D_{H, ij}}{m} := \begin{cases} \pd{R_{H, ij}}{m-1} & i = j \text{,} \\ 0 & i \neq j \text{,} \end{cases} \qquad \pd{N_{H, ij}}{m} := \begin{cases} 0 & i = j \text{,} \\ ( \imath \mf{Z} \, ( D_{ H, ii } - D_{ H, jj } ) )^{-1} \pd{R_{H, ij}}{m-1} & i \neq j \text{,} \end{cases}
\end{equation}
for all $1 \leq i, j \leq n$.
A direct computation using \eqref{eq.perf_diag_h_DNR} and \eqref{eql.perf_diag_h_1} then yields
\begin{equation}
\label{eql.perf_diag_h_3} \pd{R_H}{m-1} - \pd{D_H}{m} + [ \pd{N_H}{m}, \imath \mf{Z} \, D_H ] = 0 \text{.}
\end{equation}

Define now $( \pd{Q_H}{m}, \pd{R_H}{m} )$ by the formulas \eqref{eq.perf_diag_h_RQ}.
Recalling our inductive hypotheses, along with \eqref{eq.perf_diag_h_RQ}, \eqref{eql.perf_diag_h_0} (when $m = 1$), and \eqref{eql.perf_diag_h_3}, we then obtain
\begin{align*}
\pd{R_H}{m} &= ( \pd{R_H}{m-1} - \pd{D_H}{m} ) + [ \pd{N_H}{m}, \imath \mf{Z} \, D_H ] + \partial_z \pd{N_H}{m} \\
&\qquad + \pd{N_H}{m} ( z^{-1} B_H + R_H ) - \sum_{ k = 1 }^{ m - 1 } \pd{D_H}{k} \pd{N_H}{m} - \pd{D_H}{m} \sum_{ k = 1 }^m \pd{N_H}{k} \\
&= \partial_z \pd{N_H}{m} + \pd{N_H}{m} ( z^{-1} B_H + R_H ) - \sum_{ k = 1 }^{ m - 1 } \pd{D_H}{k} \pd{N_H}{m} - \pd{D_H}{m} \sum_{ k = 1 }^m \pd{N_H}{k} \text{.}
\end{align*}
Now, from Definition \ref{def.symbol} and \eqref{eq.perf_diag_h_DNR}, we have $ \partial_z \pd{N_H}{m} \in \mc{S}^{ -m ( \ell_\ast + 1 ) - 1 } ( Z_H; \C^n \otimes \C^n )$.
Thus, combining \eqref{eq.perf_diag_h_DNR} and \eqref{eql.perf_diag_h_1} with the above yields $\pd{R_H}{m} \in \mc{S}^{ -m ( \ell_\ast + 1 ) - 1 } ( Z_H; \C^n \otimes \C^n )$.
This completes the induction and establishes \eqref{eq.perf_diag_h_DNR}--\eqref{eq.perf_diag_h_RQ} for all $m \in \N$.

Finally, fix $m \in \N$.
By shrinking $\rho_0$, then each $| \pd{N_H}{k} |$, $1 \leq k \leq m$, can be made arbitrarily small on $Z_H$ (since $k ( \ell_\ast + 1 ) < -1$).
Thus, it follows from \eqref{eq.perf_diag_h_RQ} that $\pd{Q_H}{m}$ is invertible, with uniformly bounded inverse.
A direct computation using \eqref{eq.system_h} and \eqref{eq.perf_diag_h_RQ} then yields \eqref{eq.perf_diag_h}:
\begin{align*}
\partial_z \pd{U_H}{m} &= \sum_{ k = 1 }^m \partial_z \pd{N_H}{k} \, U_H + \pd{Q_H}{m} ( \imath \mf{Z} \, D_H + z^{-1} B_H + R_H ) U_H + \zfac^{-1} \pd{Q_H}{m} F_H \\
&= \left( \imath \mf{Z} \, D_H + \sum_{ k = 1 }^m \pd{D_H}{k} \right) \pd{U_H}{m} + \pd{R_H}{m} ( \pd{Q_H}{m} )^{-1} \pd{U_H}{m} + \zfac^{-1} \pd{F_H}{m} \text{.} \qedhere
\end{align*}
\end{proof}

\subsubsection{Estimates for Solutions}

We can now use \eqref{eq.perf_diag_h} to prove our main estimate on $Z_H$:

\begin{definition} \label{def.system_Eh}
Let $b_H: Z_H \rightarrow \C^n$ be defined as follows:
\begin{equation}
\label{eq.system_Eh_pre} b_{ H, i } ( t, \xi ) = \int_{ \zfac^{-1} \rho_0^{-1} }^t \tau^{-1} B_{ H, ii } ( \tau, \xi ) \, d \tau \text{,} \qquad 1 \leq i \leq n \text{.}
\end{equation}
We then define $\mc{E}_H: Z_H \rightarrow \C^n \otimes \C^n$ by
\begin{align}
\label{eq.system_Eh} \mc{E}_H :\!\!&= \operatorname{diag} \big( e^{ - b_{ H, 1 } }, \dots, e^{ - b_{ H, n } } \big) \text{.}
\end{align}
\end{definition}

\begin{proposition} \label{thm.energy_Eh}
$\mc{E}_H$ is $| \xi |$-independent, and there exist $C > 1$ and $c > 0$ such that
\begin{equation}
\label{eq.energy_Eh} C^{-1} z^{ -c } \leq | \mc{E}_{ H, ii } | \leq C z^c \text{,} \qquad 1 \leq i \leq n \text{.}
\end{equation}
\end{proposition}

\begin{proof}
That $\mc{E}_H$ is $| \xi |$-independent is immediate from \eqref{eq.system_Eh_pre}, since $B_{ H, ii }$ is $| \xi |$-independent.
For the estimate \eqref{eq.energy_Eh}, we use \eqref{eq.z} and \eqref{eq.z_coord} to rewrite \eqref{eq.system_Eh_pre} in $( z, \xi )$-coordinates as
\[
b_{ H, i } ( z, \xi ) = \int_{ \rho_0^{-1} }^z \zeta^{-1} B_{ H, ii } ( \zeta, \xi ) \, d \zeta \text{,} \qquad 1 \leq i \leq n \text{.}
\]
Since $B_H \in \mc{S}^0 ( Z_H; \C^n \otimes \C^n )$, then $-c \leq \operatorname{Re} B_{ H, ii } \leq c$ for some $c > 0$, so it follows that
\[
( \rho_0 z )^{ -c } \leq | e^{ - b_{ H, i } ( z, \xi ) } | \leq ( \rho_0 z )^c \text{.} \qedhere
\]
\end{proof}

\begin{proposition} \label{thm.energy_h}
Suppose \eqref{eq.system_gen} is semi-strictly hyperbolic, let $m \in \N$ be sufficiently large (depending on $B_H$ and $\ell_\ast$), and let $\rho_0$ be sufficiently small (with respect to $m$).
In addition, define
\begin{equation}
\label{eq.energy_h_U} \pd{U_{Hz}}{m} := \mc{E}_H \pd{Q_H}{m} U_H \text{,} \qquad \pd{F_{Hz}}{m} := \mc{E}_H \pd{Q_H}{m} F_H \text{.}
\end{equation}
Then, for any $\xi \in \R^d$ and $0 \leq t_0, t_1 \leq T$ such that $( t_0, \xi ), ( t_1, \xi ) \in Z_H$, we have
\begin{equation}
\label{eq.energy_h} | \pd{U_{Hz}}{m} ( t_1, \xi ) | \lesssim | \pd{U_{Hz}}{m} ( t_0, \xi ) | + \bigg| \int_{ t_0 }^{ t_1 } | \pd{F_{Hz}}{m} ( \tau, \xi ) | \, d \tau \bigg| \text{.}
\end{equation}
\end{proposition}

\begin{proof}
We adopt $( z, \xi )$-coordinates, and we set $z_0 := \zfac \, t_0$ and $z_1 := \zfac \, t_1$.
A direct computation using \eqref{eq.perf_diag_h} and Definition \ref{def.system_Eh} yields that $\pd{U_{Hz}}{m}$ satisfies
\begin{align*}
\partial_z \pd{U_{Hz}}{m} &= \mc{E}_H \bigg[ \imath \mf{Z} D_H + \partial_z \mc{E}_H + \sum_{ k = 1 }^m \pd{D_H}{k} + \pd{R_H}{m} ( \pd{Q_H}{m} )^{-1} \bigg] \mc{E}_H^{-1} \pd{U_{Hz}}{m} + \zfac^{-1} \pd{F_{Hz}}{m} \\
&= \bigg[ \imath \mf{Z} D_H + \partial_z \mc{E}_H + \sum_{ k = 1 }^m \pd{D_H}{k} \bigg] \pd{U_{Hz}}{m} + \mc{E}_H \pd{R_H}{m} ( \pd{Q_H}{m} )^{-1} \mc{E}_H^{-1} \pd{U_{Hz}}{m} + \zfac^{-1} \pd{F_{Hz}}{m} \text{,}
\end{align*}
where in the last step, we recalled \eqref{eq.system_Eh} and that $\mc{E}_H$, $D_H$, and $\pd{D_H}{1}, \dots, \pd{D_H}{m}$ are diagonal.
Moreover, expanding using \eqref{eq.perf_diag_h_DD} and Definition \ref{def.system_Eh}, we can then rewrite the above as
\begin{align}
\label{eql.energy_h_0} \partial_z \pd{U_{Hz}}{m} &= ( \imath \mf{Z} D_H + \pd{S_D}{m} + \pd{S_R}{m} ) \pd{U_{Hz}}{m} + \zfac^{-1} \pd{F_{Hz}}{m} \text{,} \\
\notag \pd{S_D}{m} :\!\!&= \operatorname{diag} ( R_{ H, 11 }, \dots, R_{ H, nn } ) + \sum_{ k = 2 }^m \pd{D_H}{k} \text{,} \\
\notag \pd{S_R}{m} :\!\!&= \mc{E}_H \pd{R_H}{m} ( \pd{Q_H}{m} )^{-1} \mc{E}_H^{-1} \text{.}
\end{align}

Multiplying \eqref{eql.energy_h_0} by $\pd{\bar{U}_{Hz}}{m}$, we obtain the estimate
\begin{equation}
\label{eql.energy_h_1} \tfrac{1}{2} \partial_z ( | \pd{U_{Hz}}{m} |^2 ) \leq ( | \pd{S_D}{m} | + | \pd{S_R}{m} | ) | \pd{U_{Hz}}{m} |^2 + \zfac^{-1} | \pd{F_{Hz}}{m} | | \pd{U_{Hz}}{m} | \text{,}
\end{equation}
where in the above, we used that
\[
\operatorname{Re} ( \imath \mf{Z} D_H \pd{U_{Hz}}{m} \cdot \pd{\bar{U}_{Hz}}{m} ) = 0
\]
due to $D_H$ being real-valued and $\imath \mf{Z} D_H$ hence being skew-Hermitian.

Now, from Assumption \ref{ass.system_h}, \eqref{eq.perf_diag_h_DNR}, and \eqref{eql.energy_h_0}, we have that
\[
\pd{S_D}{m} \in \mc{S}^{ a_H } ( Z_H; \C^n \otimes \C^n ) + \mc{S}^{ -1 - ( \ell_\ast + 1 ) } ( Z_H; \C^n \otimes \C^n ) \text{,}
\]
which implies the following uniform bound:
\begin{equation}
\label{eql.energy_h_10} \int_{ \rho_0^{-1} }^{ \zfac T } | \pd{S_D}{m} ( \zeta, \xi ) | \, d \zeta \lesssim 1 \text{.}
\end{equation}
Moreover, from Lemma \ref{thm.perf_diag_h} and \eqref{eq.energy_Eh}, we have for any $1 \leq i, j \leq n$ that
\[
| \pd{S_{ R, ij }}{m} ( z, \xi ) | \lesssim z^{ -m ( \ell_\ast + 1 ) - 1 + 2c } \text{.}
\]
Thus, as long as $m$ is sufficiently large, we also have
\begin{equation}
\label{eql.energy_h_11} \int_{ \rho_0^{-1} }^{ \zfac T } | \pd{S_R}{m} ( \zeta, \xi ) | \, d \zeta \lesssim 1 \text{.}
\end{equation}
Combining \eqref{eql.energy_h_1}--\eqref{eql.energy_h_11}, we obtain the following, which immediately implies \eqref{eq.energy_h}:
\[
| \pd{U_{Hz}}{m} ( z_1, \xi ) | \lesssim | \pd{U_{Hz}}{m} ( z_0, \xi ) | + \zfac^{-1} \bigg| \int_{ z_0 }^{ z_1 } | \pd{F_{Hz}}{m} ( \zeta, \xi ) | \, d \zeta \bigg| \text{.} \qedhere
\]
\end{proof}

\begin{remark}
From a closer inspection of its proof, we see that Proposition \ref{thm.energy_h} is applicable with any $m \in \N$ such that $[ \mc{E}_H \pd{R_H}{m} ( \pd{Q_H}{m} )^{-1} \mc{E}_H^{-1} ] ( \cdot, \xi )$ is integrable in $t$ for all large $\xi \in \R^d$.
\end{remark}

In many cases, one has a simpler characterization of the size of $\pd{U_{Hz}}{m}$:

\begin{proposition} \label{thm.energy_h_simple}
Assume the setting of Proposition \ref{thm.energy_h}, and suppose in addition there exists a function $B_{ H, 0 }: \Sph^{ d - 1 } \rightarrow \C^n$ satisfying the following bound for some $\delta > 0$:
\begin{equation}
\label{eq.energy_h_simple_ass} | B_{ H, \infty, ii } ( t, \omega ) - B_{ H, 0, i } ( \omega ) | \lesssim t^\delta \text{,} \qquad \omega \in \Sph^{ d - 1 } \text{,} \quad 1 \leq i \leq n \text{.}
\end{equation}
Then, the following comparison holds everywhere on $Z_H$,
\begin{equation}
\label{eq.energy_h_simple} | \pd{U_{Hz}}{m} ( t, \xi ) | \simeq | ( \mc{E}_{ H\ast } \pd{Q_H}{m} U_H ) ( t, \xi ) | \text{,} \qquad | \pd{F_{Hz}}{m} ( t, \xi ) | \simeq | ( \mc{E}_{ H\ast } \pd{Q_H}{m} F_H ) ( t, \xi ) | \text{,}
\end{equation}
where the constants depend on $m$, $\rho_0$, and $B_H$, and where $\mc{E}_{ H\ast }$ is defined on $Z_H$ by
\begin{equation}
\label{eq.energy_h_Ehh} \mc{E}_{ H\ast } ( t, \xi ) := \operatorname{diag} \big( z^{ - B_{ H, 0, 1 } ( \frac{ \xi }{ | \xi | } ) }, \dots, z^{ - B_{ H, 0, n } ( \frac{ \xi }{ | \xi | } ) } \big) \text{.}
\end{equation}
\end{proposition}

\begin{proof}
Fix $1 \leq i \leq n$.
Recalling \eqref{eq.z}, Assumption \ref{ass.system_h}, and \eqref{eq.system_Eh_pre}, we can write
\begin{align}
\label{eql.energy_h_simple_1} b_{ H, i } ( t, \xi ) &= B_{ H, 0, i } \big( \tfrac{ \xi }{ | \xi | } \big) \int_{ \zfac^{-1} \rho_0^{-1} }^t \tau^{-1} d \tau + \int_{ \zfac^{-1} \rho_0^{-1} }^t \tau^{-1} \big[ B_{ H, \infty, ii } \big( \tau, \tfrac{ \xi }{ | \xi | } \big) - B_{ H, 0, i } \big( \tfrac{ \xi }{ | \xi | } \big) \big] d \tau \\
\notag &= B_{ H, 0, i } \big( \tfrac{ \xi }{ | \xi | } \big) \, \log ( \rho_0 z ) + \int_{ \zfac^{-1} \rho_0^{-1} }^t \tau^{-1} \big[ B_{ H, \infty, ii } \big( \tau, \tfrac{ \xi }{ | \xi | } \big) - B_{ H, 0, i } \big( \tfrac{ \xi }{ | \xi | } \big) \big] \, d \tau \text{.}
\end{align}
Moreover, note from \eqref{eq.energy_h_simple_ass} that
\begin{align*}
\bigg| \int_{ \zfac^{-1} \rho_0^{-1} }^t \tau^{-1} \big[ B_{ H, \infty, ii } \big( \tau, \tfrac{ \xi }{ | \xi | } \big) - B_{ H, 0, i } \big( \tfrac{ \xi }{ | \xi | } \big) \big] \, d \tau \bigg| &\lesssim \int_0^T \tau^{ -1 + \delta } \, d \tau \\
&\lesssim 1 \text{.}
\end{align*}
Combining \eqref{eql.energy_h_simple_1} and the above, and recalling that $B_H \in \mc{S}^0 ( Z_H; \C^n \otimes \C^n )$, we conclude
\[
e^{ - b_{ H, i } ( t, \xi ) } \simeq z^{ - B_{ H, 0, i } ( \frac{ \xi }{ | \xi | } ) } \text{.}
\]
The desired \eqref{eq.energy_h_simple} now follows from \eqref{eq.system_Eh}, \eqref{eq.energy_h_U}, \eqref{eq.energy_h_Ehh}, and the above.
\end{proof}

\begin{remark}
When \eqref{eq.energy_h_simple_ass} holds, one could alternatively define
\[
\pd{U_{Hz}}{m} := \mc{E}_{ H\ast } \pd{Q_H}{m} U_H \text{,} \qquad \pd{F_{Hz}}{m} \simeq \mc{E}_{ H\ast } \pd{Q_H}{m} F_H \text{.}
\]
However, we elect not to adopt the above definition, since there are physically interesting settings (e.g.\ wave equations on Kasner spacetimes) for which \eqref{eq.energy_h_simple_ass} fails to hold.
\end{remark}

\subsection{Scattering and Asymptotics} \label{sec.system_results}

In this subsection, we combine the individual estimates obtained on $Z_P$, $Z_I$, and $Z_H$ to derive the main result of this paper---global estimates for (the proper renormalization of) $U$, along with scattering and asymptotics results for \eqref{eq.system_gen} at $t \searrow 0$.

\subsubsection{The Main Results}

The precise statements of our main scattering and asymptotics (and by extension, loss of regularity) results are summarized in the subsequent theorem:

\begin{theorem} \label{thm.energy_main}
Suppose Assumptions \ref{ass.system}, \ref{ass.system_i}, \ref{ass.system_p}, and \ref{ass.system_h} hold, and suppose \eqref{eq.system_gen} is semi-strictly hyperbolic.
Also, let $m_P \in \N_0$ and $m_H \in \N$ be sufficiently large (with $m_P$ depending on $a_P$, $B_P$; and $m_H$ depending on $\ell_\ast$, $B_H$), and let $\rho_0$ be sufficiently small (with respect to $m_P$, $m_H$).
\begin{itemize}
\item \emph{Asymptotics:} Given $u: \R^d_\xi \rightarrow \C^n$ and $F$ as in Assumption \ref{ass.system}, there exists a unique solution $U$ of \eqref{eq.system_gen} satisfying $U ( T, \cdot ) = u$.
Moreover, if we define $U_A, F_A: ( 0, T ]_t \times \R_\xi^d \rightarrow \C^n$ by
\begin{align}
\label{eq.system_UA} U_A |_{ Z_P \setminus Z_I } := \pd{U_{Pz}}{m_P} \text{,} \qquad U_A |_{ Z_I } &:= U \text{,} \qquad U_A |_{ Z_H \setminus Z_I } := \pd{U_{Hz}}{m_H} \text{,} \\
\notag F_A |_{ Z_P \setminus Z_I } := \pd{F_{Pz}}{m_P} \text{,} \qquad F_A |_{ Z_I } &:= F \text{,} \qquad F_A |_{ Z_H \setminus Z_I } := \pd{F_{Hz}}{m_H} \text{,}
\end{align}
and if we assume that $F_A$ is $t$-integrable,
\begin{equation}
\label{eq.energy_main_ass} \int_0^T | F_A ( \tau, \xi ) | \, d \tau < \infty \text{,} \qquad \xi \in \R^d \text{,}
\end{equation}
then the following asymptotic limits are both well-defined and finite:
\begin{equation}
\label{eq.energy_asymp_limit} u_A ( \xi ) := \lim_{ \tau \searrow 0 } U_A ( \tau, \xi ) \text{,} \qquad \xi \in \R^d \text{.}
\end{equation}
In addition, the following estimate holds for any $\xi \in \R$ and $t_0 \in ( 0, T ]$:
\begin{equation}
\label{eq.energy_asymp} | u_A ( \xi ) | \lesssim | U_A ( t_0, \xi ) | + \int_0^{ t_0 } | F_A ( \tau, \xi ) | \, d \tau \text{.}
\end{equation}

\item \emph{Scattering:} Given $u_A: \R^d_\xi \rightarrow \C^n$ and $F$ as in Assumption \ref{ass.system}, if $F_A$ is defined as in \eqref{eq.system_UA}, and if \eqref{eq.energy_main_ass} holds (i.e.\ $F_A$ is $t$-integrable), then there is a unique solution $U$ of \eqref{eq.system_gen} such that, defining $U_A$ as in \eqref{eq.system_UA}, the following asymptotic limits hold:
\begin{equation}
\label{eq.energy_scatter_limit} \lim_{ \tau \searrow 0 } U_A ( \tau, \xi ) = u_A ( \xi ) \text{,} \qquad \xi \in \R^d \text{.}
\end{equation}
In addition, the following estimate holds for any $\xi \in \R^d$ and $t_0 \in ( 0, T ]$:
\begin{equation}
\label{eq.energy_scatter} | U_A ( t_0, \xi ) | \lesssim | u_A ( \xi ) | + \int_0^{ t_0 } | F_A ( \tau, \xi ) | \, d \tau \text{.}
\end{equation}
\end{itemize}
\end{theorem}

\begin{proof}
The first step is to solve for $U$.
For the asymptotics case, one can simply solve \eqref{eq.system_gen} (or more accurately, \eqref{eq.system_i}, \eqref{eq.system_p}, and \eqref{eq.system_h}) directly as a linear system of differential equations at every $\xi \in \R^d$.
Moreover, if \eqref{eq.energy_main_ass} holds, then Proposition \ref{thm.energy_p} yields the limits \eqref{eq.energy_asymp_limit}.

Conversely, for the scattering case, one first solves for $\pd{U_{Pz}}{m_P}$ via Proposition \ref{thm.energy_p}.
(More specifically, as \eqref{eq.energy_main_ass} holds, one can solve \eqref{eql.energy_p_0} for $\pd{U_{Pz}}{m_P}$ with initial data $u_A$ at $t \searrow 0$; see the remarks after Proposition \ref{thm.energy_p}.)
Reversing all the transformations yields $U$ on $Z_P$.
Away from $Z_P$, one can then solve \eqref{eq.system_gen} directly for $U$ given data at a point within $Z_P$.

As a result, it remains only to prove the estimates \eqref{eq.energy_asymp} and \eqref{eq.energy_scatter}.
For convenience, we set
\begin{equation}
\label{eql.energy_scatter_0} t_P := \rho_0 \zfac^{-1} \text{,} \qquad t_H := \rho_0^{-1} \zfac^{-1} \text{.}
\end{equation}
Note in particular $( t_P, \xi ) \in Z_P \cap Z_I$ (provided $t_P \leq T$), while $( t_H, \xi ) \in Z_I \cap Z_H$ (if $t_H \leq T$).

First, if $( t_0, \xi ) \in Z_P \setminus Z_I$, then Proposition \ref{thm.energy_p} and \eqref{eq.system_UA} yield
\begin{align*}
| U_A ( t_0, \xi ) | &= | \pd{U_{Pz}}{m_P} ( t_0, \xi ) | \\
&\lesssim | \pd{U_{Pz}}{m_P} ( 0, \xi ) | + \int_0^{ t_0 } | \pd{F_{Pz}}{m_P} ( \tau, \xi ) | \, d \tau \\
&= | u_A ( \xi ) | + \int_0^{ t_0 } | F_A ( \tau, \xi ) | \, d \tau \text{,}
\end{align*}
proving \eqref{eq.energy_scatter} in this case.
Similarly, by Proposition \ref{thm.energy_p} and \eqref{eq.system_UA},
\begin{align*}
| u_A ( \xi ) | &= | \pd{U_{Pz}}{m_P} ( 0, \xi ) | \\
&\lesssim | U_A ( t_0, \xi ) | + \int_0^{ t_0 } | F_A ( \tau, \xi ) | \, d \tau \text{,}
\end{align*}
which proves \eqref{eq.energy_asymp} in this particular case.

Next, suppose $( t_0, \xi ) \in Z_I$.
Then, by Assumption \ref{ass.system_p}, \eqref{eq.system_p}, Proposition \ref{thm.energy_i}, and \eqref{eq.system_UA},
\begin{align*}
| U_A ( t_0, \xi ) | &= | U ( t_0, \xi ) | \\
&\lesssim | U_P ( t_P, \xi ) | + \int_{ t_P }^{ t_0 } | F_P ( \tau, \xi ) | \, d \tau \\
&\lesssim | \pd{U_{Pz}}{m_P} ( t_P, \xi ) | + \int_{ t_P }^{ t_0 } | F_P ( \tau, \xi ) | \, d \tau \text{,}
\end{align*}
where in the last step, we used that $z \simeq 1$ on $Z_I$ (see \eqref{eq.zones}) and that $( \pd{Q_P}{m} )^{-1}$ is bounded on $Z_P$ (see Lemma \ref{thm.perf_diag_p}).
Proceeding now as in the previous case, we estimate
\begin{align*}
| U_A ( t_0, \xi ) | &\lesssim | \pd{U_{Pz}}{m_P} ( 0, \xi ) | + \int_{ t_P }^{ t_0 } | F_P ( \tau, \xi ) | \, d \tau + \int_0^{ t_P } | \pd{F_{Pz}}{m_P} ( \tau, \xi ) | \, d \tau \\
&= | u_A ( \xi ) | + \int_0^{ t_0 } | F_A ( \tau, \xi ) | \, d \tau \text{,}
\end{align*}
which is \eqref{eq.energy_scatter} in this setting.
Again, all the above steps can be reversed:
\begin{align*}
| u_A ( \xi ) | &\lesssim | \pd{U_{Pz}}{m_P} ( t_P, \xi ) | + \int_0^{ t_P } | \pd{F_{Pz}}{m_P} ( \tau, \xi ) | \, d \tau \\
&\lesssim | U ( t_0, \xi ) | + \int_0^{ t_P } | \pd{F_{Pz}}{m_P} ( \tau, \xi ) | \, d \tau + \int_{ t_P }^{ t_0 } | F_P ( \tau, \xi ) | \, d \tau \\
&= | U_A ( t_0, \xi ) | + \int_0^{ t_0 } | F_A ( \tau, \xi ) | \, d \tau \text{.}
\end{align*}
This establishes \eqref{eq.energy_asymp} for this setting.

Finally, suppose $( t_0, \xi ) \in Z_H \setminus Z_I$.
Then, by Assumption \ref{ass.system_h}, Proposition \ref{thm.energy_h}, and \eqref{eq.system_UA},
\begin{align*}
| U_A ( t_0, \xi ) | &= | \pd{U_{Hz}}{m_H} ( t_0, \xi ) | \\
&\lesssim | \pd{U_{Hz}}{m_H} ( t_H, \xi ) | + \int_{ t_H }^{ t_0 } | \pd{F_{Hz}}{m_H} ( \tau, \xi ) | \, d \tau \\
&\lesssim | U ( t_H, \xi ) | + \int_{ t_H }^{ t_0 } | \pd{F_{Hz}}{m_H} ( \tau, \xi ) | \, d \tau \text{,}
\end{align*}
where in the last step, we used (recalling Definition \ref{def.system_Eh}) that $\mc{E}_H ( t_H, \xi ) = I_n$, that $z \simeq 1$ on $Z_I$ (from \eqref{eq.zones}), and that $\pd{Q_H}{m}$ is bounded on $Z_H$ (see Lemma \ref{thm.perf_diag_h}).
From the above, we proceed once again as in the previous case to continue the estimate to $t \searrow 0$,
\begin{align*}
| U_A ( t_0, \xi ) | &\lesssim | u_A ( \xi ) | + \int_0^{ t_H } | F_A ( \tau, \xi ) | \, d \tau + \int_{ t_H }^{ t_0 } | \pd{F_{Hz}}{m_H} ( \tau, \xi ) | \, d \tau \\
&= | u_A ( \xi ) | + \int_0^{ t_0 } | F_A ( \tau, \xi ) | \, d \tau \text{,}
\end{align*}
resulting in \eqref{eq.energy_scatter}.
Finally, we can again reverse the above steps to obtain \eqref{eq.energy_asymp}:
\begin{align*}
| u_A ( \xi ) | &\lesssim | \pd{U_{Hz}}{m_H} ( t_H, \xi ) | + \int_0^{ t_H } | F_A ( \tau, \xi ) | \, d \tau \\
&\lesssim | \pd{U_{Hz}}{m_H} ( t_0, \xi ) | + \int_{ t_H }^{ t_0 } | \pd{F_{Hz}}{m_H} ( \tau, \xi ) | \, d \tau + \int_0^{ t_H } | F_A ( \tau, \xi ) | \, d \tau \\
&= | U_A ( t_0, \xi ) | + \int_0^{ t_0 } | F_A ( \tau, \xi ) | \, d \tau \text{.} \qedhere
\end{align*}
\end{proof}

\begin{remark}
Note $U_A$ and $F_A$, as defined in \eqref{eq.system_UA}, fails to be continuous at the boundaries between $Z_P$ and $Z_I$ and between $Z_I$ and $Z_H$.
One could, alternatively, more smoothly connect the values of $U_A$ on $Z_P$, $Z_I$, $Z_H$ using appropriate cutoff functions on $( 0, T ]_t \times \R^d_\xi$.
In particular, these can be chosen such that the ``smoothed" unknown $U_A$ still satisfies \eqref{eq.energy_asymp} and \eqref{eq.energy_scatter}.
\end{remark}

\subsubsection{Sobolev Convergence}

While we have demonstrated the existence of a pointwise limit $u_A$ as $t \searrow 0$ for the renormalized quantity $U_A$ at each frequency $\xi$, what is more relevant, in terms of the analysis of the original system \eqref{eq.system_actual}, is the convergence in terms of Sobolev spaces.

In the Sobolev sense, there is in fact a slight loss of derivatives in the convergence.
This arises from the fact that the uniformity of our estimates for the system \eqref{eq.system_gen} is in terms of the rescaled $z$, which is frequency-dependent.
We quantify this loss in the following:

\begin{proposition} \label{thm.energy_conv}
Assume the setting of Theorem \ref{thm.energy_main}, and suppose
\begin{equation}
\label{eq.energy_conv_ass} \int_{ \R^d } \langle \xi \rangle^{ 2s } | u_A ( \xi ) |^2 \, d \xi < \infty \text{,} \qquad \int_{ Z_P } \langle \xi \rangle^{ 2s } | \pd{F_A}{ m_P } ( \tau, \xi ) |^2 \, d \xi d \tau < \infty \text{.}
\end{equation}
Then, for any $\delta > 0$, we have that
\begin{equation}
\label{eq.energy_conv} \lim_{ t \searrow 0 } \int_{ ( \{ t \} \times \R^d ) \cap Z_P } \langle \xi \rangle^{ 2 ( s - \delta ) } | U_A ( \tau, \xi ) - u_A ( \xi ) |^2 \, d \xi = 0 \text{.}
\end{equation}
\end{proposition}

\begin{proof}
First, by the equation \eqref{eql.energy_p_0} satisfied by $\pd{U_{Pz}}{m_P}$, \eqref{eq.system_UA}, and \eqref{eq.energy_scatter}, we have
\begin{align*}
| U_A ( t, \xi ) - u_A ( \xi ) | &\lesssim \int_0^{ \zfac t } | \partial_z \pd{U_{Pz}}{m_P} ( \zeta, \xi ) | \, d \zeta \\
&\lesssim \int_0^{ \zfac t } | \pd{S}{m_P} ( \zeta, \xi ) | | \pd{U_{Pz}}{m_P} ( \zeta, \xi ) | \, d \zeta + \zfac^{-1} \int_0^{ \zfac t } | \pd{F_{Pz}}{m_P} ( \zeta, \xi ) | \, d \zeta \\
&\lesssim | u_A ( \xi ) | \int_0^{ \zfac t } | \pd{S}{m_P} ( \zeta, \xi ) | \, d \zeta + \int_0^t | \pd{F_{Pz}}{m_P} ( \tau, \xi ) | \, d \tau \text{,}
\end{align*}
where $\pd{S}{m_P}$ was defined in \eqref{eql.energy_p_0} and is uniformly $z$-integrable.
Shrinking $\delta$ if needed and recalling also \eqref{eq.z_factor} and the estimates in Lemma \ref{thm.energy_Ep} for the terms of $\pd{S}{m_P}$, the above then implies
\begin{align*}
| U_A ( t, \xi ) - u_A ( \xi ) | &\lesssim | u_A ( \xi ) | \int_0^{ \zfac t } \zeta^{ -1 + ( \ell_\ast + 1 ) \delta } \, d \zeta + \int_0^t | \pd{F_{Pz}}{m_P} ( \tau, \xi ) | \, d \tau \\
&\lesssim \langle \xi \rangle^\delta | u_A ( \xi ) | t^{ ( \ell_\ast + 1 ) \delta } + \int_0^t | \pd{F_{Pz}}{m_P} ( \tau, \xi ) | \, d \tau \text{,}
\end{align*}
for any $( t, \xi ) \in Z_P$.
Squaring the above and integrating over $\xi$, we then have
\begin{align*}
\int_{ ( \{ t \} \times \R^d ) \cap Z_P } \langle \xi \rangle^{ 2 ( s - \delta ) } | U_A ( \tau, \xi ) - u_A ( \xi ) |^2 \, d \xi &\lesssim t^{ 2 ( \ell_\ast + 1 ) \delta } \int_{ \R^d } \langle \xi \rangle^{ 2 s } | u_A ( \xi ) | \, d \xi \\
&\qquad + \int_{ ( ( 0, t ] \times \R^d ) \cap Z_P } \langle \xi \rangle^{ 2 s } | \pd{F_{Pz}}{m_P} ( \tau, \xi ) |^2 \, d \xi d \tau \text{.}
\end{align*}
By \eqref{eq.energy_conv_ass}, the right-hand side of the above converges to $0$ as $t \searrow 0$, and \eqref{eq.energy_conv} follows.
\end{proof}

\begin{remark}
Proposition \ref{thm.energy_conv} shows a loss of $\delta$ derivatives, for any $\delta > 0$, in Sobolev convergence.
\end{remark}

\begin{remark}
Notice the limit in \eqref{eq.energy_conv} is necessarily restricted to $Z_P$, since the renormalization $\pd{U_{Pz}}{m_P}$ is only generally defined on $Z_P$, and since $U_A$ is discontinuous at the zonal boundaries $Z_P \cap Z_I$ and $Z_I \cap Z_H$.
However, one could recover a more global convergence statement on $[ 0, T ] \times \R^d$, provided the values of $U_A$ are defined to be more smoothly connected between $Z_P$, $Z_I$, and $Z_H$.
\end{remark}

\subsection{Estimates Without Semi-Strict Hyperbolicity}

In this final subsection, we consider the more general case of systems that fail to be semi-strictly hyperbolic.
In this case, one still obtains estimates, albeit with some loss, in that these do not yield converse asymptotics and scattering properties as in Theorem \ref{thm.energy_main}.
Below, we sketch the argument for deriving these estimates.

In the following, we again suppose Assumptions \ref{ass.system}, \ref{ass.system_i}, \ref{ass.system_p}, and \ref{ass.system_h} hold.
Moreover, since the analyses on $Z_P$ and $Z_I$ are identical to before, we hence now focus exclusively on $Z_H$.

\subsubsection{Partial Renormalization}

Since semi-strict hyperbolicity no longer holds by assumption, the first step is to describe which components of $D_H$ have separated speeds:

\begin{definition} \label{def.system_h_partition}
Let $\mc{G}$ be a partition of $\{ 1, \dots, n \}$.
Then:
\begin{itemize}
\item $\mc{G}$ is called a \emph{partition of speeds} of $D_H$ iff there exists $d_0 > 0$ such that for any $1 \leq i, j \leq n$ in different elements of $\mc{G}$, we have the following inequality uniformly on $Z_H$:
\begin{equation}
\label{eq.system_h_partition} | D_{ H, ii } - D_{ H, jj } | \geq d_0 \text{.}
\end{equation}

\item A matrix $W \in \C^n \otimes \C^n$ (or analogously, a matrix-valued function) is called \emph{$\mc{G}$-diagonal} iff $W_{ij} = 0$ for any $1 \leq i, j \leq n$ that lie in different elements of $\mc{G}$.
\end{itemize}
\end{definition}

\begin{remark}
In particular, if \eqref{eq.system_gen} is \emph{semi-strictly hyperbolic}, then
\begin{equation}
\label{eq.system_h_partition_full} \mc{G}_{\text{sh}} := \{ \{ 1 \}, \dots, \{ n \} \}
\end{equation}
is a partition of speeds of $D_H$.
Note $W \in \C^n \otimes \C^n$ is $\mc{G}_{\text{sh}}$-diagonal if and only if $W$ is diagonal.
\end{remark}

We now replace our previous assumption of semi-strict hyperbolicity with the following:

\begin{assumption} \label{ass.partition}
Let $\mc{G}$ be a partition of speeds of $D_H$.
\end{assumption}

The following lemma provides the analogue of Lemma \ref{thm.perf_diag_h} in settings without semi-strict hyperbolicity.
As a slight abuse of notation, we use the same symbols here as in Section \ref{sec.system_zh}.
Note that because of non-separated speeds, we can only $\mc{G}$-diagonalize our system to higher orders.

\begin{lemma} \label{thm.part_diag_h}
There exists a sequence $( \pd{D_H}{m}, \pd{N_H}{m}, \pd{R_H}{m} )_{ m \in \N }$, satisfying
\begin{align}
\label{eq.part_diag_h_DNR} \pd{D_H}{m} &\in \mc{S}^{ - ( m - 1 ) ( \ell_\ast + 1 ) - 1 } ( Z_H; \C^n \otimes \C^n ) \text{,} \\
\notag \pd{N_H}{m} &\in \mc{S}^{ -m ( \ell_\ast + 1 ) } ( Z_H; \C^n \otimes \C^n ) \text{,} \\
\notag \pd{R_H}{m} &\in \mc{S}^{ -m ( \ell_\ast + 1 ) - 1 } ( Z_H; \C^n \otimes \C^n )
\end{align}
for all $m \in \N$, such that the following properties hold on $Z_H$:
\begin{itemize}
\item $\pd{D_{ H, ij }}{m}$ is $\mc{G}$-diagonal for any $m \in \N$.
In particular,
\begin{equation}
\label{eq.part_diag_h_DD} \pd{D_{ H, ij }}{1} = \begin{cases} z^{-1} B_{ H, ij } + R_{ H, ij } & \text{$i, j$ lie in the same element of $\mc{G}$,} \\ 0 & \text{$i, j$ lie in different elements of $\mc{G}$.} \end{cases}
\end{equation}

\item Each $\pd{R_H}{m}$, $m \in \N$, is given by the following:
\begin{align}
\label{eq.part_diag_h_RQ} \pd{R_H}{m} &= \sum_{ k = 1 }^m \partial_z \pd{N_H}{k} + \pd{Q_H}{m} ( \imath \mf{Z} \, D_H + z^{-1} B_H + R_H ) - \bigg( \imath \mf{Z} \, D_H + \sum_{ k = 1 }^m \pd{D_H}{k} \bigg) \pd{Q_H}{m} \text{,} \\
\notag \pd{Q_H}{m} :\!\!&= I_n + \sum_{ k = 1 }^m \pd{N_H}{k} \in \mc{S}^0 ( Z_H; \C^n \otimes \C^n ) \text{.}
\end{align}
\end{itemize}
Furthermore, given $m \in \N$, if $\rho_0$ is sufficiently small (with respect to $m$), then:
\begin{itemize}
\item $\pd{Q_H}{m}$ is invertible on $Z_H$, and both $\pd{Q_H}{m}$, $( \pd{Q_H}{m} )^{-1}$ are uniformly bounded on $Z_H$.

\item The following system of differential equations holds on $Z_H$:
\begin{align}
\label{eq.part_diag_h} \partial_z \pd{U_H}{m} &= \bigg( \imath \mf{Z} \, D_H + \sum_{ k = 1 }^m \pd{D_H}{k} \bigg) \pd{U_H}{m} + \pd{R_H}{m} ( \pd{Q_H}{m} )^{-1} \pd{U_H}{m} + \zfac^{-1} \pd{Q_H}{m} F_H \text{,} \\
\notag ( \pd{U_H}{m}, \pd{F_H}{m} ) :\!\!&= ( \pd{Q_H}{m} U_H, \pd{Q_H}{m} F_H ) \text{.}
\end{align}
\end{itemize}
\end{lemma}

\begin{proof}
We again begin by setting
\begin{equation}
\label{eql.part_diag_h_0} \pd{R_H}{0} := z^{-1} B_H + R_H \in \mc{S}^{-1} ( Z_H; \C^n \otimes \C^n ) \text{.}
\end{equation}
Fix $m \in \N$, and suppose we have defined $( \pd{D_H}{k}, \pd{N_H}{k} )$ for all $1 \leq k < m$ and $\pd{R_H}{k}$ for all $0 \leq k < m$ satisfying \eqref{eq.part_diag_h_DNR} and \eqref{eq.part_diag_h_RQ}, and with $\pd{D_H}{k}$ $\mc{G}$-diagonal for any $1 \leq k < m$.

The key difference with Lemma \ref{thm.perf_diag_h} is we now define $\pd{D_H}{m} \in \mc{S}^{ - ( m - 1 ) ( \ell_\ast + 1 ) - 1 } ( Z_H; \C^n \otimes \C^n )$ and $\pd{N_H}{m} \in \mc{S}^{ -m ( \ell_\ast + 1 ) } ( Z_H; \C^n \otimes \C^n )$ by the following formulas for all $1 \leq i, j \leq n$:
\begin{align}
\label{eql.part_diag_h_1} \pd{D_{H, ij}}{m} :\!\!&= \begin{cases} \pd{R_{H, ij}}{m-1} & \text{$i, j$ lie in the same element of $\mc{G}$,} \\ 0 & \text{$i, j$ lie in different elements of $\mc{G}$,} \end{cases} \\
\notag \pd{N_{H, ij}}{m} :\!\!&= \begin{cases} 0 & \text{$i, j$ lie in the same element of $\mc{G}$,} \\ ( \imath \mf{Z} \, ( D_{ H, ii } - D_{ H, jj } ) )^{-1} \pd{R_{H, ij}}{m-1} & \text{$i, j$ lie in different elements of $\mc{G}$.} \end{cases}
\end{align}
(Note $\pd{D_H}{m}$ is by definition $\mc{G}$-diagonal.)
Observe that \eqref{eql.part_diag_h_1} again implies
\begin{equation}
\label{eql.part_diag_h_3} \pd{R_H}{m-1} - \pd{D_H}{m} + [ \pd{N_H}{m}, \imath \mf{Z} \, D_H ] = 0 \text{.}
\end{equation}

The remainder of the proof is identical to that of Lemma \ref{thm.perf_diag_h}.
We again define $\pd{Q_H}{m}$ and $\pd{R_H}{m}$ via \eqref{eq.part_diag_h_RQ}.
The same computations as in Lemma \ref{thm.perf_diag_h} yield $\pd{R_H}{m} \in \mc{S}^{ -m ( \ell_\ast + 1 ) - 1 } ( \C^n \otimes \C^n )$, completing the induction and proving \eqref{eq.perf_diag_h_DNR}--\eqref{eq.perf_diag_h_RQ} for all $m \in \N$.

Furthermore, given any $m \in \N$, choosing $\rho_0 \ll_m 1$ makes each $| \pd{N_H}{k} |$, $1 \leq k \leq m$, arbitrarily small, yielding the invertibility and boundedness of $\pd{Q_H}{m}$.
Finally, that \eqref{eq.part_diag_h} holds is a consequence of the same computation as in the proof of Lemma \ref{thm.perf_diag_h}.
\end{proof}

\subsubsection{Estimates on $Z_H$}

Next, we prove our key estimates on $Z_H$, which will be weaker versions of those proved in Proposition \ref{thm.energy_h} under semi-strict hyperbolicity.

\begin{definition} \label{def.part_diag_decomp}
Let $W \in \C^n \otimes \C^n$ be $\mc{G}$-diagonal.
\begin{itemize}
\item Let $W_s \in \C^n \otimes \C^n$ denote the Hermitian part of $W$:
\begin{equation}
\label{eq.part_diag_decomp_herm} W_s := \tfrac{1}{2} ( W + W^\ast ) \text{.}
\end{equation}

\item Given $\sigma \in \mc{G}$, we let $W_s [ \sigma ]$ denote the corresponding $\sigma$-indexed matrix,
\begin{equation}
\label{eq.part_diag_decomp_sigma} W_s [ \sigma ] := [ W_{ s, ij } ]_{ i, j \in \sigma } \text{.}
\end{equation}

\item Given $\sigma \in \mc{G}$, we let $W_+ [ \sigma ]$ and $W_- [ \sigma ]$ denote the largest and smallest eigenvalues of $W_s [ \sigma ]$, respectively.
(Note $W_s [ \sigma ]$ is Hermitian and hence has real eigenvalues.)

\item Let $W_+, W_- \in \C^n \otimes \C^n$ denote the following diagonal matrices:
\begin{equation}
\label{eq.part_diag_decomp} W_{ \pm, ij } := \begin{cases} W_\pm [ \sigma ] & i = j \in \sigma \text{,} \\ 0 & \text{otherwise,} \end{cases} \qquad 1 \leq i, j \leq n \text{,} \quad \sigma \in \mc{G} \text{.}
\end{equation}
\end{itemize}
The above constructions can be analogously applied to $\mc{G}$-diagonal matrix-valued functions.
\end{definition}

\begin{remark}
When \eqref{eq.system_gen} is semi-strictly hyperbolic, i.e.\ $\mc{G} := \mc{G}_{\text{sh}}$, we have, with $W$ as above,
\[
W_+ = W_- = \operatorname{Re} W \text{.}
\]
\end{remark}

\begin{lemma} \label{thm.part_diag_est}
Let $W \in \C^n \otimes \C^n$ be $\mc{G}$-diagonal.
Then, for any $X \in \C^n$,
\begin{equation}
\label{eq.part_diag_est} W_- X \cdot \bar{X} \leq W X \cdot \bar{X} \leq W_+ X \cdot \bar{X} \text{.}
\end{equation}
\end{lemma}

\begin{proof}
First, since $W_s$ is the Hermitian part of $X$, then
\begin{equation}
\label{eql.part_diag_est_0} W X \cdot \bar{X} = W_s X \cdot \bar{X} \text{.}
\end{equation}
Now, since $W_s - W_+$ is negative-semidefinite, and since $W_s - W_-$ is positive-semidefinite, then
\[
W_- X \cdot \bar{X} \leq W_s X \cdot \bar{X} \leq W_+ X \cdot \bar{X} \text{,}
\]
and the desired \eqref{eq.part_diag_est} now follows from \eqref{eql.part_diag_est_0} and the above.
\end{proof}

Next, for convenience, we define the following quantities, which will be crucial for constructing our higher-order renormalizations and quantifying the ensuing loss of regularity:

\begin{definition} \label{def.system_Ph}
We define the following $\mc{G}$-diagonal matrices on $Z_H$---for $1 \leq i, j \leq n$:
\begin{align}
\label{eq.system_Ph_BR} B_{ \mc{G}, ij } &= \begin{cases} B_{ H, ij } & \text{$i, j$ lie in the same element of $\mc{G}$,} \\ 0 & \text{$i, j$ lie in different elements of $\mc{G}$,} \end{cases} \\
\notag R_{ \mc{G}, ij } &= \begin{cases} R_{ H, ij } & \text{$i, j$ lie in the same element of $\mc{G}$,} \\ 0 & \text{$i, j$ lie in different elements of $\mc{G}$.} \end{cases}
\end{align}
We now define the vector-valued functions $b_H^\pm: Z_H \rightarrow \C^n$ as follows:
\begin{equation}
\label{eq.system_Ph_pre} b^\pm_{ H, i } ( t, \xi ) = \int_{ \zfac^{-1} \rho_0 }^t \tau^{-1} B_{ \mc{G}, \pm, ii } ( \tau, \xi ) \, d \tau \text{,} \qquad 1 \leq i \leq n \text{.}
\end{equation}
We then define $\mc{E}_H^\pm: Z_H \rightarrow \C^n \otimes \C^n$ by
\begin{align}
\label{eq.system_Ph} \mc{E}_H^\pm :\!\!&= \operatorname{diag} \big( e^{ - b^\pm_{ H, 1 } }, \dots, e^{ - b^\pm_{ H, n } } \big) \text{.}
\end{align}
\end{definition}

We can now state and prove our main estimate on $Z_H$, i.e.\ the analogue of Proposition \ref{thm.energy_h}:

\begin{proposition} \label{thm.energy_Ph}
Let $m \in \N$ be sufficiently large (depending on $B_H$, $\ell_\ast$), and let $\rho_0$ be sufficiently small (with respect to $m$).
In addition, we define the following quantities:
\begin{equation}
\label{eq.energy_Ph_U} \pd{U_{Hz, \pm}}{m} := \mc{E}_H^\pm \pd{Q_H}{m} U_H \text{,} \qquad \pd{F_{Hz, \pm}}{m} := \mc{E}_H^\pm \pd{Q_H}{m} F_H \text{.}
\end{equation}
Then, for any $\xi \in \R^d$ and $0 < t_0 < t_1 \leq T$ such that $( t_0, \xi ), ( t_1, \xi ) \in Z_H$, we have
\begin{align}
\label{eq.energy_Ph} | \pd{U_{Hz, +}}{m} ( t_1, \xi ) | &\lesssim | \pd{U_{Hz, +}}{m} ( t_0, \xi ) | + \int_{ t_0 }^{ t_1 } | \pd{F_{Hz, +}}{m} ( \tau, \xi ) | \, d \tau \text{,} \\
\notag | \pd{U_{Hz, -}}{m} ( t_0, \xi ) | &\lesssim | \pd{U_{Hz, -}}{m} ( t_1, \xi ) | + \int_{ t_0 }^{ t_1 } | \pd{F_{Hz, -}}{m} ( \tau, \xi ) | \, d \tau \text{.}
\end{align}
\end{proposition}

\begin{proof}
Similar to Proposition \ref{thm.energy_h}, we use $( z, \xi )$-coordinates and let $z_0 := \zfac \, t_0$, $z_1 := \zfac \, t_1$.
By \eqref{eq.part_diag_h} and Definition \ref{def.system_Ph}, we see that $\pd{U_{Hz, \pm}}{m}$ satisfies
\begin{align*}
\partial_z \pd{U_{Hz, \pm}}{m} &= \mc{E}_H^\pm \bigg[ \imath \mf{Z} D_H + \partial_z \mc{E}_H^\pm + \sum_{ k = 1 }^m \pd{D_H}{k} + \pd{R_H}{m} ( \pd{Q_H}{m} )^{-1} \bigg] ( \mc{E}_H^\pm )^{-1} \pd{U_{Hz, \pm}}{m} + \zfac^{-1} \pd{F_{Hz, \pm}}{m} \\
&= \bigg[ \imath \mf{Z} D_H + \partial_z \mc{E}_H^\pm + \sum_{ k = 1 }^m \pd{D_H}{k} \bigg] \pd{U_{Hz, \pm}}{m} + \mc{E}_H^\pm \pd{R_H}{m} ( \pd{Q_H}{m} )^{-1} ( \mc{E}_H^\pm )^{-1} \pd{U_{Hz, \pm}}{m} + \zfac^{-1} \pd{F_{Hz, \pm}}{m} \text{,}
\end{align*}
where in the last step, we used that $\smash{ \mc{E}_H^\pm W ( \mc{E}_H^\pm )^{-1} = W }$ whenever $W$ is diagonal or $\mc{G}$-diagonal (the latter since $\mc{E}_H^\pm$ is diagonal and has identical values along all the diagonal entries corresponding to a single $\sigma \in \mc{G}$).
Moreover, by \eqref{eq.part_diag_h_DD} and Definition \ref{def.system_Ph}, the above can be written as
\begin{align}
\label{eql.energy_Ph_0} \partial_z \pd{U_{Hz, \pm}}{m} &= [ \imath \mf{Z} D_H + z^{-1} ( B_{ \mc{G} } - B_{ \mc{G}, \pm } ) + \pd{S_D}{m} + \pd{S_R}{m} ] \pd{U_{Hz, \pm}}{m} + \zfac^{-1} \pd{F_{Hz, \pm}}{m} \text{,} \\
\notag \pd{S_D}{m} :\!\!&= R_{ \mc{G} } + \sum_{ k = 2 }^m \pd{D_H}{k} \text{,} \\
\notag \pd{S_R}{m} :\!\!&= \mc{E}_H^\pm \pd{R_H}{m} ( \pd{Q_H}{m} )^{-1} ( \mc{E}_H^\pm )^{-1} \text{.}
\end{align}
Multiplying \eqref{eql.energy_Ph_0} by $\pd{\bar{U}_{Hz, \pm}}{m}$ and applying Lemma \ref{thm.part_diag_est}, we obtain
\begin{align}
\label{eql.energy_Ph_1} \tfrac{1}{2} \partial_z ( | \pd{U_{Hz, +}}{m} |^2 ) &\leq ( | \pd{S_D}{m} | + | \pd{S_R}{m} | ) | \pd{U_{Hz, +}}{m} |^2 + \zfac^{-1} | \pd{F_{Hz, +}}{m} | | \pd{U_{Hz, +}}{m} | \text{,} \\
\notag \tfrac{1}{2} \partial_z ( | \pd{U_{Hz, -}}{m} |^2 ) &\geq ( | \pd{S_D}{m} | + | \pd{S_R}{m} | ) | \pd{U_{Hz, -}}{m} |^2 + \zfac^{-1} | \pd{F_{Hz, -}}{m} | | \pd{U_{Hz, -}}{m} | \text{.}
\end{align}
(Again, the term containing $\imath \mf{Z} D_H$ disappears, since $\imath \mf{Z} D_H$ is skew-Hermitian.)

Now, from Assumption \ref{ass.system_h}, Lemma \ref{thm.part_diag_h}, \eqref{eq.energy_Ph}, and \eqref{eql.energy_Ph_0}, and proceeding analogously to the proof of Proposition \ref{thm.energy_h}, we obtain the integral bounds
\begin{equation}
\label{eql.energy_Ph_10} \int_{ \rho_0^{-1} }^{ \zfac T } | \pd{S_D}{m} ( \zeta, \xi ) | \, d \zeta \lesssim 1 \text{,} \qquad \int_{ \rho_0^{-1} }^{ \zfac T } | \pd{S_R}{m} ( \zeta, \xi ) | \, d \zeta \lesssim 1 \text{,}
\end{equation}
provided $m$ is sufficiently large.
Integrating \eqref{eql.energy_Ph_1} and applying \eqref{eql.energy_Ph_10}, we then obtain
\begin{align*}
| \pd{U_{Hz, +}}{m} ( z_1, \xi ) | &\lesssim | \pd{U_{Hz, +}}{m} ( z_0, \xi ) | + \zfac^{-1} \int_{ z_0 }^{ z_1 } | \pd{F_{Hz, +}}{m} ( \zeta, \xi ) | \, d \zeta \text{,} \\
| \pd{U_{Hz, -}}{m} ( z_0, \xi ) | &\gtrsim | \pd{U_{Hz, -}}{m} ( z_1, \xi ) | + \zfac^{-1} \int_{ z_0 }^{ z_1 } | \pd{F_{Hz, -}}{m} ( \zeta, \xi ) | \, d \zeta \text{,}
\end{align*}
which immediately implies \eqref{eq.energy_Ph} after converting back to $( t, \xi )$-coordinates.
\end{proof}

\begin{remark}
Note that the proof of Proposition \ref{thm.energy_Ph} above implies the following:
\begin{itemize}
\item The derivation of the second part of \eqref{eq.energy_Ph} yields that a solution $\pd{U_{Hz, -}}{m}$ of the ``$-$" system of \eqref{eql.energy_Ph_0} has a finite asymptotic limit as $z \searrow 0$.

\item Similarly, the first part of \eqref{eq.energy_Ph} yields that given asymptotic data at $z = 0$, there exists a solution $\pd{U_{Hz, +}}{m}$ of the ``$+$" system of \eqref{eql.energy_Ph_0} attaining this asymptotic data.
\end{itemize}
The above will be crucial for obtaining asymptotics and scattering for Theorem \ref{thm.energy_ex} below.
\end{remark}

An analogue of Proposition \ref{thm.energy_h_simple} also holds in the current setting:

\begin{proposition} \label{thm.energy_Ph_simple}
Assume the setting of Proposition \ref{thm.energy_Ph}, and suppose in addition that there exist functions $B_{ \mc{G}, \pm, 0 }: \Sph^{ d - 1 } \rightarrow \C^n$ satisfying the following bound for some $\delta > 0$:
\begin{equation}
\label{eq.energy_Ph_simple_ass} | B_{ \mc{G}, \pm, ii } ( t, \omega ) - B_{ \mc{G}, \pm, 0, i } ( \omega ) | \lesssim t^\delta \text{,} \qquad \omega \in \Sph^{ d - 1 } \text{,} \quad 1 \leq i \leq n \text{.}
\end{equation}
Then, the following comparisons hold everywhere on $Z_H$,
\begin{equation}
\label{eq.energy_Ph_simple} | \pd{U_{Hz, \pm}}{m} ( t, \xi ) | \simeq | ( \mc{E}_{ H\ast }^\pm \pd{Q_H}{m} U_H ) ( t, \xi ) | \text{,} \qquad | \pd{F_{Hz, \pm}}{m} ( t, \xi ) | \simeq | ( \mc{E}_{ H\ast }^\pm \pd{Q_H}{m} F_H ) ( t, \xi ) | \text{,}
\end{equation}
where the constants depend on $m$, $\rho_0$, and $B_H$, and where $\mc{E}_{ H\ast }$ is defined on $Z_H$ by
\begin{equation}
\label{eq.energy_Ph_Ehh} \mc{E}_{ H\ast }^\pm ( t, \xi ) := \operatorname{diag} \big( z^{ - B_{ \mc{G}, \pm, 0, 1 } ( \frac{ \xi }{ | \xi | } ) }, \dots, z^{ - B_{ \mc{G}, \pm, 0, n } ( \frac{ \xi }{ | \xi | } ) } \big) \text{.}
\end{equation}
\end{proposition}

\begin{proof}
The proof is analogous to that of Proposition \ref{thm.energy_h_simple}; the key difference is that $B_{ \mc{G}, \pm }$ now plays the same role as the diagonal elements of $B_{ H, \infty }$ in Proposition \ref{thm.energy_h_simple}.
\end{proof}

\subsubsection{The Main Estimates}

Combining Proposition \ref{thm.energy_Ph} with the preceding estimates on $Z_I$ and $Z_P$, we then obtain our main result in the absence of semi-strict hyperbolicity:

\begin{theorem} \label{thm.energy_ex}
Suppose Assumptions \ref{ass.system}, \ref{ass.system_i}, \ref{ass.system_p}, \ref{ass.system_h}, and \ref{ass.partition} hold.
Let $m_P \in \N \cup \{ 0 \}$ and $m_H \in \N$ be large enough, and let $\rho_0$ be sufficiently small (in the same senses as in Theorem \ref{thm.energy_main}).
Moreover, given $U: ( 0, T ]_t \times \R_\xi \rightarrow \C^n$ and $F$ from Assumption \ref{ass.system}, we set
\begin{align}
\label{eq.system_UA_ex} U_{ A, \pm } |_{ Z_P \setminus Z_I } := \pd{U_{Pz}}{m_P} \text{,} \qquad U_{ A, \pm } |_{ Z_I } &:= U \text{,} \qquad U_{ A, \pm } |_{ Z_H \setminus Z_I } := \pd{U_{Hz, \pm}}{m_H} \text{,} \\
\notag F_{ A, \pm } |_{ Z_P \setminus Z_I } := \pd{F_{Pz}}{m_P} \text{,} \qquad F_{ A, \pm } |_{ Z_I } &:= F \text{,} \qquad F_{ A, \pm } |_{ Z_H \setminus Z_I } := \pd{F_{Hz, \pm}}{m_H} \text{.}
\end{align}
\begin{itemize}
\item \emph{Asymptotics:} Given any $u: \R^d_\xi \rightarrow \C^n$, there exists a unique solution $U$ of \eqref{eq.system_gen} satisfying $U ( T, \cdot ) = u$.
If we also assume $F_{ A, - }$ is $t$-integrable,
\begin{equation}
\label{eq.energy_ex_ass_minus} \int_0^T | F_{ A, - } ( \tau, \xi ) | \, d \tau < \infty \text{,} \qquad \xi \in \R^d \text{,}
\end{equation}
then the following asymptotic limits are both well-defined and finite,
\begin{equation}
\label{eq.energy_ex_asymp_limit} u_{ A, - } ( \xi ) := \lim_{ \tau \searrow 0 } U_{ A, - } ( \tau, \xi ) \text{,} \qquad \xi \in \R^d \text{,}
\end{equation}
and the following estimate holds for any $\xi \in \R$ and $t_0 \in ( 0, T ]$:
\begin{equation}
\label{eq.energy_ex_asymp} | u_{ A, - } ( \xi ) | \lesssim | U_{ A, - } ( t_0, \xi ) | + \int_0^{ t_0 } | F_{ A, - } ( \tau, \xi ) | \, d \tau \text{.}
\end{equation}

\item \emph{Scattering:} Suppose $F_{ A, + }$ is $t$-integrable:
\begin{equation}
\label{eq.energy_ex_ass_plus} \int_0^T | F_{ A, + } ( \tau, \xi ) | \, d \tau < \infty \text{,} \qquad \xi \in \R^d \text{.}
\end{equation}
Then, given any $u_{ A, + }: \R^d_\xi \rightarrow \C^n$, there is a unique solution $U$ of \eqref{eq.system_gen} such that
\begin{equation}
\label{eq.energy_ex_scatter_limit} \lim_{ \tau \searrow 0 } U_{ A, + } ( \tau, \xi ) = u_{ A, + } ( \xi ) \text{,} \qquad \xi \in \R^d \text{.}
\end{equation}
In addition, the following estimate holds for any $\xi \in \R$ and $t_0 \in ( 0, T ]$:
\begin{equation}
\label{eq.energy_ex_scatter} | U_{ A, + } ( t_0, \xi ) | \lesssim | u_{ A, + } ( \xi ) | + \int_0^{ t_0 } | F_{ A, + } ( \tau, \xi ) | \, d \tau \text{.}
\end{equation}
\end{itemize}
\end{theorem}

\begin{proof}
The proof is analogous to that of \eqref{eq.energy_asymp} and \eqref{eq.energy_scatter}, except we use $U_{ A, - }$ in the place of $U_A$ for the asymptotics part, and $U_{ A, + }$ in place of $U_A$ for the scattering part.
\end{proof}

\begin{remark}
The existence of the asymptotic limits \eqref{eq.energy_ex_asymp_limit} and \eqref{eq.energy_ex_scatter_limit} is a direct consequence of the estimates of Proposition \ref{thm.energy_Ph}; see also the remark following Proposition \ref{thm.energy_Ph}.
\end{remark}

Thus, in contrast to Theorem \ref{thm.energy_main}, the asymptotics and scattering theories arising from Theorem \ref{thm.energy_ex} will generally differ.
Below, we discuss in more detail the implications of this.

First, let us note the following:
\begin{itemize}
\item By \eqref{eq.system_UA_ex}, the renormalized quantities $U_{ A, \pm }$ are equal on $Z_P$.

\item Since $\bar{Z}_H$ (where $U_{ A, \pm }$ may differ) is disjoint from $t = 0$, we already trivially have existence and uniqueness for the original system \eqref{eq.system_gen} on $Z_H$.
\end{itemize}
The preceding observations immediately imply the following:
\begin{itemize}                                                                                                                                       
\item If we solve forward from asymptotic data $u_A$ at $t = 0$ to $t = T$ using the ``scattering" part of Theorem \ref{thm.energy_ex}, and we then solve backward toward $t = 0$ using the ``asymptotics" part of Theorem \ref{thm.energy_ex}, then the resulting limit at $t = 0$ will simply be the same quantity $u_A$.

\item Conversely, if we solve backward from data $u$ at $t = T$ to $t = 0$ using the ``asymptotics" part of Theorem \ref{thm.energy_ex}, and we then solve forward toward $t = T$ using the ``scattering" part of Theorem \ref{thm.energy_ex}, then we will once again obtain $u$ at $t = T$.
\end{itemize}
Thus, \textit{the asymptotics and scattering parts of Theorem \ref{thm.energy_ex} ``produce the same solutions of \eqref{eq.system_gen}"}.

On the other hand, if $U_{ A, \pm }$ differ on $Z_H$, then Theorem \ref{thm.energy_ex} yields genuinely distinct asymptotics and scattering estimates.
Thus, the upshot---for the homogeneous setting $F \equiv 0$---is that in contrast to the semi-strictly hyperbolic setting of Theorem \ref{thm.energy_main}, \textit{we no longer obtain in general a normed space isomorphism between appropriate data for solutions of \eqref{eq.system_gen} at $t = 0$ and $t = T$}.

\begin{remark}
Observe that one has a scattering isomorphism between $t = 0$ and $t = T$, as in Theorem \ref{thm.energy_ex}, only in the exceptional case that $U_{ A, + } = U_{ A, - }$.
In particular, the above is true if and only if $B_{ \mc{G}, + } = B_{ \mc{G}, - }$.
Note that this property trivially holds when \eqref{eq.system_gen} is semi-strictly hyperbolic.
It will also hold nontrivially for the linearized Einstein-scalar system in Section \ref{sec.einstein}.
\end{remark}

\begin{remark}
Since the estimates on $Z_P$ remain unchanged, the Sobolev norm convergence result of Proposition \ref{thm.energy_conv} also holds for $U_{ A, \pm }$ in the setting of Theorem \ref{thm.energy_ex}.
\end{remark}

\section{Applications} \label{sec.wave}

In this section, we apply our theory from Section \ref{sec.system} to various equations of interest.
These include not only the weakly hyperbolic wave equations \eqref{eq.intro_wave_2}--\eqref{eq.intro_wave_3} traditionally studied in the literature, but also those with critically singular lower-order coefficients, \eqref{eq.intro_wave_4} and \eqref{eq.intro_wave_5}, as well as those with anisotropic degeneracies (e.g.\ arising from Kasner backgrounds).
For all these equations, we quantify the precise loss of regularity and the asymptotics at $t \searrow 0$.
Finally, we briefly discuss how our theory applies to higher-order weakly hyperbolic and singular equations.

\subsection{Wave Equations} \label{sec.wave_basic}

We begin by considering the degenerate-singular wave equations \eqref{eq.intro_wave_5}, in any spatial dimension.
More precisely, we consider the following setting:

\begin{assumption} \label{ass.wave}
Fix $T > 0$, and consider the wave equation on $( 0, T ]_t \times \R_x^d$,
\begin{align}
\label{eq.wave} &\partial_t^2 \phi - t^{ 2 \ell } \sum_{ i, j = 1 }^d a_{ ij } (t) \, \partial_{ x_i x_j }^2 \phi + 2 t^\ell \sum_{ i = 1 }^d b_i (t) \, \partial_{ t x_i }^2 \phi \\
\notag &\quad + t^{ \ell - 1 } \sum_{ i = 1 }^d c_i (t) \, \partial_{ x_i } \phi + t^{-1} g (t) \, \partial_t \phi + t^{-2} h (t) \, \phi = 0 \text{,}
\end{align}
where the unknown is a function $\phi: ( 0, T ]_t \times \R_x^d \rightarrow \C$, where $\ell \in ( -1, \infty )$, and where the coefficients satisfy $a \in \mc{T}^\beta ( [ 0, T ]_t; \R^d \otimes \R^d )$; $b, c \in \mc{T}^\beta ( [ 0, T ]_t; \C^d )$; and $g, h \in \mc{T}^\beta ( [ 0, T ]_t; \C )$, with $\beta > 0$, and
\begin{equation}
\label{eq.wave_elliptic} \sum_{ i, j = 1 }^d a_{ij} (t) \, \xi_i \xi_j \geq \lambda_0 | \xi |^2 \text{,} \qquad ( t, \xi ) \in [ 0, T ] \times \R^d \text{,} \quad \lambda_0 > 0 \text{.}
\end{equation}
\end{assumption}

\begin{remark}
In the following, we will generally assume $\phi$ is sufficiently well-behaved, so that its spatial Fourier transform $\hat{\phi}$ is well-defined, and we can hence study the equation satisfied by $\hat{\phi}$.
\end{remark}

\begin{remark}
Note in the special case when the coefficients $a, b, c, g, h$ are smooth on $[ 0, T ]$, it follows from Definition \ref{def.coeff} that $a, b, c, g, h$ all lie in the spaces $\mc{T}^1$ (i.e.\ $\beta = 1$).
\end{remark}

\begin{definition} \label{def.wave_aux}
We define the auxiliary functions $\mf{a}, \mf{b}, \mf{c}, \alpha: [ 0, T ]_t \times \Sph^{ d - 1 }_\omega \rightarrow \C$ by
\begin{align}
\label{eq.wave_auxf} \mf{a} ( t, \omega ) := \sum_{ i, j = 1 }^d a_{ ij } (t) \, \omega_i \omega_j \text{,} &\qquad \mf{b} ( t, \omega ) := \sum_{ i = 1 }^d b_i (t) \, \omega_i \text{,} \\ 
\notag \mf{c} ( t, \omega ) := \sum_{ i = 1 }^d c_i (t) \, \omega_i \text{,} &\qquad \alpha ( t, \omega ) := \sqrt{ \mf{a} ( t, \omega ) + \mf{b}^2 ( t, \omega ) } \text{.}
\end{align}
In addition, we define the following auxiliary parameters:
\begin{align}
\label{eq.wave_auxp} \gamma &:= \sqrt{ ( g (0) - 1 )^2 - 4 h (0) } \in \C \text{,} \\
\notag \gamma_\pm &:= \tfrac{ 1 + \operatorname{Re} ( g (0) \pm \gamma ) }{ 2 ( \ell + 1 ) } \in \R \text{,} \\
\notag \delta_0 &:= \tfrac{ \operatorname{Re} g (0) - \ell }{ 2 ( \ell + 1 ) } \in \R \text{,} \\
\notag \delta_+ &:= \tfrac{1}{ 2 ( \ell + 1 ) } \sup_{ \omega \in \Sph^{d-1} } \tfrac{ | \operatorname{Re} [ \mf{c} ( 0, \omega ) - ( \ell + g (0) ) \mf{b} ( 0, \omega ) ] | }{ \alpha ( 0, \omega ) } \text{.}
\end{align}
\end{definition}

\begin{remark}
Since $\mf{a} + \mf{b}^2 > 0$ everywhere by \eqref{eq.wave_elliptic}, then $\alpha$ is everywhere smoothly defined.
Moreover, when $g (0)$ or $h (0)$ is complex-valued, then either choice of square root suffices for defining $\gamma$.
\end{remark}

Observe that taking the Fourier transform of \eqref{eq.wave} results in the following equation for $\hat{\phi}$:
\begin{align}
\label{eq.wavef} &\partial_t^2 \hat{\phi} + t^{ 2 \ell } | \xi |^2 \mf{a} \big( t, \tfrac{ \xi }{ | \xi | } \big) \, \hat{\phi} + 2 \imath t^\ell | \xi | \mf{b} \big( t, \tfrac{ \xi }{ | \xi | } \big) \, \partial_t \hat{\phi} \\
\notag &\qquad + \imath t^{ \ell - 1 } | \xi | \mf{c} \big( t, \tfrac{ \xi }{ | \xi | } \big) \, \hat{\phi} + t^{-1} g (t) \, \partial_t \hat{\phi} + t^{-2} h (t) \, \hat{\phi} = 0 \text{.}
\end{align}
To formulate \eqref{eq.wavef} in terms of the development in Section \ref{sec.system}, we first set
\begin{equation}
\label{eq.wave_degen} n := 2 \text{,} \qquad \mf{H} := t^\ell | \xi | = \bigg( \sum_{ i = 1 }^d t^{ 2 \ell } \xi_i^2 \bigg)^\frac{1}{2} \text{,}
\end{equation}
which is consistent with \eqref{eq.system_degen} and \eqref{eq.system_ellip}, with $\ell_i := \ell$ for all $1 \leq i \leq d$.
Next, following Definition \ref{def.z}, we see that the corresponding rescaled time $z$ is given by
\begin{equation}
\label{eq.wave_z} \zfac \simeq \langle \xi \rangle^\frac{1}{ \ell + 1 } \text{,} \qquad z \simeq t \langle \xi \rangle^\frac{1}{ \ell + 1 } \text{.}
\end{equation}
We also fix $\rho_0 \ll_T 1$, and we define $Z_P$, $Z_I$, $Z_H$ using this $\rho_0$ as in Definition \ref{def.zones}.

\begin{definition} \label{def.wave_U}
We define our unknown $U: ( 0, T ]_t \times \R_\xi^{ d } \rightarrow \C^2$ by
\begin{equation}
\label{eq.wave_U_ex} U := \left[ \begin{matrix} t^{-1} [ ( 1 - \chi (z) ) + \imath t^{ \ell + 1 } | \xi | \, \chi (z) ] \, \hat{\phi} \\ \partial_t \hat{\phi} \end{matrix} \right] \text{,}
\end{equation}
where $\chi: \R \rightarrow [ 0, 1 ]$ is a smooth cutoff function satisfying
\begin{equation}
\label{eq.wave_chi} \chi |_{ ( -\infty, ( 2 \rho_0 )^{-1} ] } \equiv 0 \text{,} \qquad \chi |_{ [ \rho_0^{-1}, \infty ) } \equiv 1 \text{.}
\end{equation}
\end{definition}

Observe, in particular, that $U$ takes the following values on $Z_P$ and $Z_H$:
\begin{equation}
\label{eq.wave_U} U |_{ Z_P } := \left[ \begin{matrix} t^{-1} \hat{\phi} \\ \partial_t \hat{\phi} \end{matrix} \right] \text{,} \qquad U |_{ Z_H } := \left[ \begin{matrix} \imath t^\ell | \xi | \hat{\phi} \\ \partial_t \hat{\phi} \end{matrix} \right] \text{.}
\end{equation}
In the following, we will show that the above $U$ satisfies Assumptions \ref{ass.system}, \ref{ass.system_p}, and \ref{ass.system_h}.

\begin{remark}
For conciseness, we consider only a trivial forcing term here, i.e.\ $F \equiv 0$.
However, the analysis extends to nontrivial forcing terms with the obvious modifications.
\end{remark}

\subsubsection{The Intermediate Zone}

Let us first consider $Z_I$.
For convenience, we set
\begin{equation}
\label{eq.wave_i_Xi} \Xi := ( 1 - \chi (z) ) + \imath t^{ \ell + 1 } | \xi | \chi (z) \text{.}
\end{equation}
Note the smallness of $\rho_0$ and \eqref{eq.wave_chi} imply $\chi (z) \equiv 0$ whenever $| \xi | \lesssim 1$.
As a result, using \eqref{eq.wave_z}, \eqref{eq.wave_i_Xi}, and the fact that $z \simeq 1$ on $Z_I$, we can estimate, on $Z_I$,
\begin{align}
\label{eq.wave_i_cutoff} \Xi &\simeq 1 \text{,} \\
\notag \zfac^{-1} | \partial_t \Xi | &\lesssim | \chi' (z) | ( 1 + t^{ \ell + 1 } \langle \xi \rangle ) + ( \ell + 1 ) | \chi (z) | t^\ell \langle \xi \rangle^{\frac{\ell}{\ell + 1} } \\
\notag &\lesssim 1 \text{.}
\end{align}

Thus, combining \eqref{eq.wavef} and \eqref{eq.wave_U_ex}, we obtain the following system on $Z_I$:
\begin{align*}
\partial_t U &= \mc{A} U \text{,} \\
\mc{A} &= \left[ \begin{matrix} - t^{-1} + \Xi^{-1} \partial_t \Xi & t^{-1} \Xi \\ \Xi^{-1} \big[ \mathord{-} t^{ 2 \ell + 1 } | \xi |^2 \mf{a} \big( t, \tfrac{ \xi }{ | \xi | } \big) - \imath t^\ell | \xi | \mf{c} \big( t, \tfrac{ \xi }{ | \xi | } \big) - t^{-1} h (t) \big] & - 2 \imath t^\ell | \xi | \mf{b} \big( t, \tfrac{ \xi }{ | \xi | } \big) - t^{-1} g (t) \end{matrix} \right] \text{.}
\end{align*}
Note by \eqref{eq.wave_z}, \eqref{eq.wave_i_cutoff}, and the above, we see that on $Z_I$,
\[
\zfac^{-1} | \mc{A} | \lesssim 1 \text{.}
\]
In particular, recalling Definition \ref{def.symbol}, the above implies:

\begin{proposition} \label{thm.wave_i_ass}
Assumption \ref{ass.system_i} holds for the wave equation of Assumption \ref{ass.wave}.
\end{proposition}

\subsubsection{The Pseudodifferential Zone}

Next, we consider $Z_P$.
Since $\chi \equiv 0$, then \eqref{eq.wavef} and \eqref{eq.wave_U} yield
\[
\partial_t U = \left[ \begin{matrix} - t^{-1} & t^{-1} \\ - t^{ 2 \ell + 1 } | \xi |^2 \mf{a} \big( t, \tfrac{ \xi }{ | \xi | } \big) - \imath t^\ell | \xi | \mf{c} \big( t, \tfrac{ \xi }{ | \xi | } \big) - t^{-1} h (t) & - 2 \imath t^\ell | \xi | \mf{b} \big( t, \tfrac{ \xi }{ | \xi | } \big) - t^{-1} g (t) \end{matrix} \right] U \text{.} 
\]
In particular, on $Z_P$, we can write $\mc{A}$ as
\begin{align}
\label{eq.wave_p_A} \mc{A} &= t^{-1} \left[ \begin{matrix} -1 & 1 \\ - h(0) & - g(0) \end{matrix} \right] - t^{-1} \left[ \begin{matrix} 0 & 0 \\ h (t) - h (0) & g (t) - g (0) \end{matrix} \right] \\
\notag &\qquad + \left[ \begin{matrix} 0 & 0 \\ - t^{ 2 \ell + 1 } | \xi |^2 \mf{a} \big( t, \tfrac{ \xi }{ | \xi | } \big) - \imath t^\ell | \xi | \mf{c} \big( t, \frac{ \xi }{ | \xi | } \big) & - 2 \imath t^\ell | \xi | \mf{b} \big( t, \tfrac{ \xi }{ | \xi | } \big) \end{matrix} \right] \text{.}
\end{align}
Observe the last two terms on the right-hand side of \eqref{eq.wave_p_A} are remainder terms.
More specifically, recalling \eqref{eq.wave_z} and the smoothness of $g$ and $h$ at $0$, we see that
\begin{equation}
\label{eq.wave_p_AA} \mc{A} = t^{-1} \left[ \begin{matrix} -1 & 1 \\ - h(0) & - g(0) \end{matrix} \right] + \zfac \mc{S}^{ a_P } ( Z_P; \C^2 \otimes \C^2 ) \text{,} \qquad a_P := \min ( \ell, -1 + \alpha ) \text{.}
\end{equation}

The next task is to find a transformation matrix $M_P$ such that
\begin{equation}
\label{eq.wave_Ap} M_P \mc{A} M_P^{-1} + ( \partial_t M_P ) M_P^{-1} = t^{-1} B_P + \zfac R_P \text{,}
\end{equation}
with $M_P, B_P, R_P$ satisfying the conditions in Assumption \ref{ass.system_p}.
The specific choice of $M_P$, and the resulting values of $B_P$ and $R_P$, depend on the value of $\gamma$ from \eqref{eq.wave_auxp}:
\begin{itemize}
\item \emph{Case $\gamma \neq 0$:}
Here, we can choose
\begin{equation}
\label{eq.wave_Mp_1} M_P := \left[ \begin{matrix} \frac{ g (0) - 1 - \gamma }{2} & 1 \\ \frac{ g (0) - 1 + \gamma }{2} & 1 \end{matrix} \right] \text{,} \qquad M_P^{-1} := \frac{1}{ \gamma } \left[ \begin{matrix} -1 & 1 \\ \frac{ g (0) - 1 + \gamma }{2} & - \frac{ g (0) - 1 - \gamma }{2} \end{matrix} \right] \text{.}
\end{equation}
A direct computation yields \eqref{eq.wave_Ap} with the following values:
\begin{equation}
\label{eq.wave_BR_1} B_P := \left[ \begin{matrix} - \frac{ g (0) + 1 + \gamma }{2} & 0 \\ 0 & - \frac{ g (0) + 1 - \gamma }{2} \end{matrix} \right] \text{,} \qquad R_P \in \mc{S}^{ a_P } ( Z_P; \C^2 \otimes \C^2 ) \text{.}
\end{equation}
Also, taking $m_P \in \N_0$ sufficiently large, recalling $\pd{Q_P}{ m_P } \in \mc{S}^0 ( Z_P; \C^2 \otimes \C^2 )$ from Lemma \ref{thm.perf_diag_p}, and recalling Proposition \ref{thm.energy_p}, then the asymptotic quantity on $Z_P$ can be written
\begin{equation}
\label{eq.wave_Upz_1} \pd{U_{Pz}}{ m_P } = \left[ \begin{matrix} \langle \xi \rangle^\frac{ g (0) + 1 + \gamma }{ 2 ( \ell + 1 ) } \hat{\varphi}_+ \\ \langle \xi \rangle^\frac{ g (0) + 1 - \gamma }{ 2 ( \ell + 1 ) } \hat{\varphi}_- \end{matrix} \right] \text{,} \qquad \left[ \begin{matrix} \hat{\varphi}_+ \\ \hat{\varphi}_- \end{matrix} \right] := \left[ \begin{matrix} t^\frac{ g (0) + 1 + \gamma }{2} & \!\!\!\!\!\!\! 0 \\ 0 & \!\!\!\!\!\!\! t^\frac{ g (0) + 1 - \gamma }{2} \end{matrix} \right] \pd{Q_P}{ m_P } \left[ \begin{matrix} \frac{ g (0) - 1 - \gamma }{ 2 t } \hat{\phi} + \partial_t \hat{\phi} \\ \frac{ g (0) - 1 + \gamma }{ 2 t } \hat{\phi} + \partial_t \hat{\phi} \end{matrix} \right] \text{.}
\end{equation}

\item \emph{Case $\gamma = 0$:}
In this setting, we can choose instead
\begin{equation}
\label{eq.wave_Mp_2} M_P := \left[ \begin{matrix} \frac{ 3 - g (0) }{2} & -1 \\ \frac{ g (0) - 1 }{2} & 1 \end{matrix} \right] \text{,} \qquad M_P^{-1} := \left[ \begin{matrix} 1 & 1 \\ \frac{ 1 - g (0) }{2} & \frac{ 3 - g (0) }{2} \end{matrix} \right] \text{.}
\end{equation}
Then, another direct computation with the above $M_P$ yields \eqref{eq.wave_Ap}, but now with
\begin{equation}
\label{eq.wave_BR_2} B_P := \left[ \begin{matrix} - \frac{ g (0) + 1 }{2} & 1 \\ 0 & - \frac{ g (0) + 1 }{2} \end{matrix} \right] \text{,} \qquad R_P \in \mc{S}^{ a_P } ( Z_P; \C^2 \otimes \C^2 ) \text{.}
\end{equation}
Moreover, the asymptotic quantity on $Z_P$ is (again with $m_P \in \N_0$ large enough)
\begin{equation}
\label{eq.wave_Upz_2} \pd{U_{Pz}}{ m_P } = \langle \xi \rangle^\frac{ g (0) + 1 }{ 2 ( \ell + 1 ) } \left[ \begin{matrix} \hat{\varphi}_+ - ( \log z ) \hat{\varphi}_- \\ \hat{\varphi}_- \end{matrix} \right] \text{,} \qquad \left[ \begin{matrix} \hat{\varphi}_+ \\ \hat{\varphi}_- \end{matrix} \right] := t^\frac{ g (0) + 1 }{2} \pd{Q_P}{ m_P } \left[ \begin{matrix} \frac{ 3 - g (0) }{ 2 t } \hat{\phi} - \partial_t \hat{\phi} \\ \frac{ g (0) - 1 }{ 2 t } \hat{\phi} + \partial_t \hat{\phi} \end{matrix} \right] \text{.}
\end{equation}
\end{itemize}

\begin{remark}
Note in particular that $\gamma$ represents the separation between the eigenvalues of $B_P$, or equivalently, the matrix in \eqref{eq.wave_p_AA} that is multiplied to $t^{-1}$.
\end{remark}

Combining all the above, we can then conclude:

\begin{proposition} \label{thm.wave_p_ass}
Assumption \ref{ass.system_p} holds for the wave equation of Assumption \ref{ass.wave}.
\end{proposition}

\begin{remark}
In the special case of a \emph{non-singular} weakly hyperbolic wave equation, namely \eqref{eq.intro_wave_2}, we can derive a more explicit formula for the asymptotic quantities $\hat{\varphi}_\pm$.
Note in this setting,
\[
\ell \in \N \text{,} \qquad g (0) = h (0) = 0 \text{,} \qquad g (t) = O (t) \text{,} \qquad h (t) = O ( t^2 ) \text{,} \qquad \gamma = 1 \text{.}
\]
By computing more carefully the value of $R_P$ from \eqref{eq.wave_p_A} and \eqref{eq.wave_Mp_1},
\begin{align*}
R_P &= \zfac^{-1} M_P \left[ \begin{matrix} 0 & 0 \\ - t^{ 2 \ell + 1 } | \xi |^2 \mf{a} \big( t, \tfrac{ \xi }{ | \xi | } \big) - \imath t^\ell | \xi | \mf{c} \big( t, \tfrac{ \xi }{ | \xi | } \big) - t^{-1} h (t) & - 2 \imath t^\ell | \xi | \mf{b} \big( t, \tfrac{ \xi }{ | \xi | } \big) - t^{-1} g (t) \end{matrix} \right] M_P^{-1} \\
&= \left[ \begin{matrix} -1 & 1 \\ 0 & 1 \end{matrix} \right] \left[ \begin{matrix} 0 & 0 \\ O ( z ) & O (1) \end{matrix} \right] \left[ \begin{matrix} -1 & 1 \\ 0 & 1 \end{matrix} \right] \\
&= \left[ \begin{matrix} O (z) & O (1) \\ O (z) & O (1) \end{matrix} \right] \text{.}
\end{align*}

Recalling \eqref{eq.explicit_Ep_diag}, we then have
\begin{align*}
\mc{E}_P R_P \mc{E}_P^{-1} &= \operatorname{diag} ( z, 1 ) \, R_P \, \operatorname{diag} ( z^{-1}, 1 ) \\
&= \left[ \begin{matrix} O (z) & O (z) \\ O (1) & O (1) \end{matrix} \right] \text{.}
\end{align*}
In particular, the above is $z$-integrable on $Z_P$, so by the remark below Proposition \ref{thm.energy_p}, the results of Proposition \ref{thm.energy_p} hold with $m_P := 0$.
Thus, by \eqref{eq.wave_Upz_1}, the asymptotic quantity on $Z_P$ here is
\begin{equation}
\label{eq.wave_Upz_3} \pd{U_{Pz}}{0} = \left[ \begin{matrix} \langle \xi \rangle^\frac{1}{ \ell + 1 } \hat{\varphi}_+ \\ \hat{\varphi}_- \end{matrix} \right] \text{,} \qquad \left[ \begin{matrix} \hat{\varphi}_+ \\ \hat{\varphi}_- \end{matrix} \right] := \left[ \begin{matrix} - \hat{\phi} + t \partial_t \hat{\phi} \\ \partial_t \hat{\phi} \end{matrix} \right] \text{.}
\end{equation}
\end{remark}

\subsubsection{The Hyperbolic Zone}

Next, on $Z_H$, using \eqref{eq.wavef}, \eqref{eq.wave_U}, and that $\chi \equiv 1$, we obtain
\[
 \partial_t U = \begin{bmatrix} \ell t^{-1} & \imath t^{\ell} | \xi | \\ ( \imath t^\ell |\xi| )^{-1} \big[ \mathord{-} t^{ 2 \ell } | \xi |^2 \mf{a} \big( t, \tfrac{ \xi }{ | \xi | } \big) - \imath t^{ \ell - 1 } | \xi | \mf{c} \big( t, \tfrac{ \xi }{ | \xi | } \big) - t^{-2} h (t) \big] & - 2 \imath t^\ell | \xi | \mf{b} \big( t, \tfrac{ \xi }{ | \xi | } \big) - t^{-1} g (t) \end{bmatrix} U \text{.}
\]
In particular, recalling \eqref{eq.wave_z}, we can then write $\mc{A}$ on $Z_H$ as
\begin{align}
\label{eq.wave_h_A} \mc{A} &= \imath \mf{H} \begin{bmatrix} 0 & 1 \\ \mf{a} \big( t, \tfrac{ \xi }{ | \xi | } \big) & -2 \mf{b} \big( t, \tfrac{ \xi }{ | \xi | } \big) \end{bmatrix} + t^{-1} \begin{bmatrix} \ell & 0 \\ - \mf{c} \big( t, \tfrac{ \xi }{ | \xi | } \big) & - g (t) \end{bmatrix} + \begin{bmatrix} 0 & 0 \\ - ( t^{ \ell + 2 } | \xi | )^{-1} \, \imath h (t) & 0 \end{bmatrix} \\
\notag &= \imath \mf{H} \begin{bmatrix} 0 & 1 \\ \mf{a} \big( t, \tfrac{ \xi }{ | \xi | } \big) & -2 \mf{b} \big( t, \tfrac{ \xi }{ | \xi | } \big) \end{bmatrix} + t^{-1} \begin{bmatrix} \ell & 0 \\ - \mf{c} \big( t, \tfrac{ \xi }{ | \xi | } \big) & - g (t) \end{bmatrix} + \zfac \mc{S}^{ - ( \ell + 2 ) } ( Z_H; \C^2 \otimes \C^2 ) \text{.}
\end{align}

To transform the above, we choose (recalling $\alpha$ from \eqref{eq.wave_auxf})
\begin{align}
\label{eq.wave_Mh} M_H :\!\!&= \frac{1}{2} \begin{bmatrix} ( \mf{b}^2 - \alpha^2 ) \big( t, \tfrac{ \xi }{ | \xi | } \big) & ( \alpha + \mf{b} ) \big( t, \tfrac{ \xi }{ | \xi | } \big) \\ ( \alpha^2 - \mf{b}^2 ) \big( \tfrac{ \xi }{ | \xi | } \big) & ( \alpha - \mf{b} ) \big( t, \tfrac{ \xi }{ | \xi | } \big) \end{bmatrix} \text{,} \\
\notag M_H^{-1} &= \tfrac{1}{ \alpha ( t, \frac{ \xi }{ | \xi | } ) } \begin{bmatrix} - \frac{1}{ ( \alpha + \mf{b} ) ( t, \frac{ \xi }{ | \xi | } ) } & \frac{1}{ ( \alpha - \mf{b} ) ( t, \frac{ \xi }{ | \xi | } ) } \\ 1 & 1 \end{bmatrix} \text{.}
\end{align}
Notice that $M_H, M_H^{-1} \in \mc{S}^0 ( Z_H; \C^2 \otimes \C^2 )$, since \eqref{eq.wave_elliptic} and \eqref{eq.wave_auxf} imply
\[
\alpha \gtrsim 1 \text{,} \qquad \alpha \pm b \gtrsim 1 \text{.}
\]
A direct computation using \eqref{eq.wave_h_A} and \eqref{eq.wave_Mh} yields
\begin{equation}
\label{eq.wave_Ah} M_H \mc{A} M_H^{-1} + ( \partial_t M_H ) M_H^{-1} = \imath \mf{H} D_H + t^{-1} B_H + \zfac R_H \text{,}
\end{equation}
where the quantities on the right-hand side are given by
\begin{align}
\label{eq.wave_Dh} D_H &= \begin{bmatrix} ( - \alpha - \mf{b} ) \big( t, \tfrac{ \xi }{ | \xi | } \big) & 0 \\ 0 & ( \alpha - \mf{b} ) \big( t, \frac{ \xi }{ | \xi | } \big) \end{bmatrix} \text{,} \\
\notag B_H &= \begin{bmatrix} \frac{ [ \ell ( \alpha - \mf{b} ) + \mf{c} ] ( t, \frac{ \xi }{ | \xi | } ) - g (t) ( \alpha + \mf{b} ) ( t, \frac{ \xi }{ | \xi | } ) }{ 2 \alpha ( t, \frac{ \xi }{ | \xi | } ) } & - \frac{ [ \mf{c} ( \alpha + \mf{b} ) ] ( t, \frac{ \xi }{ | \xi | } ) }{ 2 [ \alpha ( \alpha - \mf{b} ) ] ( t, \frac{ \xi }{ | \xi | } ) } - \frac{ ( \ell + g(t) ) ( \alpha + \mf{b} ) ( t, \frac{ \xi }{ | \xi | } ) }{ 2 \alpha ( t, \frac{ \xi }{ | \xi | } ) } \\ \frac{ [ \mf{c} ( \alpha - \mf{b} ) ( t, \frac{ \xi }{ | \xi | } ) }{ 2 [ \alpha ( \alpha + \mf{b} ) ] ( t, \frac{ \xi }{ | \xi | } ) } - \frac{ ( \ell + g(t) ) ( \alpha - \mf{b} ) ( t, \frac{ \xi }{ | \xi | } ) }{ 2 \alpha ( t, \frac{ \xi }{ | \xi | } ) } & \frac{ [ \ell ( \alpha + \mf{b} ) - \mf{c} ] ( t, \frac{ \xi }{ | \xi | } ) - g(t) ( \alpha - \mf{b} ) ( t, \frac{ \xi }{ | \xi | } ) } { 2 \alpha ( t, \frac{ \xi }{ | \xi | } ) } \end{bmatrix} \\
\notag &\qquad + t ( \partial_t M_H ) M_H^{-1} \text{,} \\
\notag R_H &\in \mc{S}^{ -1 - ( \ell + 1 ) } ( Z_H; \C^2 \otimes \C^2 ) \text{.}
\end{align}

Note the eigenvalues of $D_H$ are uniformly separated,
\begin{equation}
\label{eq.wave_ssh} | D_{ H, 22 } ( t, \xi ) - D_{ H, 11 } ( t, \xi ) | = 2 \alpha \big( t, \tfrac{ \xi }{ | \xi | } \big) \gtrsim 1 \text{.}
\end{equation}
Also, since $B_H$ is $| \xi |$-independent, then in terms of \eqref{eq.B_h}, we have, for any $\omega \in \Sph^{ d - 1 }$,
\[
B_{ H, \infty } ( t, \omega ) = \begin{bmatrix} \frac{ [ \ell ( \alpha - \mf{b} ) + \mf{c} ] ( t, \omega ) - g (t) ( \alpha + \mf{b} ) ( t, \omega ) }{ 2 \alpha ( t, \omega ) } & - \frac{ [ \mf{c} ( \alpha + \mf{b} ) ] ( t, \omega ) }{ 2 [ \alpha ( \alpha - \mf{b} ) ] ( t, \omega ) } - \frac{ ( \ell + g(t) ) ( \alpha + \mf{b} ) ( t, \omega ) }{ 2 \alpha ( t, \omega ) } \\ \frac{ [ \mf{c} ( \alpha - \mf{b} ) ( t, \omega ) }{ 2 [ \alpha ( \alpha + \mf{b} ) ] ( t, \omega ) } - \frac{ ( \ell + g(t) ) ( \alpha - \mf{b} ) ( t, \omega ) }{ 2 \alpha ( t, \omega ) } & \frac{ [ \ell ( \alpha + \mf{b} ) - \mf{c} ] ( t, \omega ) - g(t) ( \alpha - \mf{b} ) ( t, \omega ) } { 2 \alpha ( t, \omega ) } \end{bmatrix} + O (t) \text{.}
\]
In particular, from all the above, along with the smoothness of $a, b, c, g$, we have:

\begin{proposition} \label{thm.wave_h_ass}
Assumption \ref{ass.system_h} holds for the wave equation of Assumption \ref{ass.wave}.
Moreover:
\begin{itemize}
\item The system for $U$ is semi-strictly hyperbolic.

\item The assumptions of Proposition \ref{thm.energy_h_simple} hold, with $B_{ H, 0 }$ given by
\begin{align}
\label{eq.wave_Bh_diag} B_{ H, 0, 1 } ( \omega ) &= \tfrac{ \ell - g (0) }{2} + \tfrac{ \mf{c} ( 0, \omega ) - ( \ell + g (0) ) \mf{b} ( 0, \omega ) }{ 2 \alpha ( 0, \omega ) } \text{,} \qquad B_{ H, 0, 2 } ( \omega ) = \tfrac{ \ell - g (0) }{2} - \tfrac{ \mf{c} ( 0, \omega ) - ( \ell + g (0) ) \mf{b} ( 0, \omega ) }{ 2 \alpha ( 0, \omega ) } \text{.}
\end{align}
\end{itemize}
\end{proposition}

Finally, applying Propositions \ref{thm.energy_h} and \ref{thm.energy_h_simple}, then taking $m_H \in \N$ sufficiently large (and setting $\pd{Q_H}{ m_H }$ as in Lemma \ref{thm.perf_diag_h}), the controlled quantity on $Z_H$ satisfies
\begin{align}
\label{eq.wave_Uhz} | \pd{U_{Hz}}{ m_H } | &\simeq \Big| \operatorname{diag} \Big( z^{ - B_{ H, 0, 1 } ( \frac{ \xi }{ | \xi | } ) }, z^{ - B_{ H, 0, 2 } ( \frac{ \xi }{ | \xi | } ) } \Big) \pd{Q_H}{ m_H } U_H \Big| \text{,} \\
\notag U_H &= \frac{1}{2} \begin{bmatrix} \imath t^\ell ( \mf{b}^2 - \alpha^2 ) \big( t, \tfrac{ \xi }{ | \xi | } \big) \, \hat{\phi} + ( \alpha + \mf{b} ) \big( t, \tfrac{ \xi }{ | \xi | } \big) \, \partial_t \hat{\phi} \\ \imath t^\ell ( \alpha^2 - \mf{b}^2 ) \big( t, \tfrac{ \xi }{ | \xi | } \big) \, \hat{\phi} + ( \alpha - \mf{b} ) \big( t, \tfrac{ \xi }{ | \xi | } \big) \, \partial_t \hat{\phi} \end{bmatrix} \text{,}
\end{align}
In particular, noting from \eqref{eq.wave_auxp} and \eqref{eq.wave_Bh_diag} that
\[
\langle \xi \rangle^{ \delta_0 - \delta_+ } t^{ ( \ell + 1 ) ( \delta_0 - \delta_+ ) } \lesssim \Big| z^{ - B_{ H, 0, i } ( \frac{ \xi }{ | \xi | } ) } \Big| \lesssim \langle \xi \rangle^{ \delta_0 + \delta_+ } t^{ ( \ell + 1 ) ( \delta_0 + \delta_+ ) } \text{,} \qquad 1 \leq i \leq 2 \text{,}
\]
we then conclude the following inequalities from \eqref{eq.wave_Uhz} on $Z_H$:
\begin{align}
\label{eq.wave_Uhz_ex} | \pd{U_{Hz}}{ m_H } ( t, \xi ) | &\lesssim \langle \xi \rangle^{ \delta_0 + \delta_+ } \, t^{ ( \ell + 1 ) ( \delta_0 + \delta_+ ) } [ t^\ell \langle \xi \rangle | \hat{\phi} ( t, \xi ) | + | \partial_t \hat{\phi} ( t, \xi ) | ] \text{,} \\
\notag | \pd{U_{Hz}}{ m_H } ( t, \xi ) | &\gtrsim \langle \xi \rangle^{ \delta_0 - \delta_+ } \, t^{ ( \ell + 1 ) ( \delta_0 - \delta_+ ) } [ t^\ell \langle \xi \rangle | \hat{\phi} ( t, \xi ) | + | \partial_t \hat{\phi} ( t, \xi ) | ] \text{.}
\end{align}

\subsubsection{Loss of Regularity}

Finally, we combine the preceding developments with Theorem \ref{thm.energy_main} to obtain our main asymptotics, scattering, and loss of regularity result for the wave equation \eqref{eq.wave}.
To connect with previous literature, we will state less precise estimates on $Z_H$ in terms of $U$ rather than $\pd{U_{Hz}}{m_H}$, with the imprecision manifesting as a loss of regularity for $\phi$.

\begin{theorem} \label{thm.wave_main}
Consider the setting of Assumption \ref{ass.wave}, and let $\gamma, \gamma_\pm, \delta_0, \delta_+$ be as in \eqref{eq.wave_auxp}.
\begin{itemize}
\item \emph{Asymptotics:} Given any $\hat{\phi}_0, \hat{\phi}_1: \smash{ \R_\xi^d } \rightarrow \C$, there exists a unique solution $\hat{\phi}$ of \eqref{eq.wavef} satisfying $\smash{ ( \hat{\phi}, \partial_t \hat{\phi} ) ( T, \cdot ) = ( \hat{\phi}_0, \hat{\phi}_1 ) }$.
Moreover, letting $\hat{\varphi}_\pm$ be as in \eqref{eq.wave_Upz_1} (if $\gamma \neq 0$) or \eqref{eq.wave_Upz_2} (if $\gamma = 0$), then the following asymptotic limits are finite for all $\xi \in \R^d$:
\begin{equation}
\label{eq.wave_asymp_limit} \begin{cases} \hat{\varphi}_{ +, 0 } ( \xi ) := \lim_{ \tau \searrow 0 } \hat{\varphi}_+ ( \tau, \xi ) \text{,} \quad \hat{\varphi}_{ -, 0 } ( \xi ) := \lim_{ \tau \searrow 0 } \hat{\varphi}_- ( \tau, \xi ) & \quad \gamma \neq 0 \text{,} \\ \hat{\varphi}_{ +, 0 } ( \xi ) := \lim_{ \tau \searrow 0 } [ \hat{\varphi}_+ - ( \log z ) \hat{\varphi}_- ] ( \tau, \xi ) \text{,} \quad \hat{\varphi}_{ -, 0 } ( \xi ) := \lim_{ \tau \searrow 0 } \hat{\varphi}_- ( \tau, \xi ) & \quad \gamma = 0 \text{.} \end{cases}
\end{equation}
In addition, the following estimate holds for any $\xi \in \R^{ d }$:
\begin{equation}
\label{eq.wave_asymp} \langle \xi \rangle^{ \gamma_+ } | \hat{\varphi}_{ +, 0 } ( \xi ) | + \langle \xi \rangle^{ \gamma_- } | \hat{\varphi}_{ -, 0 } ( \xi ) | \lesssim \langle \xi \rangle^{ \delta_+ } \big[ \langle \xi \rangle^{ 1 + \delta_0 } | \hat{\phi} ( T, \xi ) | + \langle \xi \rangle^{ \delta_0 } | \partial_t \hat{\phi} ( T, \xi ) | \big] \text{.}
\end{equation}

\item \emph{Scattering:} Given $\hat{\varphi}_{ \pm, 0 }: \R_\xi^d \rightarrow \C$, there is a unique solution $\hat{\phi}$ of \eqref{eq.wavef} such that, defining $\hat{\varphi}_\pm$ as in \eqref{eq.wave_Upz_1} ($\gamma \neq 0$) or \eqref{eq.wave_Upz_2} ($\gamma = 0$), the following limits hold for all $\xi \in \R^d$:
\begin{equation}
\label{eq.wave_scatter_limit} \begin{cases} \lim_{ \tau \searrow 0 } \hat{\varphi}_+ ( \tau, \xi ) = \hat{\varphi}_{ +, 0 } ( \xi ) \text{,} \quad \lim_{ \tau \searrow 0 } \hat{\varphi}_- ( \tau, \xi ) = \hat{\varphi}_{ -, 0 } ( \xi ) & \quad \gamma \neq 0 \text{,} \\ \lim_{ \tau \searrow 0 } [ \hat{\varphi}_+ - ( \log z ) \hat{\varphi}_- ] ( \tau, \xi ) = \hat{\varphi}_{ +, 0 } ( \xi ) \text{,} \quad \lim_{ \tau \searrow 0 } \hat{\varphi}_- ( \tau, \xi ) = \hat{\varphi}_{ -, 0 } ( \xi ) & \quad \gamma = 0 \text{.} \end{cases}
\end{equation}
In addition, the following estimate holds for any $\xi \in \R^{ d }$:
\begin{equation}
\label{eq.wave_scatter} \langle \xi \rangle^{ 1 + \delta_0 } | \hat{\phi} ( T, \xi ) | + \langle \xi \rangle^{ \delta_0 } | \partial_t \hat{\phi} ( T, \xi ) | \lesssim \langle \xi \rangle^{ \delta_+ } \big[ \langle \xi \rangle^{ \gamma_+ } | \hat{\varphi}_{ +, 0 } ( \xi ) | + \langle \xi \rangle^{ \gamma_- } | \hat{\varphi}_{ -, 0 } ( \xi ) | \big] \text{.}
\end{equation}
\end{itemize}
\end{theorem}

\begin{proof}
Since Assumptions \ref{ass.system_i}, \ref{ass.system_p}, \ref{ass.system_h} hold in our setting, and since our system is also semi-strictly hyperbolic (see Propositions \ref{thm.wave_i_ass}, \ref{thm.wave_p_ass}, \ref{thm.wave_h_ass}), then Theorem \ref{thm.energy_main}---for sufficiently large $m_P, m_H$ and small $\rho_0$---and \eqref{eq.wave_U_ex} imply the existence and uniqueness of $\smash{ \hat{\phi} }$ for both the asymptotics and scattering settings.
Moreover, from \eqref{eq.system_UA}, along with \eqref{eq.wave_Upz_1} and \eqref{eq.wave_Upz_2} relating $\hat{\varphi}_\pm$ to $\pd{U_{Pz}}{m_P}$, we see that the limits \eqref{eq.wave_asymp_limit} and \eqref{eq.wave_scatter_limit} hold.
Thus, it remains only to prove the estimates \eqref{eq.wave_asymp} and \eqref{eq.wave_scatter}.

By Proposition \ref{thm.zones_ass}, we can choose $\rho_0$ small enough so that $\{ t = T \}$ lies entirely within $Z_H \cup Z_I$.
First, by \eqref{eq.system_UA}, \eqref{eq.wave_Upz_1}, and \eqref{eq.wave_Upz_2}, we have (using the language of Theorem \ref{thm.energy_main})
\begin{equation}
\label{eql.wave_main_1} | u_A ( \xi ) | \simeq \langle \xi \rangle^{ \gamma_+ } | \hat{\varphi}_{ +, 0 } ( \xi ) | + \langle \xi \rangle^{ \gamma_- } | \hat{\varphi}_{ -, 0 } ( \xi ) | \text{.}
\end{equation}
Next, if $( T, \xi ) \in Z_I$, then \eqref{eq.wave_U_ex}, \eqref{eq.wave_chi}, and that $z \simeq 1$ on $Z_I$ yield
\begin{equation}
\label{eql.wave_main_2} | U ( T, \xi ) | \simeq_q \langle \xi \rangle^q [ \langle \xi \rangle | \hat{\phi} ( T, \xi ) | + | \partial_t \hat{\phi} ( T, \xi ) | ] \text{,} \qquad q \in \R \text{.}
\end{equation}
(In particular, note that the smallness of $\rho_0$ implies $\chi (z) \equiv 0$ whenever $| \xi | \lesssim 1$.)
Furthermore, if $( T, \xi ) \in Z_H$, then by \eqref{eq.system_UA} and \eqref{eq.wave_Uhz_ex}, we have the estimates
\begin{align}
\label{eql.wave_main_3} | U_A ( T, \xi ) | &\lesssim \langle \xi \rangle^{ \delta_0 + \delta_+ } [ \langle \xi \rangle | \hat{\phi} ( T, \xi ) | + | \partial_t \hat{\phi} ( T, \xi ) | ] \text{,} \\
\notag | U_A ( T, \xi ) | &\gtrsim \langle \xi \rangle^{ \delta_0 - \delta_+ } [ \langle \xi \rangle | \hat{\phi} ( T, \xi ) | + | \partial_t \hat{\phi} ( T, \xi ) | ] \text{.}
\end{align}
Since \eqref{eq.energy_asymp} and \eqref{eq.energy_scatter} imply
\[
| U_A ( T, \xi ) | \simeq | u_A ( \xi ) | \text{,}
\]
then both \eqref{eq.wave_asymp} and \eqref{eq.wave_scatter} follow from \eqref{eql.wave_main_1}--\eqref{eql.wave_main_3} and the above.
\end{proof}

As mentioned earlier, for non-singular wave equations, we can obtain a more precise statement matching the optimal results of the existing literature:

\begin{corollary} \label{thm.wave_nonsing}
Assume the setting of Theorem \ref{thm.wave_main}, and suppose in addition that
\begin{equation}
\label{eq.wave_nonsing} \ell \in \N \text{,} \qquad g (0) = 0 \text{,} \qquad h (0) = h' (0) = 0 \text{,}
\end{equation}
that is, we consider non-singular equations of the form \eqref{eq.intro_wave_2}.
In addition, let
\begin{equation}
\label{eq.wave_nonsing_aux} \mu := \tfrac{1}{ 2 ( \ell + 1 ) } \sup_{ \omega \in \Sph^{ d - 1 } } \tfrac{ | \operatorname{Re} [ \mf{c} ( 0 , \omega ) - \ell \, \mf{b} ( 0, \omega ) ] | }{ \alpha ( 0, \omega ) } \text{.}
\end{equation}
Then, the following statements hold:
\begin{itemize}
\item \emph{Asymptotics:} For any $\hat{\phi}_0, \hat{\phi}_1: \smash{ \R_\xi^d } \rightarrow \C$, there exists a unique solution $\hat{\phi}$ of \eqref{eq.wavef} satisfying $\smash{ ( \hat{\phi}, \partial_t \hat{\phi} ) ( T, \cdot ) = ( \hat{\phi}_0, \hat{\phi}_1 ) }$.
Furthermore, the following asymptotic limits are finite:
\begin{equation}
\label{eq.wave_nonsing_asymp_limit} \hat{\phi} ( 0, \xi ) := \lim_{ \tau \searrow 0 } \hat{\phi} ( \tau, \xi ) \text{,} \qquad \partial_t \hat{\phi} ( 0, \xi ) := \lim_{ \tau \searrow 0 } \partial_t \hat{\phi} ( \tau, \xi ) \text{,} \qquad \xi \in \R^d \text{.}
\end{equation}
In addition, the following estimate holds for any $\xi \in \R^d$:
\begin{equation}
\label{eq.wave_nonsing_asymp} \langle \xi \rangle^\frac{1}{ \ell + 1 } | \hat{\phi} ( 0, \xi ) | + | \partial_t \hat{\phi} ( 0, \xi ) | \lesssim \langle \xi \rangle^\mu \big[ \langle \xi \rangle^{ 1 - \frac{ \ell }{ 2 ( \ell + 1 ) } } | \hat{\phi} ( T, \xi ) | + \langle \xi \rangle^{ - \frac{ \ell }{ 2 ( \ell + 1 ) } } | \partial_t \hat{\phi} ( T, \xi ) | \big] \text{.}
\end{equation}

\item \emph{Scattering:} Given $\hat{\phi}_0, \hat{\phi}_1: \R_\xi^d \rightarrow \C$, there is a unique solution $\hat{\phi}$ of \eqref{eq.wavef} such that
\begin{equation}
\label{eq.wave_nonsing_scatter_limit} \hat{\phi} ( 0, \xi ) := \lim_{ \tau \searrow 0 } \hat{\phi} ( \tau, \xi ) = \hat{\phi}_0 ( \xi ) \text{,} \qquad \partial_t \hat{\phi} ( 0, \xi ) := \lim_{ \tau \searrow 0 } \partial_t \hat{\phi} ( \tau, \xi ) = \hat{\phi}_1 ( \xi ) \text{,} \qquad \xi \in \R^d \text{.}
\end{equation}
In addition, the following estimate holds for any $\xi \in \R^{d}$:
\begin{equation}
\label{eq.wave_nonsing_scatter} \langle \xi \rangle^{ 1 - \frac{ \ell }{ 2 ( \ell + 1 ) } } | \hat{\phi} ( T, \xi ) | + \langle \xi \rangle^{ - \frac{ \ell }{ 2 ( \ell + 1 ) } } | \partial_t \hat{\phi} ( T, \xi ) | \lesssim \langle \xi \rangle^\mu \big[ \langle \xi \rangle^\frac{1}{ \ell + 1 } | \hat{\phi} ( 0, \xi ) | + | \partial_t \hat{\phi} ( 0, \xi ) | \big] \text{.}
\end{equation}
\end{itemize}
\end{corollary}

\begin{proof}
Note that in this setting, the quantities from \eqref{eq.wave_auxp} reduce to
\begin{equation}
\label{eql.wave_nonsing_1} \gamma = 1 \text{,} \qquad \gamma_\pm = \tfrac{ 1 \pm 1 }{ 2 ( \ell + 1 ) } \text{,} \qquad \delta_0 = - \tfrac{ \ell }{ 2 ( \ell + 1 ) } \text{,} \qquad \delta_+ = \mu \text{.}
\end{equation}
Moreover, from the remark following Proposition \ref{thm.wave_p_ass}, we can apply Theorem \ref{thm.wave_main} with $m_P := 0$.
In particular, recalling \eqref{eq.wave_Upz_3}, we can characterize our asymptotic quantities as
\[
| \hat{\varphi}_+ | + | \hat{\varphi}_- | \simeq | \hat{\phi} | + | \partial_t \hat{\phi} | \text{.}
\]
The desired results now follow from Theorem \ref{thm.wave_main}, \eqref{eql.wave_nonsing_1}, and the above.
\end{proof}

\begin{remark}
Observe in particular that for the model wave equation \eqref{eq.intro_wave_1} in one spatial dimension, the estimates \eqref{eq.wave_nonsing_asymp} and \eqref{eq.wave_nonsing_scatter} reduce to the optimal \eqref{eq.intro_wave_1scat} and \eqref{eq.intro_wave_1asymp}.
\end{remark}

\begin{remark}
\eqref{eq.wave_nonsing_asymp}, \eqref{eq.wave_nonsing_scatter} also recover the regularity loss obtained in \cite{dreher_reissig_2000, dreher_witt_2002, dreher_witt_2005}, with a slight improvement in the estimates on $Z_H$---here, we measure the loss of regularity only from the infinite-frequency limit of $B_H$ at $t = 0$, rather than a supremum of $B_H$ over all frequencies.
\end{remark}

\subsection{Anisotropic Settings} \label{sec.wave_kasner}

Next, we turn our attention to scalar wave equations with anisotropic degeneracies, such as \eqref{eq.intro_wave_kasner}.
Here, we restrict to only a subclass of coefficients in order to focus on the anisotropy and simplify computations.
More specifically, we consider the following setting:

\begin{assumption} \label{ass.kasner}
Fix $T > 0$, and consider the wave equation on $(0,T]_t \times \R^d_x$,
\begin{equation}
\label{eq.kasner} \partial_t^2 \phi - \sum_{i=1}^d t^{ 2 \ell_i } \partial_{x_i}^2 \phi + \sum_{i=1}^d c_i t^{ \ell_i - 1 } \partial_{x_i} \phi + g t^{-1} \partial_t \phi = 0 \text{,}
\end{equation}
where the unknown is $\phi : (0,T]_t \times \R^d_x \rightarrow \C$, and where $\ell \in ( -1, \infty )^d$, $c \in \C^d$, and $g \in \C$.
\end{assumption}

\begin{remark}
The analysis here extends readily to more general wave equations of the form
\begin{align}
\label{eq.kasner_gen} &\partial_t^2 \phi - \sum_{ i, j = 1 }^d t^{ \ell_i + \ell_j } \lambda_{ij} (t) \, \partial_{ x_i x_j }^2 \phi + 2 \sum_{ i = 1 }^d t^{ \ell_i } b_i (t) \, \partial_{ t x_i } \phi \\
\notag &\quad + \sum_{i=1}^d t^{ \ell_i - 1 } c_i (t) \, \partial_{x_i} \phi + t^{-1} g (t) \, \partial_t \phi + t^{-2} h (t) \, \phi = 0 \text{,}
\end{align}
using similar computations as in Section \ref{sec.wave_basic}, provided the condition \eqref{eq.system_ellip} holds.
\end{remark}

\begin{remark}
Taking $c \equiv 0$, $g = 1$ in \eqref{eq.kasner} results in the wave equation \eqref{eq.intro_wave_kasner} on Kasner backgrounds.
Therefore, \eqref{eq.kasner} extends the Kasner wave equation \eqref{eq.intro_wave_kasner} to encompass a general catalog of critical weakly hyperbolic and singular behaviors at $t \searrow 0$, characterized by the parameters $\ell$, $c$, $g$.
\end{remark}

To connect Assumption \ref{ass.kasner} to Section \ref{sec.system}, we first set
\begin{equation}
\label{eq.kasner_degen} n := 2 \text{,} \qquad \ell_\ast := \min ( \ell_1, \dots, \ell_d ) \text{,} \qquad \mf{H} := \bigg[ \sum_{ i = 1 }^d t^{ 2 \ell_i } \xi_i^2 \bigg]^\frac{1}{2} \text{,} \qquad \mf{C} := \sum_{ i = 1 }^d c_i t^{ \ell_i - 1 } \xi_i \text{.}
\end{equation}
Note $\mf{H}$ is consistent with \eqref{eq.system_degen} and trivially satisfies \eqref{eq.system_ellip}.
Observe also from \eqref{eq.kasner_degen} that
\begin{equation}
\label{eq.kasner_degen_deriv} | t \partial_t \mf{H} | \lesssim \mf{H} \text{,} \qquad | t \mf{C} | \lesssim \mf{H} \text{.}
\end{equation}
Moreover, the rescaled $t$ is given directly by Definition \ref{def.z}:
\begin{equation}
\label{eq.kasner_z} \zfac := \max_{ 1 \leq i \leq d } \langle \xi_i \rangle^\frac{1}{ \ell_i + 1 } \text{,} \qquad z := \zfac t \text{.}
\end{equation}
We again fix $\rho_0 \ll_T 1$ and define $Z_P$, $Z_I$, $Z_H$ using this $\rho_0$ as in Definition \ref{def.zones}.

Assuming $\phi$ is well-behaved, then its Fourier transform satisfies the following on $(0,T]_t \times \R_\xi$:
\begin{equation}
\label{eq.kasnerf} \partial_t^2 \hat{\phi} + ( \mf{H}^2 + \imath \mf{C} ) \, \hat{\phi} + g t^{-1} \, \partial_t \hat{\phi} = 0 \text{.}
\end{equation}
To formulate \eqref{eq.kasnerf} as a first-order system, we set the following:

\begin{definition} \label{def.kasner_U}
	We define our unknown $U: ( 0, T ]_t \times \R_\xi^{ d } \rightarrow \C^2$ by
\begin{equation}
\label{eq.kasner_U} U := \begin{bmatrix} t^{-1} [ ( 1 - \chi (z) ) + \imath t \mf{H} \, \chi (z) ] \, \hat{\phi} \\ \partial_t \hat{\phi} \end{bmatrix} \text{,}
\end{equation}
where $\chi: \R \rightarrow [ 0, 1 ]$ is a smooth cutoff function satisfying \eqref{eq.wave_chi}.
\end{definition}

\subsubsection{The Intermediate Zone}

The analysis on $Z_I$ proceeds analogously as in Section \ref{sec.wave_basic}.
We set
\begin{equation}
\label{eq.kasner_i_Xi} \Xi := ( 1 - \chi (z) ) + \imath t \mf{H} \chi (z) \text{.}
\end{equation}
From Proposition \ref{thm.degen_est}, \eqref{eq.kasner_degen_deriv}, \eqref{eq.kasner_z}, \eqref{eq.kasner_i_Xi}, and that $z \simeq 1$ on $Z_I$, we obtain
\begin{equation}
\label{eq.kasner_i_cutoff} \Xi \simeq 1 \text{,} \qquad \zfac^{-1} | \partial_t \Xi | \lesssim 1 \text{.}
\end{equation}

Combining \eqref{eq.kasnerf} and \eqref{eq.kasner_U} yields the following system on $Z_I$:
\[
\partial_t U = \mc{A} U \text{,} \qquad \mc{A} = \begin{bmatrix} - t^{-1} + \Xi^{-1} \partial_t \Xi & t^{-1} \Xi \\ \Xi^{-1} ( - t \mf{H}^2 - \imath t \mf{C} ) & - g t^{-1} \end{bmatrix} \text{.}
\]
Applying Proposition \ref{thm.degen_est}, \eqref{eq.kasner_degen_deriv}, and \eqref{eq.kasner_i_cutoff} yields $\zfac^{-1} | \mc{A} | \lesssim 1$ on $Z_I$, hence:

\begin{proposition} \label{thm.kasner_i_ass}
Assumption \ref{ass.system_i} holds for the wave equation of Assumption \ref{ass.kasner}.
\end{proposition}

\subsubsection{The Pseudodifferential Zone}

On $Z_P$, we see from \eqref{eq.kasnerf} and \eqref{eq.kasner_U} that
\begin{equation}
\label{eq.kasner_U_p} \partial_t U = \begin{bmatrix} -t^{-1} & t^{-1} \\ - t \mf{H}^2 - \imath t \mf{C} & - g t^{-1} \end{bmatrix} U \text{,} \qquad U = \begin{bmatrix} t^{-1} \hat{\phi} \\ \partial_t \hat{\phi} \end{bmatrix} \text{,}
\end{equation}
so that we can write $\mc{A}$ on $Z_P$ as
\begin{align}
\label{eq.kasner_p_A} \mc{A} &= t^{-1} \begin{bmatrix} -1 & 1 \\ 0 & -g \end{bmatrix} + \begin{bmatrix} 0 & 0 \\ - t \mf{H}^2 - \imath t \mf{C} & 0 \end{bmatrix} \\
\notag &= t^{-1} \begin{bmatrix} -1 & 1 \\ 0 & -g \end{bmatrix} + \zfac \mc{S}^{ \ell_\ast } ( Z_P; \C^2 \otimes \C^2 ) \text{,}
\end{align}
where we recalled Proposition \ref{thm.degen_est}, \eqref{eq.kasner_degen_deriv}, and \eqref{eq.kasner_z} in the last step.

To transform $\mc{A}$ into the form \eqref{eq.wave_Ap}, we again split into cases:
\begin{itemize}
\item \emph{Case $g \neq 1$:} Similar to Section \ref{sec.wave_basic}, here we can take
\begin{equation}
\label{eq.kasner_BR_1} M_P := \begin{bmatrix} g - 1 & 1 \\ 0 & 1 \end{bmatrix} \text{,} \qquad B_P := \begin{bmatrix} -1 & 0 \\ 0 & -g \end{bmatrix} \text{,} \qquad R_P \in \mc{S}^{ \ell_\ast } ( Z_P; \C^2 \otimes \C^2 ) \text{.}
\end{equation}
Taking $m_P \in \N_0$ large enough and recalling Proposition \ref{thm.energy_p}, we then have 
\begin{equation}
\label{eq.kasner_Upz_1} \pd{U_{Pz}}{ m_P } = \begin{bmatrix} \zfac \hat{\varphi}_+ \\ \zfac^g \hat{\varphi}_- \end{bmatrix} \text{,} \qquad \begin{bmatrix} \hat{\varphi}_+ \\ \hat{\varphi}_- \end{bmatrix} := \begin{bmatrix} t & 0 \\ 0 & t^g \end{bmatrix} \pd{Q_P}{ m_P } \left[ \begin{matrix} \frac{ g - 1 }{ t } \hat{\phi} + \partial_t \hat{\phi} \\ \partial_t \hat{\phi} \end{matrix} \right] \text{.}
\end{equation}

\item \emph{Case $g = 1$:} The leading part is already in Jordan normal form, so we can take
\begin{equation}
\label{eq.kasner_BR_2} M_P := \begin{bmatrix} 1 & 0 \\ 0 & 1 \end{bmatrix} \text{,} \qquad B_P := \begin{bmatrix} -1 & 1 \\ 0 & -1 \end{bmatrix} \text{,} \qquad R_P \in \mc{S}^{ \ell_\ast } ( Z_P; \C^2 \otimes \C^2 ) \text{.}
\end{equation}
In fact, from \eqref{eq.system_Ep} and \eqref{eq.kasner_p_A}, we can more explicitly compute
\begin{align*}
\mc{E}_P R_P \mc{E}_P^{-1} &= \begin{bmatrix} z & - z ( \log z ) \\ 0 & z \end{bmatrix} \begin{bmatrix} 0 & 0 \\ \zfac^{-1} ( - t \mf{H}^2 - \imath t \mf{C} ) & 0 \end{bmatrix} \begin{bmatrix} z^{-1} & z^{-1} ( \log z ) \\ 0 & z^{-1} \end{bmatrix} \\
&= O ( z^{ \ell_\ast } ( \log z )^2 ) \text{,}
\end{align*}
which is $z$-integrable on $Z_P$, so (see the remark after Proposition \ref{thm.energy_p}) we can take $m_P := 0$:
\begin{equation}
\label{eq.kasner_Upz_2} \pd{U_{Pz}}{0} = \zfac \begin{bmatrix} \hat{\varphi}_+ - ( \log z ) \hat{\varphi}_- \\ \hat{\varphi}_- \end{bmatrix} \text{,} \qquad \begin{bmatrix} \hat{\varphi}_+ \\ \hat{\varphi}_- \end{bmatrix} := \left[ \begin{matrix} \hat{\phi} \\ t \partial_t \hat{\phi} \end{matrix} \right] \text{.}
\end{equation}
\end{itemize}
From all the above, we hence conclude:

\begin{proposition} \label{thm.kasner_p_ass}
Assumption \ref{ass.system_p} holds for the wave equation of Assumption \ref{ass.kasner}.
\end{proposition}

\begin{remark}
For small enough $| g - 1 |$ (depending on $\ell_\ast$), we can again take $m_P := 0$ in the above.
\end{remark}

\subsubsection{The Hyperbolic Zone}

On $Z_H$, we obtain from \eqref{eq.kasnerf} and \eqref{eq.kasner_U} that
\begin{equation}
\label{eq.kasner_U_h} \partial_t U = \begin{bmatrix} \mf{H}^{-1} \partial_t \mf{H} & \imath \mf{H} \\ ( \imath \mf{H} )^{-1} ( - \mf{H}^2 - \imath \mf{C} ) & - g t^{-1} \end{bmatrix} U \text{,} \qquad U := \begin{bmatrix} \imath \mf{H} \hat{\phi} \\ \partial_t \hat{\phi} \end{bmatrix} \text{.}
\end{equation}
In particular, we can write $\mc{A}$ on $Z_H$ as 
\begin{equation}
\label{eq.kasner_h_A} \mc{A} = \imath \mf{H} \begin{bmatrix} 0 & 1 \\ 1 & 0 \end{bmatrix} + t^{-1} \begin{bmatrix} \mf{H}^{-1} ( t \partial_t \mf{H} ) & 0 \\ - \mf{H}^{-1} ( t \mf{C} ) & -g \end{bmatrix} \text{.}
\end{equation}

To transform $\mc{A}$ into the form \eqref{eq.wave_Ah}, we can proceed as in Section \ref{sec.wave_basic} and take
\begin{align}
\label{eq.kasner_Dh} M_H &:= \frac{1}{2} \begin{bmatrix} -1 & 1 \\ 1 & 1 \end{bmatrix} \text{,} \\
\notag D_H &:= \begin{bmatrix} -1 & 0 \\ 0 & 1 \end{bmatrix}, \\
\notag B_H &:= \frac{1}{2} \begin{bmatrix}
\mf{H}^{-1} ( t \partial_t \mf{H} ) + \mf{H}^{-1} ( t \mf{C} ) - g & - \mf{H}^{-1} ( t \partial_t \mf{H} ) - \mf{H}^{-1} ( t \mf{C} ) - g \\ - \mf{H}^{-1} ( t \partial_t \mf{H} ) + \mf{H}^{-1} ( t \mf{C} ) - g & \mf{H}^{-1} ( t \partial_t \mf{H} ) - \mf{H}^{-1} ( t \mf{C} ) - g \end{bmatrix} \text{.}
\end{align}
Observe in particular that the eigenvalues of $D_H$ are uniformly separated, and that $B_H$ from \eqref{eq.kasner_Dh} is already $| \xi |$-independent.
Combining the above, we then conclude the following:

\begin{proposition} \label{thm.kasner_h_ass}
Assumption \ref{ass.system_h} holds for the wave equation of Assumption \ref{ass.wave}.
Furthermore, we have that the system for $U$ is semi-strictly hyperbolic.
\end{proposition}

To further analyze $B_H$, let us consider the vector-valued function
\begin{equation}
\label{eq.kasner_Bh_theta} \vec{\mf{H}} := ( t^{ \ell_1 } \xi_1, \dots, t^{ \ell_d } \xi_d ) \text{,}
\end{equation}
and we let $\theta$ denote the angle between the vectors $c$ and $\smash{\vec{\mf{H}}}$.
Then, by \eqref{eq.kasner_degen} and \eqref{eq.kasner_Bh_theta},
\begin{equation}
\label{eq.kasner_Bh_1} \mf{H}^{-1} ( t \partial_t \mf{H} ) = t \partial_t ( \log \mf{H} ) \text{,} \qquad \mf{H}^{-1} ( t \mf{C} ) = \big| \vec{ \mf{H} } \big|^{-1} \big( c \cdot \vec{\mf{H}} \big) = | c | \cos \theta \text{.}
\end{equation}
As a result, from \eqref{eq.kasner_Dh} and \eqref{eq.kasner_Bh_1}, we obtain
\[
B_{ H, 11 } = \tfrac{1}{2} [ t \partial_t ( \log \mf{H} ) + | c | \cos \theta - g ] \text{,} \qquad B_{ H, 22 } = \tfrac{1}{2} [ t \partial_t ( \log \mf{H} ) - | c | \cos \theta - g ] \text{,}
\]
and Definition \ref{def.system_Eh}, applied to this setting, then yields, for any $( t, \xi ) \in Z_H$,
\begin{align}
\label{eq.kasner_Bh_2} b_{ H, 1 } ( t, \xi ) &= \tfrac{1}{2} \log \tfrac{ \mf{H} ( t, \xi ) }{ \mf{H} ( \zfac^{-1} \rho_0^{-1}, \xi ) } + \tfrac{1}{2} \int_{ \zfac^{-1} \rho_0^{-1} }^t \tfrac{ | c | \cos \theta ( \tau, \xi ) }{ \tau } \, d \tau - \tfrac{1}{2} g \log ( \rho_0^{-1} z ) \text{,} \\
\notag b_{ H, 2 } ( t, \xi ) &= \tfrac{1}{2} \log \tfrac{ \mf{H} ( t, \xi ) }{ \mf{H} ( \zfac^{-1} \rho_0^{-1}, \xi ) } - \tfrac{1}{2} \int_{ \zfac^{-1} \rho_0^{-1} }^t \tfrac{ | c | \cos \theta ( \tau, \xi ) }{ \tau } \, d \tau - \tfrac{1}{2} g \log ( \rho_0^{-1} z ) \text{,}
\end{align}
and Definition \ref{def.system_Eh} and Proposition \ref{thm.energy_h} then imply the following for large enough $m_H \in \N$:
\begin{align}
\label{eq.kasner_Uhz} \pd{U_{Hz}}{m_H} = \operatorname{diag} \big( e^{ - b_{ H, 1 } }, e^{ - b_{ H, 2 } } \big) \pd{Q_H}{m_H} \begin{bmatrix} \tfrac{1}{2} ( \partial_t \hat{\phi} - \imath \mf{H} \hat{\phi} ) \\ \tfrac{1}{2} ( \partial_t \hat{\phi} + \imath \mf{H} \hat{\phi} ) \end{bmatrix} \text{.}
\end{align}

Noting in particular that (see \eqref{eq.kasner_degen} and \eqref{eq.kasner_z})
\[
- \tfrac{ | c | }{t} \leq \tfrac{ | c | \cos \theta ( t, \xi ) }{t} \leq \tfrac{ | c | }{t} \text{,} \qquad \mf{H} ( \zfac^{-1} \rho_0^{-1}, \xi ) \simeq \zfac \text{,} \qquad ( t, \xi ) \in Z_H \text{,}
\]
we then conclude from \eqref{eq.kasner_Bh_2} and the above that
\begin{equation}
\label{eq.kasner_Bh_3} \sqrt{ \tfrac{ \zfac }{ \mf{H} ( t, \xi ) } } \, z^\frac{ g - | c | }{2} \lesssim e^{ - b_{ H, i } ( t, \xi ) } \lesssim \sqrt{ \tfrac{ \zfac }{ \mf{H} ( t, \xi ) } } \, z^\frac{ g + | c | }{2} \text{.}
\end{equation}
Finally, combining \eqref{eq.kasner_Uhz} with \eqref{eq.kasner_Bh_3} results in the following estimates for $\pd{U_{Hz}}{m_H}$:
\begin{align}
\label{eq.kasner_Uhz_ex} | \pd{U_{Hz}}{m_H} ( t, \xi ) | &\lesssim \zfac^\frac{ 1 + g + |c| }{2} t^\frac{ g + |c| }{2} \big[ \mf{H}^\frac{1}{2} | \hat{\phi} ( t, \xi ) | + \mf{H}^{ -\frac{1}{2} } | \partial_t \hat{\phi} ( t, \xi ) | \big] \text{,} \\
\notag| \pd{U_{Hz}}{m_H} ( t, \xi ) | &\gtrsim \zfac^\frac{ 1 + g - |c| }{2} t^\frac{ g - |c| }{2} \big[ \mf{H}^\frac{1}{2} | \hat{\phi} ( t, \xi ) | + \mf{H}^{ - \frac{1}{2} } | \partial_t \hat{\phi} ( t, \xi ) | \big] \text{.}
\end{align} 

\subsubsection{Loss of Regularity}

Applying Theorem \ref{thm.energy_main} to the preceding development results in the key asymptotics, scattering, and loss of regularity result for the anisotropic wave equations \eqref{eq.kasner}.
The precise statements of the above are summarized in the following theorem:

\begin{theorem} \label{thm.kasner_main}
Consider the setting of Assumption \ref{ass.kasner}.
\begin{itemize}
\item \emph{Asymptotics:} Given $\hat{\phi}_0, \hat{\phi}_1: \smash{ \R_\xi^d } \rightarrow \C$, there exists a unique solution $\hat{\phi}$ of \eqref{eq.kasnerf} satisfying $\smash{ ( \hat{\phi}, \partial_t \hat{\phi} ) ( T, \cdot ) = ( \hat{\phi}_0, \hat{\phi}_1 ) }$.
Moreover, letting $\hat{\varphi}_\pm$ be as in \eqref{eq.kasner_Upz_1} (if $g \neq 1$) or \eqref{eq.kasner_Upz_2} (if $g = 1$), then the following asymptotic limits are finite for all $\xi \in \R^d$:
\begin{equation}
\label{eq.kasner_asymp_limit} \begin{cases} \hat{\varphi}_{ +, 0 } ( \xi ) := \lim_{ \tau \searrow 0 } \hat{\varphi}_+ ( \tau, \xi ) \text{,} \quad \hat{\varphi}_{ -, 0 } ( \xi ) := \lim_{ \tau \searrow 0 } \hat{\varphi}_- ( \tau, \xi ) & \quad g \neq 1 \text{,} \\ \hat{\varphi}_{ +, 0 } ( \xi ) := \lim_{ \tau \searrow 0 } [ \hat{\varphi}_+ - ( \log z ) \hat{\varphi}_- ] ( \tau, \xi ) \text{,} \quad \hat{\varphi}_{ -, 0 } ( \xi ) := \lim_{ \tau \searrow 0 } \hat{\varphi}_- ( \tau, \xi ) & \quad g = 1 \text{.} \end{cases}
\end{equation}
In addition, the following estimate holds for any $\xi \in \R^{ d }$:
\begin{equation}
\label{eq.kasner_asymp} \zfac^\frac{ 1 - g }{2} | \hat{\varphi}_{ +, 0 } ( \xi ) | + \zfac^\frac{ g - 1 }{2} | \hat{\varphi}_{ -, 0 } ( \xi ) | \lesssim \zfac^\frac{ |c| }{2} \big[ \langle \xi \rangle^\frac{1}{2} | \hat{\phi} ( T, \xi ) | + \langle \xi \rangle^{ -\frac{1}{2} } | \partial_t \hat{\phi} ( T, \xi ) | \big] \text{.}
\end{equation}

\item \emph{Scattering:} Given $\hat{\varphi}_{ \pm, 0 }: \smash{ \R_\xi^d } \rightarrow \C$, there is a unique solution $\hat{\phi}$ of \eqref{eq.kasnerf} such that, defining $\hat{\varphi}_\pm$ as in \eqref{eq.kasner_Upz_1} ($g \neq 1$) or \eqref{eq.kasner_Upz_2} ($g = 1$), the following limits hold for all $\xi \in \R^d$:
\begin{equation}
\label{eq.kasner_scatter_limit} \begin{cases} \lim_{ \tau \searrow 0 } \hat{\varphi}_+ ( \tau, \xi ) = \hat{\varphi}_{ +, 0 } ( \xi ) \text{,} \quad \lim_{ \tau \searrow 0 } \hat{\varphi}_- ( \tau, \xi ) = \hat{\varphi}_{ -, 0 } ( \xi ) & \quad g \neq 1 \text{,} \\ \lim_{ \tau \searrow 0 } [ \hat{\varphi}_+ - ( \log z ) \hat{\varphi}_- ] ( \tau, \xi ) = \hat{\varphi}_{ +, 0 } ( \xi ) \text{,} \quad \lim_{ \tau \searrow 0 } \hat{\varphi}_- ( \tau, \xi ) = \hat{\varphi}_{ -, 0 } ( \xi ) & \quad g = 1 \text{.} \end{cases}
\end{equation}
In addition, the following estimate holds for any $\xi \in \R^{ d }$:
\begin{equation}
\label{eq.kasner_scatter} \langle \xi \rangle^\frac{1}{2} | \hat{\phi} ( T, \xi ) | + \langle \xi \rangle^{ -\frac{1}{2} } | \partial_t \hat{\phi} ( T, \xi ) | \lesssim \zfac^\frac{ |c| }{2} \big[ \zfac^\frac{ 1 - g }{2} | \hat{\varphi}_{ +, 0 } ( \xi ) | + \zfac^\frac{ g - 1 }{2} | \hat{\varphi}_{ -, 0 } ( \xi ) | \big] \text{.}
\end{equation}
\end{itemize}
\end{theorem}

\begin{proof}
That the limits $\hat{\varphi}_{ \pm, 0 }$ in \eqref{eq.kasner_asymp_limit} and \eqref{eq.kasner_scatter_limit} are finite is established in the same manner as in the proof of Theorem \ref{thm.wave_main}.
Thus, we focus on the estimates \eqref{eq.kasner_asymp} and \eqref{eq.kasner_scatter}.
Also, we only consider $( T, \xi ) \in Z_H$, as the case $( T, \xi ) \in Z_I$ holds trivially via the same argument as in Theorem \ref{thm.wave_main}.

For $( T, \xi ) \in Z_H$, using \eqref{eq.system_UA}, \eqref{eq.kasner_Uhz_ex}, and the observation that $\mf{H} ( T, \xi ) \simeq \langle \xi \rangle$, we obtain
\begin{align}
\label{eql.kasner_main_1} | U_A ( T, \xi ) | &\lesssim \zfac^\frac{ 1 + g + |c| }{2} \big[ \langle \xi \rangle^\frac{1}{2} | \hat{\phi} ( T, \xi ) | + \langle \xi \rangle^{ -\frac{1}{2} } | \partial_t \hat{\phi} ( T, \xi ) | \big] \text{,} \\
\notag| U_A ( T, \xi ) | &\gtrsim \zfac^\frac{ 1 + g - |c| }{2} \big[ \langle \xi \rangle^\frac{1}{2} | \hat{\phi} ( T, \xi ) | + \langle \xi \rangle^{ - \frac{1}{2} } | \partial_t \hat{\phi} ( T, \xi ) | \big] \text{.}
\end{align}
Moreover, by \eqref{eq.system_UA}, \eqref{eq.kasner_Upz_1}, and \eqref{eq.kasner_Upz_2}, we have (by the language of Theorem \ref{thm.energy_main})
\begin{equation}
\label{eql.kasner_main_2} | u_A ( \xi ) | \simeq \zfac | \hat{\varphi}_{ +, 0 } ( \xi ) | + \zfac^g | \hat{\varphi}_{ -, 0 } ( \xi ) | \text{.}
\end{equation}
Since \eqref{eq.energy_asymp} and \eqref{eq.energy_scatter} imply
\[
| U_A ( T, \xi ) | \simeq | u_A ( \xi ) | \text{,}
\]
then both \eqref{eq.kasner_asymp} and \eqref{eq.kasner_scatter} follow from \eqref{eql.kasner_main_1}--\eqref{eql.kasner_main_2} and the above.
\end{proof}

\begin{remark}
Note the estimates \eqref{eq.kasner_asymp} and \eqref{eq.kasner_scatter} now contain powers of $\zfac$.
These can be viewed as anisotropic weights, corresponding to different numbers of derivatives in different directions.
\end{remark}

\begin{remark}
Observe that in the special case of the scalar wave equation \eqref{eq.intro_wave_kasner} on Kasner backgrounds, i.e.\ $c := 0$ and $g = 1$, the estimates \eqref{eq.kasner_asymp} and \eqref{eq.kasner_scatter} reduce to the following:
\begin{equation}
\label{eq.kasner_ex} \langle \xi \rangle^\frac{1}{2} | \hat{\phi} ( T, \xi ) | + \langle \xi \rangle^{ -\frac{1}{2} } | \partial_t \hat{\phi} ( T, \xi ) | \simeq \big[ | \hat{\varphi}_{ +, 0 } ( \xi ) | + | \hat{\varphi}_{ -, 0 } ( \xi ) | \big] \text{.}
\end{equation}
In particular, \eqref{eq.kasner_ex} implies \eqref{eq.intro_kasner_est}, which was precisely the estimate obtained in \cite{li_2024} for \eqref{eq.intro_wave_kasner}.
\end{remark}

\subsection{Higher-Order Equations} \label{sec.wave_higher}

Finally, we briefly discuss how Section \ref{sec.system} can be applied to higher-order critically weakly hyperbolic and singular equations.
Our precise setting is as follows:

\begin{assumption} \label{ass.higher}
Fix $T > 0$, and consider the following $n$-th order equation on $( 0, T ]_t \times \R_x^d$:
\begin{equation}
\label{eq.higher} \partial_t^n \phi - \sum_{ j = 0 }^{ n - 1 } \sum_{ | \alpha | = 0 }^{ n - j } t^{ j - n + ( \ell + 1 ) | \alpha | } a_{ j, \alpha } (t) \partial_x^\alpha \partial_t^j \phi = 0 \text{,}
\end{equation}
where the unknown is $\phi: (0,T]_t \times \R_x^d \rightarrow \C$, where $\ell \in ( -1, \infty )$, and where the coefficients satisfy
\begin{equation}
\label{eq.higher_coeff} a_{ j, \alpha } \in C^\infty ( [ 0, T ]_t; \underbrace{ \C^d \otimes \dots \otimes \C^d }_{ | \alpha | \text{ times} } ) \text{,} \qquad 0 \leq j < n \text{,} \quad | \alpha | \leq n - j \text{.}
\end{equation}
\end{assumption}

\begin{remark}
In Assumption \ref{ass.higher}, we used $\alpha \in \N_0^d$ to denote $d$-dimensional multi-indices.
We will use standard notations regarding multi-indices throughout this subsection, e.g.
\[
| \alpha | := \alpha_1 + \dots + \alpha_d \text{,} \qquad \partial_x^\alpha := \partial_{ x_1 }^{ \alpha_1 } \dots \partial_{ x_d }^{ \alpha_d } \text{,} \qquad \xi^\alpha := \xi_1^{ \alpha_1 } \cdot \dots \cdot \xi_d^{ \alpha_d } \text{.}
\]
\end{remark}

It will be convenient to also define the following auxiliary quantities:

\begin{definition} \label{def.higher_aux}
Let $\mf{a}_{ j, k } \in C^\infty ( [ 0, T ]_t \times \Sph^{ d - 1 }_\omega; \C )$, for $0 \leq j < n$ and $0 \leq k \leq n - j$, be given by
\begin{equation}
\label{eq.higher_aux} \mf{a}_{ j, k } ( t, \omega ) := \sum_{ | \alpha | = k } a_{ j, \alpha } (t) \, \omega^\alpha \text{.}
\end{equation}
\end{definition}

Observe the Fourier transform of $\phi$ (when exists) satisfies the following on $( 0, T ]_t \times \R_\xi^d$: 
\begin{equation}
\label{eq.higherf} \partial_t^n \hat{\phi} - \sum_{ j = 0 }^{ n - 1 } A_j ( t, \xi ) \, t^{ j - n } \partial_t^j \hat{\phi} = 0 \text{,} \qquad A_j ( t, \xi ) := \sum_{ k = 0 }^{ n - j } ( \imath t^{ \ell + 1 } | \xi | )^k \mf{a}_{ j, k } \big( t, \tfrac{ \xi }{ | \xi | } \big) \text{.}
\end{equation}
In order for \eqref{eq.higher} and \eqref{eq.higherf} to be appropriately hyperbolic, we must also assume the following:

\begin{assumption} \label{ass.higher_hyp}
Suppose the polynomials $\mf{p}_{ t, \omega }: \C \rightarrow \C$, given by
\begin{equation}
\label{eq.higher_hyp} \mf{p}_{ t, \omega } ( \lambda ) := \lambda^n + \mf{a}_{ n - 1, 1 } ( t, \omega ) \, \lambda^{ n - 1 } + \dots + \mf{a}_{ 1, n - 1 } ( t, \omega ) \, \lambda + \mf{a}_{ 0, n } ( t, \omega ) \text{,}
\end{equation}
have real, distinct roots for all $( t, \omega ) \in [ 0, T ] \times \Sph^{ d - 1 }$.
\end{assumption}

The basic setup, in relation to Section \ref{sec.system}, is identical to that of Section \ref{sec.wave_basic}:
\begin{equation}
\label{eq.higher_z} \mf{H} := t^\ell | \xi | \text{,} \qquad \zfac \simeq \langle \xi \rangle^\frac{1}{ \ell + 1 } \text{,} \qquad z \simeq t \langle \xi \rangle^\frac{1}{ \ell + 1 } \text{.}
\end{equation}
As before, we fix $\rho_0 \ll_T 1$ and set $Z_P$, $Z_I$, $Z_H$ as in Definition \ref{def.zones}.
To write \eqref{eq.higherf} as a first-order system, we define our unknown $U: ( 0, T ]_t \times \R^d_\xi \rightarrow \C^n$ on $Z_P$ and $Z_H$ as
\begin{equation}
\label{eq.higher_U} U_j |_{ Z_P } := t^{ - ( n - j ) } \, \partial_t^{j-1} \hat{\phi} \text{,} \qquad U_j |_{ Z_H } := ( \imath t^\ell |\xi| )^{ n - j } \, \partial_t^{ j - 1 } \hat{\phi} \text{,} \qquad 1 \leq j \leq n \text{.}
\end{equation}

$U$ can then be defined on $Z_I$ by interpolating between its values on $Z_P$ and $Z_H$ in \eqref{eq.higher_U} using a cutoff function, as in Sections \ref{sec.wave_basic} and \ref{sec.wave_kasner}.
That Assumption \ref{ass.system_i} holds for our setting is then a trivial consequence of the fact that $z \simeq 1$ on $Z_I$.
As a result, we henceforth focus on demonstrating that our setup satisfies all the conditions in Assumptions \ref{ass.system_p} for $Z_P$ and \ref{ass.system_h} for $Z_H$.

\subsubsection{The Pseudodifferential Zone}

First, from \eqref{eq.higherf} and \eqref{eq.higher_U}, we see $U$ satisfies, on $Z_P$,
\begin{align*}
\partial_t U = \begin{bmatrix}
- (n-1) t^{-1} & t^{-1} & 0 & \cdots & 0 & 0 \\
0 & - (n-2) t^{-1} & t^{-1} & \cdots & 0 & 0 \\
\vdots & \vdots & \vdots & \ddots & \vdots & \vdots \\
0 & 0 & 0 & \cdots & - t^{-1} & t^{-1} \\
t^{-1} A_0 & t^{-1} A_1 & t^{-1} A_2 & \cdots & t^{-1} A_{ n - 2 } & t^{-1} A_{ n - 1 }
\end{bmatrix} U \text{.}
\end{align*}
In particular, by \eqref{eq.higherf}, we can write the above coefficients $\mc{A}$ on $Z_P$ as 
\begin{align}
\label{eq.higher_p_A} \mc{A} &= t^{-1} \mc{B}^\ast_P + \zfac \mc{S}^{ a_P } ( Z_P; \C^n \otimes \C^n ) \text{,} \\
\notag \mc{B}^\ast_P :\!&= \begin{bmatrix}
  -(n-1) & 1 & \cdots & 0 & 0 \\
	0 & -(n-2) & \cdots & 0 & 0 \\
	\vdots & \vdots & \ddots & \vdots & \vdots \\
	0 & 0 & \cdots & -1 & 1 \\
	\mf{a}_{ 0, 0 } \big( 0, \tfrac{ \xi }{ | \xi | } \big) & \mf{a}_{ 1, 0 } \big( 0, \tfrac{ \xi }{ | \xi | } \big) & \cdots & \mf{a}_{ n-2, 0 } \big( 0, \tfrac{ \xi }{ | \xi | } \big) & \mf{a}_{ n-1, 0 } \big( 0, \tfrac{ \xi }{ | \xi | } \big)
\end{bmatrix} \text{,}
\end{align}
where $a_P > -1$.
Note that $\mc{B}^\ast_P$ is a \emph{constant} matrix (and is $t$-independent in particular), hence one can find a constant matrix $M_P$ such that $M_P \mc{B}^\ast_P M_P^{-1}$ is in Jordan normal form.
Thus, by taking the above $M_P$, we can write $\mc{A}$ in the form \eqref{eq.wave_Ap}, and we hence conclude:

\begin{proposition} \label{thm.higher_p_ass}
Assumption \ref{ass.system_p} holds for the higher-order equation from Assumption \ref{ass.higher}.
In particular, the entries of $B_P$ are given by the eigenvalues of $\mc{B}^\ast_P$.
\end{proposition}

\subsubsection{The Hyperbolic Zone}

On $Z_H$, we again apply \eqref{eq.higherf} and \eqref{eq.higher_U} to obtain
\begin{align*}
\partial_{t} U = \begin{bmatrix}
	\ell ( n - 1 ) t^{-1} & \imath \mf{H} & 0 & \cdots & 0 & 0 \\
	0 & \ell ( n - 2 ) t^{-1} & \imath \mf{H} & \cdots & 0 & 0 \\
	\vdots & \vdots & \vdots & \ddots & \vdots & \vdots \\
  0 & 0 & 0 & \cdots & \ell t^{-1} & \imath \mf{H} \\
	\tfrac{ \imath \mf{H} \, A_0 }{ ( \imath t^{ \ell + 1 } | \xi | )^n } & \tfrac{ \imath \mf{H} \, A_1 }{ ( \imath t^{ \ell + 1 } | \xi | )^{ n - 1 } } & \tfrac{ \imath \mf{H} \, A_2 }{ ( \imath t^{ \ell + 1 } | \xi | )^{ n - 2 } } & \cdots & \tfrac{ \imath \mf{H} \, A_{ n - 2 } }{ ( \imath t^{ \ell + 1 } | \xi | )^2 } & \tfrac{ \imath \mf{H} \, A_{ n - 1 } }{ \imath t^{ \ell + 1 } | \xi | }
\end{bmatrix} U \text{.}
\end{align*}
Then, by \eqref{eq.higherf} and \eqref{eq.higher_z}, we can express the above coefficients $\mc{A}$ on $Z_H$ as
\begin{equation}
\label{eq.higher_h_A} \mc{A} = \imath \mf{H} \mc{D}^\ast_H + t^{-1} \mc{B}^\ast_H + \zfac \mc{S}^{ a_H } ( Z_H; \C^n \otimes \C^n ) \text{,}
\end{equation}
where $a_H < -1$, and where the matrices $\mc{D}^\ast_H$ and $\mc{B}^\ast_H$ are given by
\begin{align}
\label{eq.higher_h_DB} \mc{D}^\ast_H &:= \begin{bmatrix}
	 0 & 1 & \cdots & 0 & 0 \\
	 0 & 0 & \cdots & 0 & 0 \\
	 \vdots & \vdots & \ddots & \vdots & \vdots \\
	 0 & 0 & \cdots & 0 & 1 \\ 
	 \mf{a}_{ 0, n } \big( t, \tfrac{ \xi }{ | \xi | } \big) & \mf{a}_{ 1, n - 1 } \big( t, \tfrac{ \xi }{ | \xi | } \big) & \cdots & \mf{a}_{ n - 2, 2 } \big( t, \tfrac{ \xi }{ | \xi | } \big) & \mf{a}_{ n - 1, 1 } \big( t, \tfrac{ \xi }{ | \xi | } \big)
\end{bmatrix} \text{,} \\
\notag \mc{B}^\ast_H &:= \begin{bmatrix}
	\ell (n-1) & 0 & \cdots & 0 & 0 \\
	0 & \ell (n-2) & \cdots & 0 & 0 \\
	\vdots & \vdots & \ddots & \vdots & \vdots \\
	0 & 0 & \cdots & \ell & 0 \\
	\mf{a}_{ 0, n - 1 } \big( t, \tfrac{ \xi }{ | \xi | } \big) & \mf{a}_{ 1, n - 2 } \big( t, \tfrac{ \xi }{ | \xi | } \big) & \cdots & \mf{a}_{ n - 2, 1 } \big( t, \tfrac{ \xi }{ | \xi | } \big) & \mf{a}_{ n - 1, 0 } \big( t, \tfrac{ \xi }{ | \xi | } \big)
\end{bmatrix} \text{.}
\end{align}

To satisfy Assumption \ref{ass.system_h}, we first need to find a matrix $M_H$ that diagonalizes $\mc{D}^\ast_H$, and such that its resulting diagonal entries are everywhere real.
For this, we note that for any $\lambda \in \C$,
\[
\det ( \lambda I_n - \mc{D}^\ast_H ) = \mf{p}_{ t, \frac{ \xi }{ | \xi | } } ( \lambda ) \text{,}
\]
with $\mf{p}_{ t, \omega }$ defined as in \eqref{eq.higher_hyp}.
Since Assumption \ref{ass.higher_hyp} ensures the roots of the above are everywhere real and distinct, it follows there is a matrix-valued function $M_H$ such that $\smash{ M_H \mc{D}^\ast_H M_H^{-1} }$ is everywhere diagonal, with its diagonal entries being everywhere real and distinct.

Thus, taking the above $M_H$, we can express $\mc{A}$ in the form \eqref{eq.wave_Ah}, with $D_H$ being diagonal, and with everywhere real and distinct entries.
The compactness of $[ 0, T ] \times \Sph^{ d - 1 }$ then implies our system is semi-strictly hyperbolic (see Definition \ref{def.system_h_strict}).
Lastly, since $\mc{D}^\ast_H$ and $\mc{B}^\ast_H$ are both $| \xi |$-independent, then so is $M_H$, as well as $B_H = \smash{ M_H \mc{B}^\ast_H M_H^{-1} }$.
As a result, we obtain the following:

\begin{proposition} \label{thm.higher_h_ass}
Assumption \ref{ass.system_h} holds for the higher-order equation from Assumption \ref{ass.higher}.
Furthermore, the system for $U$ is semi-strictly hyperbolic.
\end{proposition}

\begin{remark}
In fact, $M_H^{-1}$ can be written as a Vandermonde matrix whose entries depend on the roots of the $\mf{p}_{ t, \omega }$'s.
Thus, $B_H$ depends on $n$, $\ell$, the roots of the $\mf{p}_{ t, \omega }$'s, and the values of $\mf{a}_{ 0, n - 1 }, \dots, \mf{a}_{ n - 1, 0 }$.
Furthermore, like in Section \ref{sec.wave_basic}, the assumptions of Proposition \ref{thm.energy_h_simple} hold here, hence the loss of regularity depends only on $n$, $\ell$, the roots of the $\mf{p}_{ 0, \omega }$'s, and the $\mf{a}_{ i, n - i }$'s at $t = 0$.
\end{remark}

\subsubsection{Loss of Regularity}

Since Assumptions \ref{ass.system_i}, \ref{ass.system_p}, \ref{ass.system_h} hold in our setting and the system is semi-strictly hyperbolic, we can now apply Theorem \ref{thm.energy_main} to derive precise asymptotics, scattering, and loss of regularity results for \eqref{eq.higher}.
As the general formulas become quite involved when $n > 2$, we omit the computations and leave the details to the interested reader.

\section{Linearized Einstein-Scalar Equations on Kasner Backgrounds} \label{sec.einstein}

In this section, we apply the analysis in Section \ref{sec.system} to derive energy estimates for the linearized Einstein-scalar system about Kasner spacetimes.
From this, we not only \emph{recover the asymptotics and scattering results of \cite{li_2024} for Kasner spacetimes whose exponents satisfy the subcriticality condition}, but we \emph{further extend the energy estimates to encompass all non-degenerate Kasner exponents}.

Relative to the applications in Section \ref{sec.wave}, studying the linearized Einstein equations leads to some new complications.
The first is that this system contains additional structures beyond being a degenerate-singular hyperbolic system---for instance, an invariance under gauge transformations, as well as additional constraint equations.
An important consequence of this is that we can no longer directly apply the main results of Section \ref{sec.system} (Theorems \ref{thm.energy_main} and \ref{thm.energy_ex}) as black boxes.
Instead, we must take apart the proofs of these theorems and make use of these additional structures within the linearized Einstein system during various intermediate steps of the proofs.

Yet another complication is that the linearized Einstein-scalar system fails to be semi-strictly hyperbolic, so we must apply the more general Theorem \ref{thm.energy_ex} rather than Theorem \ref{thm.energy_main}.
That we still have converse (i.e.\ reversible) asymptotics and scattering theories for this system is due to the special structure of the Fuchsian term within the hyperbolic zone.

\subsection{The Linearized Einstein-Scalar System}

In this subsection, we recall the linearization of the Einstein-scalar equations about Kasner backgrounds formulated in \cite{li_2024}.
For brevity, we merely state the system here; a full derivation of the system can be found in \cite[Section 4]{li_2024}.

\begin{remark}
We will work on $\R^d$ rather than the torus $\mathbb{T}^d$ as our spatial domain in order to remain within the setting of this article.
(An analogous analysis on $\mathbb{T}^d$ can be straightforwardly carried out using Fourier series rather than Fourier transforms.)
Similarly, we alter some notations from \cite{li_2024} to main consistency with our setup---e.g.\ we work with the \emph{negatives} of the Kasner exponents.
\end{remark}

Consider a fixed \emph{Kasner spacetime} with spatial slices $\R^d$ and with \emph{Kasner exponents}
\begin{equation}
\label{eq.kasner_exponent} ( p_1, \dots, p_d ) := ( -\ell_1, \dots, -\ell_d ) \text{,} \qquad \ell_1, \dots, \ell_d > -1 \text{.}
\end{equation}
Recall the metric for this spacetime is given by \eqref{eq.intro_kasner_metric}.
We also recall, as in \cite{li_2024}, some key geometric and matter quantities associated with this background (in a constant mean curvature foliation):
\begin{itemize}
\item The spatial metric $\mathring{g}$ and inverse $\mathring{g}^{-1}$, given in Cartesian coordinates by
\begin{align}
\label{eq.kasner_metric} \mathring{g}_{ij} := t^{ -2 \ell_i } \delta_{ij} \text{,} \qquad \mathring{g}^{ij} := t^{ 2 \ell_i } \delta_{ij} \text{,} \qquad 1 \leq i, j \leq d \text{.}
\end{align}

\item The second fundamental form $\mathring{k}$ (with one index raised), given in Cartesian coordinates by
\begin{equation}
\label{eq.kasner_sff} \mathring{k}_i{}^j := \ell_i t^{-1} \delta_{ij} \text{,} \qquad 1 \leq i, j \leq d \text{.}
\end{equation}

\item In addition, the lapse $\mathring{n}$ and shift $\mathring{X}$ (in Cartesian coordinates) have trivial values:
\begin{equation}
\label{eq.kasner_lapse_shift} \mathring{n} :\equiv 1 \text{,} \qquad \mathring{X}^i :\equiv 0 \text{,} \qquad 1 \leq i \leq d \text{.}
\end{equation}

\item The associated scalar field $\mathring{\phi}$ has the following values:
\begin{equation}
\label{eq.kasner_scalar} \mathring{\phi} = \ell_\phi \log t \text{,} \qquad \partial_t \mathring{\phi} = \ell_\phi t^{-1} \text{,} \qquad \ell_\phi \in \R \text{.}
\end{equation}
\end{itemize}

\begin{remark}
For the above to satisfy the Einstein-scalar equations, we also require the \emph{Kasner relations}
\begin{equation}
\label{eq.kasner_relation} \sum_{ i = 1 }^d \ell_i = -1 \text{,} \qquad \sum_{ i = 1 }^d \ell_i^2 + 2 \ell_\phi^2 = 1 \text{.}
\end{equation}
While \eqref{eq.kasner_relation} is physically fundamental, we will, however, not need to assume it for our analysis.
\end{remark}

\begin{remark}
To remain consistent with the style and conventions of the remainder of this article, \emph{we will avoid utilizing Einstein summation notation for indexed quantities}.
Instead, every summation over indices in our upcoming equations will be explicitly indicated with a ``$\Sigma$".
\end{remark}

\subsubsection{The Linearized System} \label{sec.einstein_system}

The unknowns of our linearized Einstein-scalar system are as follows:
\begin{itemize}
\item The linearized metric $\check{\eta}$ and second fundamental form $\check{\kappa}$ (with one index raised):
\begin{equation}
\label{eq.kasner_lmetric} \check{\eta}_i{}^j: ( 0, T ]_t \times \R^d_x \rightarrow \R \text{,} \qquad \check{\kappa}_i{}^j: ( 0, T ]_t \times \R^d_x \rightarrow \R \text{,} \qquad 1 \leq i, j \leq d \text{.}
\end{equation}

\item The linearized lapse $\check{\nu}$ and shift $\check{\chi}$:
\begin{equation}
\label{eq.kasner_llapse_lshift} \check{\nu}: ( 0, T ]_t \times \R^d_x \rightarrow \R \text{,} \qquad \check{\chi}^i: ( 0, T ]_t \times \R^d_x \rightarrow \R \text{,} \qquad 1 \leq i \leq d \text{.}
\end{equation}

\item The linearized scalar field $\check{\phi}$ and time derivative $\check{\psi}$:
\begin{equation}
\label{eq.kasner_lscalar} \check{\phi}, \, \check{\psi}: ( 0, T ]_t \times \R^d_x \rightarrow \R \text{.}
\end{equation}
\end{itemize}
(The indexed quantities above represent components of tensor fields in Cartesian coordinates.)

Let $\eta$, $\kappa$, $\nu$, $\chi$, $\phi$, $\psi$ denote the spatial Fourier transforms of $\check{\eta}$, $\check{\kappa}$, $\check{\nu}$, $\check{\chi}$, $\check{\phi}$, $\check{\psi}$, respectively.
Taking the equations in \cite[Proposition 4.7]{li_2024} (adapted from $\mathbb{T}^d$ to $\R^d$) and applying \eqref{eq.kasner_metric}--\eqref{eq.kasner_scalar}, we obtain:

\begin{assumption}[Fourier-Transformed Linearized Einstein-Scalar System] \label{ass.einstein}
Fix $T > 0$, and consider the following system on $( 0, T ]_t \times \R^d_\xi$ in terms of the above-mentioned unknowns:
\begin{itemize}
\item The \emph{evolution equations} satisfied by $\eta$, $\kappa$, $\phi$, $\psi$---for any $1 \leq i, j \leq d$:
\begin{align}
\label{eq.einstein_evolution} t \partial_t \eta_i{}^j &= \kappa_i{}^j + 2 ( \ell_j - \ell_i ) \eta_i{}^j + \ell_i \delta_i{}^j \nu - \tfrac{1}{2} \imath ( \xi_i \chi^j + t^{ 2 \ell_j - 2 \ell_i } \xi_j \chi^i ) \text{,} \\
\notag t \partial_t \kappa_i{}^j &= \sum_{ a = 1 }^d ( - t^{ 2 + 2 \ell_a } \xi_a^2 \eta_i{}^j - t^{ 2 + 2 \ell_j } \xi_i \xi_j \eta_a{}^a + t^{ 2 + 2 \ell_a } \xi_i \xi_a \eta_a{}^j + t^{ 2 + 2 \ell_j } \xi_j \xi_a \eta_i{}^a ) \\
\notag &\qquad + t^{ 2 + 2 \ell_j } \xi_i \xi_j \nu - \ell_i \delta_i{}^j \nu + \imath ( \ell_j - \ell_i ) \xi_i \chi^j \text{,} \\
\notag t \partial_t \phi &= \psi + \ell_\phi \nu \text{,} \\
\notag t \partial_t \psi &= - \sum_{ a = 1 }^d t^{ 2 + 2 \ell_a } \xi_a^2 \phi - \ell_\phi \nu \text{.}
\end{align}

\item The \emph{elliptic equation} satisfied by $\nu$:
\begin{equation}
\label{eq.einstein_elliptic} \bigg( 1 + \sum_{ a = 1 }^d t^{ 2 + 2 \ell_a } \xi_a^2 \bigg) \nu = 2 \sum_{ a, b = 1 }^d ( t^{ 2 + 2 \ell_a } \xi_a^2 \eta_b{}^b - t^{ 2 + 2 \ell_a } \xi_a \xi_b \eta_a{}^b ) \text{.}
\end{equation}

\item The \emph{constraint equations}---for any $1 \leq j \leq d$:
\begin{align}
\label{eq.einstein_constraint} \sum_{ a, b = 1 }^d ( t^{ 2 + 2 \ell_a } \xi_a^2 \eta_b{}^b - t^{ 2 + 2 \ell_a } \xi_a \xi_b \eta_a{}^b ) + \sum_{ a = 1 }^d \ell_a \kappa_a{}^a + 2 \ell_\phi \psi &= 0 \text{,} \\
\notag \imath \sum_{ a = 1 }^d [ \xi_a \kappa_j{}^a + ( \ell_a - \ell_j ) \xi_j \eta_a{}^a ] + 2 \imath \ell_\phi \xi_j \phi &= 0 \text{,} \\
\notag \imath \sum_{ a = 1 }^d ( t^{ 2 \ell_a } \xi_a \kappa_a{}^j + \ell_a t^{ 2 \ell_j } \xi_j \eta_a{}^a - 2 \ell_a t^{ 2 \ell_a } \xi_a \eta_a{}^j ) + 2 \imath \ell_\phi t^{ 2 \ell_j } \xi_j \phi &= 0 \text{.}
\end{align}

\item The \emph{symmetry conditions}---for any $1 \leq i, j \leq d$:
\begin{equation}
\label{eq.einstein_symmetry} t^{ 2 \ell_i } \eta_i{}^j = t^{ 2 \ell_j } \eta_j{}^i \text{,} \qquad t^{ 2 \ell_i } ( \kappa_i{}^j - 2 \ell_i \eta_i{}^j ) = t^{ 2 \ell_j } ( \kappa_j{}^i  - 2 \ell_j \eta_j{}^i ) \text{.}
\end{equation}
\end{itemize}
\end{assumption}

\subsubsection{Gauge Covariance}

Due to the (spatial) coordinate-independence of the Einstein-scalar equations, there is a residual gauge freedom for the linearized equations \eqref{eq.einstein_evolution}--\eqref{eq.einstein_symmetry}.
More specifically, following \cite[Lemma 4.3]{li_2024} (after taking a spatial Fourier transform), given any
\[
\beta^i: ( 0, T ]_t \times \R^d_\xi \rightarrow \C \text{,} \qquad 1 \leq i \leq d
\]
that are Fourier transforms of real-valued functions, the system \eqref{eq.einstein_evolution}--\eqref{eq.einstein_symmetry} is invariant under the following transformation (for all $1 \leq i, j \leq d$) induced by $\beta$:
\begin{align}
\label{eq.einstein_gauge} \eta_i{}^j &\mapsto \eta_i{}^j + \tfrac{1}{2} \imath \xi_i \beta^j + \tfrac{1}{2} \imath t^{ 2 \ell_j - 2 \ell_i } \xi_j \beta^i \text{,} \\
\notag \kappa_i{}^j &\mapsto \kappa_i{}^j + \imath ( \ell_i - \ell_j ) \xi_i \beta^j \text{,} \\
\notag \chi^j &\mapsto \chi^j - t \partial_t \beta^j \text{.}
\end{align}
In other words, we interpret any transformation under \eqref{eq.einstein_gauge} as representing the same solution.

For the present discussion, we will focus exclusively on two gauges (for detailed derivations of the gauge invariance and choices, see \cite[Proposition 4.2, Lemma 4.3]{li_2024}):

\begin{assumption}[Fourier-Transformed Linearized Einstein-Scalar Gauges] \label{ass.einstein_gauge}
In the context of Assumption \ref{ass.einstein}, we will always adopt two particular gauge choices:
\begin{itemize}
\item \emph{Zero shift gauge}: $\beta$ is chosen so that the linearized shift vanishes:
\begin{equation}
\label{eq.gauge_shift} \chi^i \equiv 0 \text{,} \qquad 1 \leq i \leq d \text{.}
\end{equation}

\item \emph{Spatial harmonic gauge}: $\beta$ is chosen so that the following relations hold:
\begin{equation}
\label{eq.gauge_harmonic} \sum_{ a = 1 }^d ( 2 \xi_a \eta_j{}^a - \xi_j \eta_a{}^a ) = 0 \text{,} \qquad 1 \leq j \leq d \text{.}
\end{equation}
In this gauge, we also have the following additional equations---for all $1 \leq j \leq d$,
\begin{align}
\label{eq.einstein_elliptic_ex} \sum_{ a = 1 }^d t^{ 2 \ell_a } \xi_a^2 \, \chi^j &= - \imath ( 1 + 2 \ell_j ) t^{ 2 \ell_j } \xi_j \nu - 4 \imath \sum_{ a = 1 }^d \ell_a t^{ 2 \ell_a } \xi_a \eta_a{}^j + 2 \imath t^{ 2 \ell_j } \xi_j \bigg( \sum_{ a = 1 }^d \ell_a \eta_a{}^a + 2 \ell_\phi \phi \bigg) \text{.}
\end{align}
\end{itemize}
\end{assumption}

\subsubsection{The Main Quantities}

We now begin connecting the equations from Assumptions \ref{ass.einstein} and \ref{ass.einstein_gauge} to the framework of Section \ref{sec.system}.
First, we set the degenerate hyperbolicity and related quantities:
\begin{equation}
\label{eq.einstein_z} \mf{H} := \bigg[ \sum_{ i = 1 }^d t^{ 2 \ell_i } \xi_i^2 \bigg]^\frac{1}{2} \text{,} \qquad \zfac := \max_{ 1 \leq i \leq d } \langle \xi_i \rangle^\frac{1}{ \ell_i + 1 } \text{,} \qquad \ell_\ast := \min ( \ell_1, \dots, \ell_d ) \text{.}
\end{equation}
We can then define the usual zones $Z_P$, $Z_I$, $Z_H$ as in Definition \ref{def.zones}, with some $\rho_0 \ll_T 1$.

For the following, it will be useful to consider the following renormalizations of $\eta$ and $\kappa$:
\begin{equation}
\label{eq.einstein_renorm} \bar{\eta}_i{}^j := t^{ \ell_i - \ell_j } \eta_i{}^j \text{,} \qquad \bar{\kappa}_i{}^j := t^{ \ell_i - \ell_j } \kappa_i{}^j \text{,} \qquad 1 \leq i, j \leq d \text{.}
\end{equation}
This can be viewed as expressing $\eta$ and $\kappa$ in terms of a $\mathring{g}$-orthonormal frame and coframe.

To obtain a system as in Assumption \ref{ass.system}, we define the unknown $U: ( 0, T ]_t \times \R^d_\xi \rightarrow \C^{ 2 d^2 + 2 }$ by
\begin{align}
\label{eq.einstein_U} U^T |_{ Z_P } &:= \begin{bmatrix} \bar{\eta}_1{}^1 & \bar{\kappa}_1{}^1 & \bar{\eta}_1{}^2 & \bar{\kappa}_1{}^2 & \dots & \bar{\eta}_d{}^d & \bar{\kappa}_d{}^d & \phi & \psi \end{bmatrix} \text{,} \\
\notag U^T |_{ Z_H } &:= \begin{bmatrix} \imath t \mf{H} \bar{\eta}_1{}^1 & \bar{\kappa}_1{}^1 & \dots & \imath t \mf{H} \bar{\eta}_d{}^d & \bar{\kappa}_d{}^d & \imath t \mf{H} \phi & \psi \end{bmatrix} \text{.}
\end{align}
Analogously to Section \ref{sec.wave}, one can define $U |_{ Z_I }$ by interpolating between the above values using a cutoff function.
A system of the form \eqref{eq.system_gen} can now be derived by taking the evolution equations \eqref{eq.einstein_evolution}.
To close this system, we replace every instance of $\nu$ on the right-hand side of \eqref{eq.einstein_evolution} using \eqref{eq.einstein_elliptic}, and we replace each $\chi$ using either \eqref{eq.gauge_shift} or \eqref{eq.einstein_elliptic_ex}, depending on the choice of gauge.

\begin{remark}
Note the first part of \eqref{eq.einstein_symmetry} implies $\bar{\eta}$ is symmetric, i.e.\ $\bar{\eta}_i{}^j = \bar{\eta}_j{}^i$.
Moreover, \eqref{eq.einstein_constraint} and \eqref{eq.einstein_symmetry} imply that various components of $U$ contain redundant information.
\end{remark}

\begin{remark}
Note the quantities in $U$ from \eqref{eq.einstein_U} are in fact analogues of $t$ times the quantities studied in Sections \ref{sec.wave_basic} and \ref{sec.wave_kasner}; this is to maintain consistency with the unknown quantities studied in \cite{li_2024}.
While this discrepancy slightly alters the ensuing analysis, it will not affect the final result.
\end{remark}

\subsection{The Pseudodifferential Zone} \label{sec.einstein_zp}

As in \cite{li_2024}, we adopt the zero shift gauge on $Z_P$.
The system, in this gauge and in terms of the the renormalized quantities, is then as follows:

\begin{proposition} \label{thm.einstein_low}
Suppose Assumptions \ref{ass.einstein} and \ref{ass.einstein_gauge} hold, and assume moreover the zero shift gauge \eqref{eq.gauge_shift} on $Z_P$.
Then, the following equations hold on $Z_P$:
\begin{align}
\label{eq.einstein_low} t \partial_t \bar{\eta}_i{}^j &= \bar{\kappa}_i{}^j - ( \ell_i - \ell_j ) \bar{\eta}_i{}^j + \ell_i \delta_i{}^j \nu \text{,} \\
\notag t \partial_t \bar{\kappa}_i{}^j &= ( \ell_i - \ell_j ) \bar{\kappa}_i{}^j + t^{ 2 + \ell_i + \ell_j } \xi_i \xi_j \nu - \ell_i \delta_i{}^j \nu - t^2 \mf{H}^2 \, \bar{\eta}_i{}^j \\
\notag &\qquad + \sum_{ a = 1 }^d ( t^{ 2 + \ell_i + \ell_a } \xi_i \xi_a \bar{\eta}_a{}^j + t^{ 2 + \ell_j + \ell_a } \xi_j \xi_a \bar{\eta}_i{}^a - t^{ 2 + \ell_i + \ell_j } \xi_i \xi_j \bar{\eta}_a{}^a ) \text{,} \\
\notag t \partial_t \phi &= \psi + \ell_\phi \nu \text{,} \\
\notag t \partial_t \psi &= - t^2 \mf{H}^2 \, \phi - \ell_\phi \nu \text{,} \\
\notag ( 1 + t^2 \mf{H}^2 ) \nu &= 2 \sum_{ a, b = 1 }^d ( t^{ 2 + 2 \ell_a } \xi_a^2 \bar{\eta}_b{}^b - t^{ 2 + \ell_a + \ell_b } \xi_a \xi_b \bar{\eta}_a{}^b ) \text{.}
\end{align}
\end{proposition}

\begin{proof}
These are immediate consequences of \eqref{eq.einstein_evolution}--\eqref{eq.einstein_elliptic}, \eqref{eq.gauge_shift}, and \eqref{eq.einstein_z}--\eqref{eq.einstein_renorm}.
\end{proof}

\subsubsection{Estimates on $Z_P$}

Next, we express our system in the form \eqref{eq.system_p}.
The following shows that the coefficients of the system satisfied by $U$ has a block-diagonal structure at leading order:

\begin{proposition} \label{thm.einstein_system_p}
Assume the setting of Proposition \ref{thm.einstein_low}.
Then, on $Z_P$:
\begin{itemize}
\item For any $1 \leq i, j \leq d$ with $\ell_i \neq \ell_j$, we have
\begin{align}
\label{eq.einstein_system_p1a} \partial_t \begin{bmatrix} \bar{\eta}_i{}^j + \frac{1}{ 2 ( \ell_j - \ell_i ) } \bar{\kappa}_i{}^j \\ \bar{\kappa}_i{}^j \end{bmatrix} &= t^{-1} \begin{bmatrix} \ell_j - \ell_i & 0 \\ 0 & \ell_i - \ell_j \end{bmatrix} \begin{bmatrix} \bar{\eta}_i{}^j + \frac{1}{ 2 ( \ell_j - \ell_i ) } \bar{\kappa}_i{}^j \\ \bar{\kappa}_i{}^j \end{bmatrix} \\
\notag &\qquad + \zfac \, O \big( z^{ -1 + 2 ( \ell_\ast + 1 ) } \big) \, U \text{.} 
\end{align}

\item For any $1 \leq i, j \leq d$ with $\ell_i = \ell_j$, we have
\begin{equation}
\label{eq.einstein_system_p1b} \partial_t \begin{bmatrix} \bar{\eta}_i{}^j \\ \bar{\kappa}_i{}^j \end{bmatrix} = t^{-1} \begin{bmatrix} 0 & 1 \\ 0 & 0 \end{bmatrix} \begin{bmatrix} \bar{\eta}_i{}^j \\ \bar{\kappa}_i{}^j \end{bmatrix} + \zfac \, O \big( z^{ -1 + 2 ( \ell_\ast + 1 ) } \big) \, U \text{.}
\end{equation}

\item The following equation also holds for $\phi$ and $\psi$:
\begin{equation}
\label{eq.einstein_system_p2} \partial_t \begin{bmatrix} \phi \\ \psi \end{bmatrix} = t^{-1} \begin{bmatrix} 0 & 1 \\ 0 & 0 \end{bmatrix} \begin{bmatrix} \phi \\ \psi \end{bmatrix} + \zfac \, O \big( z^{ -1 + 2 ( \ell_\ast + 1 ) } \big) \, U \text{.}
\end{equation}
\end{itemize}
\end{proposition}

\begin{proof}
By Proposition \ref{thm.degen_est}, we have on $Z_P$ that
\[
1 + t^2 \mf{H}^2 \simeq 1 \text{,}
\]
hence the last equation of \eqref{eq.einstein_low} yields
\begin{equation}
\label{eql.einstein_system_p_1} \nu = \sum_{ a, b = 1 }^d O ( z^{ 2 ( \ell_\ast + 1 ) } ) \, \bar{\eta}_a{}^b \text{.}
\end{equation}

In particular, by the first two parts of \eqref{eq.einstein_low} and \eqref{eql.einstein_system_p_1}, we have, for $1 \leq i, j \leq d$,
\[
\partial_t \left[ \begin{matrix} \bar{\eta}_i{}^j \\ \bar{\kappa}_i{}^j \end{matrix} \right] = t^{-1} \left[ \begin{matrix} \ell_j - \ell_i & 1 \\ 0 & \ell_i - \ell_j \end{matrix} \right] \left[ \begin{matrix} \bar{\eta}_i{}^j \\ \bar{\kappa}_i{}^j \end{matrix} \right] + \zfac \, O \big( z^{ -1 + 2 ( \ell_\ast + 1 ) } \big) \, U \text{.}
\]
Diagonalizing the leading term on the right-hand side of the above (using a constant $2 \times 2$ matrix) yields \eqref{eq.einstein_system_p1a}--\eqref{eq.einstein_system_p1b}.
One can similarly obtain \eqref{eq.einstein_system_p2} from the third and fourth parts of \eqref{eq.einstein_low}.
\end{proof}

Combining \eqref{eq.einstein_system_p1a}--\eqref{eq.einstein_system_p2} (and in effect constructing a diagonalizer $M_P \in \C^{ 2 d^2 + 2 } \otimes \C^{ 2 d^2 + 2 }$), we then obtain a system of the form \eqref{eq.system_p}, with $F_P \equiv 0$,
\begin{align}
\label{eq.einstein_p_U} \partial_z U_P &= \big[ z^{-1} B_P + O \big( z^{ -1 + 2 ( \ell_\ast + 1 ) } \big) \big] U_P \text{,} \\
\notag U_P^T :\!\!&= \begin{bmatrix} \bar{\Upsilon}_1{}^1 & \bar{\kappa}_1{}^1 & \dots & \bar{\Upsilon}_d{}^d & \bar{\kappa}_d{}^d & \phi & \psi \end{bmatrix} \text{,}
\end{align}
where $\bar{\Upsilon}$ in the right-hand side of the second part is given by
\begin{equation}
\label{eq.einstein_p_Ux} \bar{\Upsilon}_i{}^j := \begin{cases} \bar{\eta}_i{}^j & \ell_i = \ell_j \text{,} \\ \bar{\eta}_i{}^j + \tfrac{1}{ 2 ( \ell_j - \ell_i ) } \bar{\kappa}_i{}^j & \ell_i \neq \ell_j \text{,} \end{cases} \qquad 1 \leq i, j \leq d \text{,}
\end{equation}
and where $B_P \in \C^{ 2 d^2 + 2 } \otimes \C^{ 2 d^2 + 2 }$ is block-diagonal matrix of the form
\begin{equation}
\label{eq.einstein_p_B} B_P = \operatorname{diag} ( B_{ P; (11) }, B_{ P; (12) }, \dots, B_{ P; (dd) }, B_{ P; (0) } ) \text{,}
\end{equation}
with the quantities on the right-hand side given by
\begin{equation}
\label{eq.einstein_p_Bx} B_{ P; (0) } = \begin{bmatrix} 0 & 1 \\ 0 & 0 \end{bmatrix} \text{,} \qquad B_{ P; (ij) } = \begin{bmatrix} \ell_j - \ell_i & \delta_{ \ell_i \ell_j } \\ 0 & \ell_i - \ell_j \end{bmatrix} \text{,} \qquad 1 \leq i, j \leq d \text{,}
\end{equation}
with $\delta$ denoting the Kronecker delta.
Note $B_P$ fails to be diagonal but is in Jordan normal form.

\begin{remark}
In the above, the blocks $B_{ P; (ij) }$ correspond to the $( \bar{\eta}_i{}^j, \bar{\kappa}_i{}^j )$-components of $U$, while the block $B_{ P; (0) }$ corresponds to the remaining $( \phi, \psi )$-components of $U$.
\end{remark}

\begin{remark}
Note that \eqref{eq.einstein_p_U}--\eqref{eq.einstein_p_Bx} do not quite correspond to Assumption \ref{ass.system_p} being satisfied, since here we are performing not only a change of basis (given by $M_P$), but also a gauge transformation to the zero shift gauge.
Nonetheless, the analysis of Section \ref{sec.system_zp} still applies to this system.
\end{remark}

We can now apply the analysis of Section \ref{sec.system_zp} to the system \eqref{eq.einstein_p_U}--\eqref{eq.einstein_p_Bx}.
In particular, Proposition \ref{thm.energy_p} yields the asymptotic quantity $\pd{U_{Pz}}{m_P}$ (for sufficiently large $m_P \in \N_0$).
Moreover, \eqref{eq.system_Ep}, \eqref{eq.energy_p_U}, \eqref{eq.einstein_p_U}, and \eqref{eq.einstein_p_B} imply that $\pd{U_{Pz}}{m_P}$ is precisely given by
\begin{equation}
\label{eq.einstein_p_Upz} \pd{U_{Pz}}{m_P} = \mc{E}_P \pd{Q_P}{m_P} U_P \text{,}
\end{equation}
where $\pd{Q_P}{m_P}$ is bounded and invertible, and where $\mc{E}_P$ is block-diagonal and of the form
\begin{equation}
\label{eq.einstein_p_E} \mc{E}_P = \operatorname{diag} ( \mc{E}_{ P; (11) }, \mc{E}_{ P; (12) }, \dots, \mc{E}_{ P; (dd) }, \mc{E}_{ P; (0) } ) \text{,}
\end{equation}
with the quantities on the right-hand side given by
\begin{equation}
\notag \mc{E}_{ P; (0) } = \begin{bmatrix} 1 & - \log z \\ 0 & 1 \end{bmatrix} \text{,} \qquad \mc{E}_{ P; (ij) } = \begin{bmatrix} z^{ \ell_i - \ell_j } & - ( \log z ) \, \delta_{ \ell_i \ell_j } \\ 0 & z^{ \ell_j - \ell_i } \end{bmatrix} \text{,} \qquad 1 \leq i, j \leq d \text{,}
\end{equation}
Recalling Proposition \ref{thm.energy_p} again results in the following properties for $\pd{U_{Pz}}{m_P}$:

\begin{proposition} \label{thm.einstein_energy_p}
Suppose Assumptions \ref{ass.einstein} and \ref{ass.einstein_gauge} hold, and assume also that we have fixed the zero shift gauge \eqref{eq.gauge_shift} on $Z_P$.
Then, for sufficiently large $m_P \in \N_0$:
\begin{itemize}
\item The following limits are both well-defined and finite:
\begin{equation}
\label{eq.einstein_energy_p_asymp} \pd{U_{Pz}}{m_P} ( 0, \xi ) := \lim_{ \tau \searrow 0 } \pd{U_{Pz}}{m_P} ( \tau, \xi ) \text{,} \qquad \xi \in \R^d \text{.}
\end{equation}

\item For any $\xi \in \R^d$ and $0 \leq t_0, t_1 \leq T$ such that $( t_0, \xi ), ( t_1, \xi ) \in Z_P$, we have
\begin{equation}
\label{eq.einstein_energy_p} | \pd{U_{Pz}}{ m_P } ( t_1, \xi ) | \lesssim | \pd{U_{Pz}}{ m_P } ( t_0, \xi ) | \text{.}
\end{equation}
\end{itemize}
\end{proposition}

\subsubsection{Asymptotics on $Z_P$}

Next, one can ask how large $m_P$ must be for the conclusions of Proposition \ref{thm.einstein_energy_p} to hold.
Recall, from the remark below Proposition \ref{thm.energy_p}, that the key condition for $m_P$ is that $\smash{ \mc{E}_P \pd{R_P}{m_P} ( \pd{Q_P}{m_P} )^{-1} \mc{E}_P^{-1} }$ is $z$-integrable.
By a closer analysis of the error coefficients $R_P$ in \eqref{eq.einstein_p_U} (from inspecting \eqref{eq.einstein_low}), one can show, as in \cite{li_2024}, that $m_P := 0$ suffices (that is, $\smash{ \mc{E}_P R_P \mc{E}_P^{-1} }$ is $z$-integrable) whenever the Kasner exponents satisfy the \emph{subcriticality condition}:
\begin{equation}
\label{eq.kasner_subcritical} \max_{ \substack{ 1 \leq i, j, k \leq d \\ j \neq k } } ( \ell_i - \ell_j - \ell_k ) < 1 \text{.}
\end{equation}

In particular, when \eqref{eq.kasner_subcritical} holds, the following is bounded and has finite limits as $t \searrow 0$,
\begin{align}
\label{eq.einstein_energy_U0} ( \pd{U_{Pz}}{0} )^T &= \mc{E}_P U_P \\
\notag &= \begin{bmatrix} V_{11} \mid V_{12} \mid \dots \mid V_{nn} \mid V_0 \end{bmatrix} \text{,}
\end{align}
where the quantities on the right-hand side are given, for any $1 \leq i, j \leq d$, by
\begin{align}
\label{eq.einstein_energy_UV} V_{ij} &= \begin{cases} \begin{bmatrix} z^{ \ell_i - \ell_j } \big[ \bar{\eta}_i{}^j + \tfrac{1}{ 2 ( \ell_j - \ell_i ) } \bar{\kappa}_i{}^j \big] & z^{ \ell_j - \ell_i } \bar{\kappa}_i{}^j \end{bmatrix} & \quad \ell_i \neq \ell_j \text{,} \\ \begin{bmatrix} \bar{\eta}_i{}^j - ( \log z ) \, \bar{\kappa}_i{}^j & \bar{\kappa}_i{}^j \end{bmatrix} & \quad \ell_i = \ell_j \text{,} \end{cases} \\
\notag V_0 &= \begin{bmatrix} \phi - ( \log z ) \, \psi & \psi \end{bmatrix} \text{.}
\end{align}
From the above, we obtain a more explicit description of the asymptotic quantities:

\begin{corollary} \label{thm.einstein_energy_subc}
Assume the setting of Proposition \ref{thm.einstein_energy_p}, and suppose in addition \eqref{eq.kasner_subcritical} holds.
Then, the following quantities are uniformly bounded on $Z_P$ and have finite limits as $t \searrow 0$:
\begin{itemize}
\item $\kappa_i{}^j$, for all $1 \leq i, j \leq d$.

\item $\eta_i{}^j$ and $t^{ 2 \ell_i - 2 \ell_j } \eta_j{}^i$, for all $1 \leq i, j \leq d$ such that $\ell_j > \ell_i$.

\item $\eta_i{}^j - ( \log z ) \, \kappa_i{}^j$, for all $1 \leq i, j \leq d$ such that $\ell_i = \ell_j$.

\item $\phi - ( \log z ) \, \psi$ and $\psi$.
\end{itemize}
\end{corollary}

\begin{proof}
From the above discussion, Proposition \ref{thm.einstein_energy_p} yields that $\pd{U_{Pz}}{0}$, given in \eqref{eq.einstein_energy_U0}--\eqref{eq.einstein_energy_UV}, is uniformly bounded and has a finite limit as $t \searrow 0$.
The second part of \eqref{eq.einstein_energy_UV} then immediately implies boundedness and finite limits for $\phi - ( \log z ) \, \psi$ and $\psi$, while \eqref{eq.einstein_renorm} and the first part of \eqref{eq.einstein_energy_UV} yields boundedness and finite limits for $t^{ \ell_j - \ell_i } \bar{\kappa}_i{}^j = \kappa_i{}^j$ for every $1 \leq i, j \leq d$.

Fix now $1 \leq i, j \leq d$.
If $\ell_i = \ell_j$, then \eqref{eq.einstein_renorm} and the first part of \eqref{eq.einstein_energy_UV} also implies boundedness and finite limits for $\eta_i{}^j - ( \log z ) \, \kappa_i{}^j$.
On the other hand, if $\ell_j > \ell_i$, then \eqref{eq.einstein_symmetry} yields that
\[
\eta_i{}^j = \tfrac{1}{ 2 ( \ell_i - \ell_j ) } [ \kappa_i{}^j - t^{ 2 ( \ell_j - \ell_i ) } \kappa_j{}^i ] \text{,} \qquad t^{ 2 \ell_i - 2 \ell_j } \eta_j{}^i = \eta_i{}^j
\]
are bounded and has finite limits, completing the proof.
\end{proof}

\begin{remark}
Observe that whenever $\ell_i \neq \ell_j$, then \eqref{eq.einstein_symmetry} and \eqref{eq.einstein_p_Ux} imply
\[
2 ( \ell_i - \ell_j ) \bar{\Upsilon}_i{}^j = \bar{\kappa}_j{}^i \text{,}
\]
so solving \eqref{eq.einstein_p_U} only yields information on $\bar{\kappa}$.
From here, one then recovers $\bar{\eta}_i{}^j$ from \eqref{eq.einstein_symmetry}.
\end{remark}

The novel portion of our analysis, however, concerns the settings where the subcriticality condition \eqref{eq.kasner_subcritical} fails to hold, as these cases were in particular not treated in \cite{li_2024}.
Nonetheless, using the methods of Section \ref{sec.system_zp}, we can still generate asymptotic quantities $\pd{U_{Pz}}{m_P}$ as in Proposition \ref{thm.einstein_energy_p}, as long as we apply additional higher-order renormalizations to the unknowns.
Moreover, the process detailed in Lemma \ref{thm.perf_diag_p} and Proposition \ref{thm.energy_p} provides a systematic algorithm for computing these asymptotic quantities $\pd{U_{Pz}}{m_P}$ as expansions involving the entries of $U_P$.

\subsection{The Hyperbolic Zone} \label{sec.einstein_zh}

We now adopt the spatial harmonic gauge \eqref{eq.gauge_harmonic} on $Z_H$.
Similar to \cite{li_2024}, the importance of this gauge is that it reveals the hyperbolic nature of the system.

\begin{remark}
The estimates for our system on $Z_H$ can be obtained by directly adapting the analysis in \cite[Section 6]{li_2024}, which already applies to all Kasner exponents, from the torus $\mathbb{T}^d$ to $\R^d$.
However, for completeness, we derive these estimates here using the methods from Section \ref{sec.system_zh}.
\end{remark}

\subsubsection{The System on $Z_H$}

The first step is to rewrite the system on $Z_H$ in terms of the renormalized quantities \eqref{eq.einstein_renorm}, while also incorporating the spatial harmonic gauge.

\begin{proposition} \label{thm.einstein_high}
Suppose Assumptions \ref{ass.einstein} and \ref{ass.einstein_gauge} hold, and assume we have also fixed the spatial harmonic gauge \eqref{eq.gauge_harmonic} on $Z_H$.
Then, the following equations hold on $Z_H$:
\begin{itemize}
\item Evolution equations--for any $1 \leq i, j \leq d$:
\begin{align}
\label{eq.einstein_high} t \partial_t \bar{\eta}_i{}^j &= \bar{\kappa}_i{}^j - ( \ell_i - \ell_j ) \bar{\eta}_i{}^j + \ell_i \delta_i{}^j \nu - \tfrac{1}{2} \imath ( t^{ \ell_i - \ell_j } \xi_i \chi^j + t^{ \ell_j - \ell_i } \xi_j \chi^i ) \text{,} \\
\notag t \partial_t \bar{\kappa}_i{}^j &= - t^2 \mf{H}^2 \, \bar{\eta}_i{}^j + ( \ell_i - \ell_j ) \bar{\kappa}_i{}^j + t^{ 2 + \ell_i + \ell_j } \xi_i \xi_j \nu \\
\notag &\qquad - \ell_i \delta_i{}^j \nu + \imath ( \ell_j - \ell_i ) t^{ \ell_i - \ell_j } \xi_i \chi^j \text{,} \\
\notag t \partial_t \phi &= \psi + \ell_\phi \nu \text{,} \\
\notag t \partial_t \psi &= - t^2 \mf{H}^2 \, \phi - \ell_\phi \nu \text{.}
\end{align}

\item Constraint equations--for any $1 \leq j \leq d$:
\begin{align}
\label{eq.einstein_highc} \tfrac{1}{2} t^2 \mf{H}^2 \, \sum_{ a = 1 }^d \bar{\eta}_a{}^a + \sum_{ a = 1 }^d \ell_a \bar{\kappa}_a{}^a + 2 \ell_\phi \psi &= 0 \text{,} \\
\notag \imath \sum_{ a = 1 }^d [ t^{ \ell_a - \ell_j } \xi_a \bar{\kappa}_j{}^a + ( \ell_a - \ell_j ) \xi_j \bar{\eta}_a{}^a ] + 2 \imath \ell_\phi \xi_j \phi &= 0 \text{,} \\
\notag \imath \sum_{ a = 1 }^d ( t^{ \ell_a + \ell_j } \xi_a \bar{\kappa}_a{}^j + \ell_a t^{ 2 \ell_j } \xi_j \bar{\eta}_a{}^a - 2 \ell_a t^{ \ell_a + \ell_j } \xi_a \bar{\eta}_a{}^j ) + 2 \imath \ell_\phi t^{ 2 \ell_j } \xi_j \phi &= 0 \text{.}
\end{align}

\item Elliptic equations--for any $1 \leq j \leq d$:
\begin{align}
\label{eq.einstein_highe} ( 1 + t^2 \mf{H}^2 ) \nu &= t^2 \mf{H}^2 \, \sum_{ a = 1 }^d \bar{\eta}^a{}_a \text{,} \\
\notag ( 1 + t^2 \mf{H}^2 ) \nu &= - 2 \sum_{ a = 1 }^d \ell_a \bar{\kappa}_a{}^a - 4 \ell_\phi \psi \text{,} \\
\notag t^2 \mf{H}^2 \chi^j &= - \imath ( 1 + 2 \ell_j ) t^{ 2 + 2 \ell_j } \xi_j \nu - 4 \imath \sum_{ a = 1 }^d \ell_a t^{ 2 + \ell_j + \ell_a } \xi_a \bar{\eta}_a{}^j \\
\notag &\qquad + 2 \imath t^{ 2 + 2 \ell_j } \xi_j \bigg( \sum_{ a = 1 }^d \ell_a \bar{\eta}_a{}^a + 2 \ell_\phi \phi \bigg) \text{,} \\
\notag t^2 \mf{H}^2 \chi^j &= - \imath ( 1 + 2 \ell_j ) t^{ 2 + 2 \ell_j } \xi_j \nu - 2 \imath \sum_{ a = 1 }^d t^{ 2 + \ell_j + \ell_a } \xi_a \bar{\kappa}_a{}^j \text{.}
\end{align}
\end{itemize}
\end{proposition}

\begin{proof}
Both \eqref{eq.einstein_high} and \eqref{eq.einstein_highc} follow immediately from \eqref{eq.einstein_evolution}, \eqref{eq.einstein_constraint}, \eqref{eq.gauge_harmonic}, and \eqref{eq.einstein_renorm}.
Similarly, the equalities in \eqref{eq.einstein_highe} are immediate consequences of \eqref{eq.einstein_elliptic}, \eqref{eq.einstein_elliptic_ex}, \eqref{eq.einstein_renorm}, and \eqref{eq.einstein_highc}.
\end{proof}

Similar to $Z_P$, the system on $Z_H$ again has a block-diagonal structure at the top order:

\begin{proposition} \label{thm.einstein_system_h}
Assume the setting of Proposition \ref{thm.einstein_high}.
Then, on $Z_H$:
\begin{itemize}
\item The following relations hold for $\nu$:
\begin{equation}
\label{eq.einstein_system_h0a} \nu = \sum_{ a = 1 }^d O ( 1 ) \, \bar{\eta}_a{}^a \text{,} \qquad t^2 \mf{H}^2 \, \nu = \sum_{ a = 1 }^d O (1) \, \bar{\kappa}_a{}^a + O (1) \, \psi \text{.}
\end{equation}

\item The following relations hold for $\chi$, for every $1 \leq i, j \leq d$:
\begin{align}
\label{eq.einstein_system_h0b} t^{ \ell_i - \ell_j } \xi_i \chi^j &= \sum_{ a, b = 1 }^d O (1) \, \bar{\eta}_a{}^b + O (1) \, \phi \text{,} \qquad t^{ \ell_i - \ell_j } \xi_i \chi^j = \sum_{ a, b = 1 }^d O (1) \, \bar{\kappa}_a{}^b + O ( z^{ -2 ( \ell_\ast + 1 ) } ) \, \psi \text{.}
\end{align}

\item The following equation holds for any $1 \leq i, j \leq d$:
\begin{align}
\label{eq.einstein_system_h1} \partial_t \begin{bmatrix} \tfrac{1}{ \sqrt{2} } ( \imath t \mf{H} \bar{\eta}_i{}^j + \bar{\kappa}_i{}^j ) \\ \tfrac{1}{ \sqrt{2} } ( \imath t \mf{H} \bar{\eta}_i{}^j - \bar{\kappa}_i{}^j ) \end{bmatrix} &= \left( \imath \mf{H} \left[ \begin{matrix} 1 & 0 \\ 0 & -1 \end{matrix} \right] + t^{-1} B_{ H; (ij) } \right) \begin{bmatrix} \tfrac{1}{ \sqrt{2} } ( \imath t \mf{H} \bar{\eta}_i{}^j + \bar{\kappa}_i{}^j ) \\ \tfrac{1}{ \sqrt{2} } ( \imath t \mf{H} \bar{\eta}_i{}^j - \bar{\kappa}_i{}^j ) \end{bmatrix} \\
\notag &\qquad + \zfac \, O \big( z^{ -1 - ( \ell_\ast + 1 ) } \big) \, U + F_{ H; (ij) } \text{,} \\
\notag B_{ H; (ij) } :\!\!&= \tfrac{1}{2} \begin{bmatrix} 1 + \frac{ t \partial_t \mf{H} }{ \mf{H} } & 1 + \frac{ t \partial_t \mf{H} }{ \mf{H} } + 2 ( \ell_j - \ell_i ) \\ 1 + \frac{ t \partial_t \mf{H} }{ \mf{H} } + 2 ( \ell_j - \ell_i ) & 1 + \frac{ t \partial_t \mf{H} }{ \mf{H} } \end{bmatrix} \text{,} \\
\notag F_{ H; (ij) } :\!\!&= \tfrac{1}{ \sqrt{2} } \, t^{-1} \begin{bmatrix} 1 & 1 \\ 1 & -1 \end{bmatrix} \begin{bmatrix} \frac{1}{2} t \mf{H} ( t^{ \ell_i - \ell_j } \xi_i \chi^j + t^{ \ell_j - \ell_i } \xi_j \chi^i ) \\ t^{ 2 + \ell_i + \ell_j } \xi_i \xi_j \nu \end{bmatrix} \text{.}
\end{align}

\item The following equation holds for $\phi$ and $\psi$:
\begin{align}
\label{eq.einstein_system_h2} \partial_t \begin{bmatrix} \tfrac{1}{ \sqrt{2} } ( \imath t \mf{H} \phi + \psi ) \\ \tfrac{1}{ \sqrt{2} } ( \imath t \mf{H} \phi - \psi ) \end{bmatrix} &= \left( \imath \mf{H} \left[ \begin{matrix} 1 & 0 \\ 0 & -1 \end{matrix} \right] + t^{-1} B_{ H; (0) } \right) \begin{bmatrix} \tfrac{1}{ \sqrt{2} } ( \imath t \mf{H} \phi + \psi ) \\ \tfrac{1}{ \sqrt{2} } ( \imath t \mf{H} \phi - \psi ) \end{bmatrix} \\
\notag &\qquad + \zfac \, O \big( z^{ -1 - ( \ell_\ast + 1 ) } \big) \, U + F_{ H; (0) } \text{,} \\
\notag B_{ H; (0) } :\!\!&= \tfrac{1}{2} \begin{bmatrix} 1 + \frac{ t \partial_t \mf{H} }{ \mf{H} } & 1 + \frac{ t \partial_t \mf{H} }{ \mf{H} } \\ 1 + \frac{ t \partial_t \mf{H} }{ \mf{H} } & 1 + \frac{ t \partial_t \mf{H} }{ \mf{H} } \end{bmatrix} \text{,} \\
\notag F_{ H; (0) } :\!\!&= \tfrac{1}{ \sqrt{2} } \, t^{-1} \begin{bmatrix} 1 & 1 \\ 1 & -1 \end{bmatrix} \begin{bmatrix} 0 \\ 0 \end{bmatrix} = \begin{bmatrix} 0 \\ 0 \end{bmatrix} \text{.}
\end{align}
\end{itemize}
\end{proposition}

\begin{proof}
First, Proposition \ref{thm.degen_est} implies $t^2 \mf{H}^2 \gtrsim z^{ 2 ( \ell_\ast + 1 ) } \gtrsim 1$, which, once combined with \eqref{eq.einstein_highe}, yields \eqref{eq.einstein_system_h0a}--\eqref{eq.einstein_system_h0b}.
Next, from Proposition \ref{thm.degen_est}, \eqref{eq.einstein_high}, and \eqref{eq.einstein_system_h0a}--\eqref{eq.einstein_system_h0b}, we have, for $1 \leq i, j \leq d$:
\begin{align}
\label{eql.einstein_system_h_1} \partial_t \begin{bmatrix} \imath t \mf{H} \bar{\eta}_i{}^j \\ \bar{\kappa}_i{}^j \end{bmatrix} &= \left( \imath \mf{H} \begin{bmatrix} 0 & 1 \\ 1 & 0 \end{bmatrix} + t^{-1} \begin{bmatrix} 1 + \frac{ t \partial_t \mf{H} }{ \mf{H} } + \ell_j - \ell_i & 0 \\ 0 & \ell_i - \ell_j \end{bmatrix} \right) \begin{bmatrix} \imath t \mf{H} \bar{\eta}_i{}^j \\ \bar{\kappa}_i{}^j \end{bmatrix} \\
\notag &\qquad + \zfac \, O \big( z^{ -1 - ( \ell_\ast + 1 ) } \big) \, U + t^{-1} \begin{bmatrix} \frac{1}{2} t \mf{H} ( t^{ \ell_i - \ell_j } \xi_i \chi^j + t^{ \ell_j - \ell_i } \xi_j \chi^i ) \\ t^{ 2 + \ell_i + \ell_j } \xi_i \xi_j \nu \end{bmatrix} \text{,} \\
\notag \partial_t \begin{bmatrix} \imath t \mf{H} \phi \\ \psi \end{bmatrix} &= \left( \imath \mf{H} \begin{bmatrix} 0 & 1 \\ 1 & 0 \end{bmatrix} + t^{-1} \begin{bmatrix} 1 + \frac{ t \partial_t \mf{H} }{ \mf{H} } & 0 \\ 0 & 0 \end{bmatrix} \right) \begin{bmatrix} \imath t \mf{H} \phi \\ \psi \end{bmatrix} \\
\notag &\qquad + \zfac \, O \big( z^{ -1 - ( \ell_\ast + 1 ) } \big) \, U + t^{-1} \begin{bmatrix} 0 \\ 0 \end{bmatrix} \text{.}
\end{align}
(Here, the final term on the right-hand side of \eqref{eql.einstein_system_h_1} are those containing the elliptic quantities $\nu$, $\chi$ that cannot be treated as remainders.)
Now, both \eqref{eq.einstein_system_h1} and \eqref{eq.einstein_system_h2} follow immediately from \eqref{eql.einstein_system_h_1}, in particular by diagonalizing the leading coefficients in \eqref{eql.einstein_system_h_1}.
\end{proof}

In particular, from \eqref{eq.einstein_system_h1}--\eqref{eq.einstein_system_h2}, we obtain a system of the form \eqref{eq.system_h},
\begin{align}
\label{eq.einstein_h_U} \partial_z U_H &= \big[ \imath \mf{Z} D_H + z^{-1} B_H + O \big( z^{ -1 - ( \ell_\ast + 1 ) } \big) \big] U_H + F_H \text{,} \\
\notag U_H^T :\!\!&= \tfrac{1}{ \sqrt{2} } \begin{bmatrix} \imath t \mf{H} \bar{\eta}_1{}^1 + \bar{\kappa}_1{}^1 & \imath t \mf{H} \bar{\eta}_1{}^1 - \bar{\kappa}_1{}^1 & \dots & \imath t \mf{H} \bar{\eta}_d{}^d - \bar{\kappa}_d{}^d & \imath t \mf{H} \phi + \psi & \imath t \mf{H} \phi - \psi \end{bmatrix} \text{,} \\
\notag F_H^T :\!\!&= \begin{bmatrix} F_{ H; (11) }^T \mid \dots \mid F_{ H; (dd) }^T \mid F_{ H; (0) }^T \end{bmatrix} \text{,}
\end{align}
where the coefficients $D_H$ and $B_H$ are given by
\begin{align}
\label{eq.einstein_h_DB} D_H &= \operatorname{diag} \left( \begin{bmatrix} 1 & 0 \\ 0 & -1 \end{bmatrix}, \dots, \begin{bmatrix} 1 & 0 \\ 0 & -1 \end{bmatrix}, \begin{bmatrix} 1 & 0 \\ 0 & -1 \end{bmatrix} \right) \text{,} \\
\notag B_H &= \operatorname{diag} ( B_{ H; (11) }, \dots, B_{ H; (dd) }, B_{ H; (0) } ) \text{.}
\end{align}
Note that $D_H$ is purely diagonal, with only real-valued entries.

\begin{remark}
Again, \eqref{eq.einstein_h_U}--\eqref{eq.einstein_h_DB} do not quite correspond to Assumption \ref{ass.system_h} being satisfied, since here we also perform a gauge transformation to the spatial harmonic gauge.
\end{remark}

Finally, observe our system \eqref{eq.einstein_h_U} fails to be semi-strictly hyperbolic, since $D_H$ has only two distinct eigenvalues, each having multiplicity $d^2 + 1$.
In terms of Definition \ref{def.system_h_partition}, we see from \eqref{eq.einstein_h_DB} that the appropriate partition of speeds for \eqref{eq.einstein_h_U} has two elements and is given by
\begin{equation}
\label{eq.einstein_h_partition} \mc{G} := \{ \{ 1, 3, \dots, 2 d^2 + 1 \}, \{ 2, 4, \dots, 2 d^2 + 2 \} \} \text{,}
\end{equation}
corresponding to the eigenvalues $+1$ and $-1$, respectively.
Most importantly for our analysis, \emph{the $\mc{G}$-diagonal part of $B_H$} (see \eqref{eq.system_Ph_BR}) \emph{is diagonal, with the same element along each diagonal entry}:
\begin{equation}
\label{eq.einstein_h_BG} B_{ \mc{G} } := \operatorname{diag} \Big( \tfrac{1}{2} \big( 1 + \tfrac{ t \partial_t \mf{H} }{ \mf{H} } \big), \dots, \tfrac{1}{2} \big( 1 + \tfrac{ t \partial_t \mf{H} }{ \mf{H} } \big) \Big) \text{.}
\end{equation}

\subsubsection{Estimates on $Z_H$}

We now apply the higher-order diagonalization, as in Lemma \ref{thm.part_diag_h}, to the system \eqref{eq.einstein_h_U}, with $m := 1$.
In our specific case, the system \eqref{eq.part_diag_h} is given by
\begin{equation}
\label{eq.einstein_h_U1} \partial_z \pd{U_H}{1} = ( \imath \mf{Z} \, D_H + \pd{D_H}{1} ) \pd{U_H}{1} + \pd{R_H}{1} ( \pd{Q_H}{1} )^{-1} \pd{U_H}{1} + \zfac^{-1} \pd{F_H}{1} \text{,}
\end{equation}
with all the notations defined as in Lemma \ref{thm.part_diag_h}.
In particular, note that
\begin{equation}
\label{eq.einstein_h_Q1} \pd{Q_H}{1} = I + \pd{N_H}{1} = I + O ( z^{ - ( \ell_\ast + 1 ) } ) \text{,} \qquad \pd{R_H}{1} ( \pd{Q_H}{1} )^{-1} = O ( z^{ -1 - ( \ell_\ast + 1 ) } ) \text{.}
\end{equation}
Moreover, recall from \eqref{eq.part_diag_h_DD} that $\pd{D_H}{1}$ is $\mc{G}$-diagonal and is given by (see \eqref{eq.system_Ph_BR})
\begin{equation}
\label{eq.einstein_h_D1} \pd{D_H}{1} = t^{-1} B_{ \mc{G} } + R_{ \mc{G} } \text{,} \qquad R_{ \mc{G} } = O ( z^{ -1 + ( \ell_\ast + 1 ) } ) \text{.}
\end{equation}

From \eqref{eq.einstein_h_BG}, we see that the remaining quantities in Definition \ref{def.system_Ph} satisfy
\begin{align}
\label{eq.einstein_h_Eh} b_{ H, \pm } ( t, \xi ) &= \Big( \tfrac{1}{2} \log ( \tau \mf{H} ( \tau, \xi ) ) |_{ \tau = \zfac^{-1} \rho_0 }^{ \tau = t }, \dots, \tfrac{1}{2} \log ( \tau \mf{H} ( \tau, \xi ) ) |_{ \tau = \zfac^{-1} \rho_0 }^{ \tau = t } \Big) \text{,} \\
\notag \mc{E}_H^\pm ( t, \xi ) &= \operatorname{diag} \Big( \Big[ \tfrac{ \zfac^{-1} \rho_0 \, \mf{H} ( \zfac^{-1} \rho_0, \, \xi ) }{ t \mf{H} ( t, \xi ) } \Big]^\frac{1}{2}, \dots, \Big[ \tfrac{ \zfac^{-1} \rho_0 \, \mf{H} ( \zfac^{-1} \rho_0, \, \xi ) }{ t \mf{H} ( t, \xi ) } \Big]^\frac{1}{2} \Big) \text{.}
\end{align}
In particular, note that since $\mc{E}_H^\pm$ is diagonal, and with the same value in each diagonal entry, then multiplication by $\mc{E}_H^\pm$ is the same as a scalar multiplication by
\begin{equation}
\label{eq.einstein_h_eh} \mf{e}_H := \Big[ \tfrac{ \zfac^{-1} \rho_0 \, \mf{H} ( \zfac^{-1} \rho_0, \, \xi ) }{ t \mf{H} ( t, \xi ) } \Big]^\frac{1}{2} \simeq [ t \mf{H} ( t, \xi ) ]^{ -\frac{1}{2} } \text{.}
\end{equation}
As a result, the renormalized quantity on $Z_H$ is the same in both directions:
\begin{equation}
\label{eq.einstein_h_Uhz} \pd{U_{Hz, \pm}}{1} = \mf{e}_H \pd{Q_H}{1} U_H := \pd{U_{Hz}}{1} \text{,} \qquad \pd{F_{Hz, \pm}}{1} = \mf{e}_H \pd{Q_H}{1} F_H := \pd{F_{Hz}}{1} \text{.}
\end{equation}

\begin{proposition} \label{thm.einstein_energy_h}
Suppose Assumptions \ref{ass.einstein} and \ref{ass.einstein_gauge} hold, and assume the spatial harmonic gauge \eqref{eq.gauge_harmonic} on $Z_H$.
Then, the following estimate holds for all $( t_0, \xi ), ( t_1, \xi ) \in Z_H$:
\begin{align}
\label{eq.einstein_energy_h} | \pd{U_{Hz}}{1} ( t_1, \xi ) | &\lesssim | \pd{U_{Hz}}{1} ( t_0, \xi ) | \text{.}
\end{align}
\end{proposition}

\begin{proof}
For this, we proceed through the proof of Proposition \ref{thm.energy_Ph}, but with a few alterations at key steps.
(Unfortunately, we cannot apply Proposition \ref{thm.energy_Ph} directly as a black box, as we must treat the inhomogeneous term $\pd{F_{Hz}}{1}$ specially using the constraint equations.)

From \eqref{eql.energy_Ph_0} and the fact that $B_{ \mc{G} } = B_{ \mc{G}, \pm }$ (see \eqref{eq.einstein_h_BG} and Definition \ref{def.part_diag_decomp}), we have
\begin{align}
\label{eql.einstein_energy_h_0} \partial_z \pd{U_{Hz}}{1} &= ( i \mf{Z} D_H + \pd{S_D}{1} + \pd{S_R}{1} ) \pd{U_{Hz}}{1} + \zfac^{-1} \pd{F_{Hz}}{1} \text{,} \\
\notag \pd{S_D}{1} &= R_{ \mc{G} } = O ( z^{ -1 + ( \ell_\ast + 1 ) } ) \text{,} \\
\notag \pd{S_R}{1} &= \mc{E}_H^\pm \pd{R_H}{1} ( \pd{Q_H}{1} )^{-1} ( \mc{E}_H^\pm )^{-1} = O ( z^{ -1 + ( \ell_\ast + 1 ) } ) \text{.}
\end{align}
(For the last two parts of \eqref{eql.einstein_energy_h_0}, we used \eqref{eq.einstein_h_D1} and the key observation that multiplication by $\mc{E}_H^\pm$ and $( \mc{E}_H^\pm )^{-1}$ are the same as multiplying by the scalar factors $\mf{e}_H$ and $\mf{e}_H^{-1}$, from \eqref{eq.einstein_h_eh}.)
Multiplying \eqref{eql.einstein_energy_h_0} by $\pd{\bar{U}_{Hz}}{1}$ and proceeding along the proof of Proposition \ref{thm.energy_Ph} toward \eqref{eql.energy_Ph_1}, except that we do not estimate the dot product involving $\pd{F_{Hz}}{1}$, we then obtain
\begin{equation}
\label{eql.einstein_energy_h_1} \tfrac{1}{2} \partial_z ( | \pd{U_{Hz}}{1} |^2 ) = ( | \pd{S_D}{1} | + | \pd{S_R}{1} | ) | \pd{U_{Hz}}{1} |^2 + | \zfac^{-1} \pd{F_{Hz}}{1} \cdot \pd{U_{Hz}}{1} | \text{.}
\end{equation}
As $\pd{S_D}{1}$ and $\pd{S_R}{1}$ are $z$-integrable, the corresponding terms in \eqref{eql.einstein_energy_h_1} can be treated as in the proof of Proposition \ref{thm.energy_Ph}.
Thus, it remains only to treat the last term on the right-hand side of \eqref{eql.einstein_energy_h_1}.

First of all, by \eqref{eq.einstein_system_h0a}--\eqref{eq.einstein_system_h2} and \eqref{eq.einstein_h_U}, we can estimate
\begin{align}
\label{eql.einstein_energy_h_10} | \zfac^{-1} F_H | &\lesssim z^{-1} \, \sum_{ a, b = 1 }^d \big( t \mf{H} \, | t^{ \ell_a - \ell_b } \xi_a \chi^b | + | t^{ 2 + \ell_a + \ell_b } \xi_a \xi_b \nu | \big) \\
\notag &\lesssim \sum_{ a, b = 1 }^d O ( z^{-1} ) \, ( | t \mf{H} \bar{\eta}_a{}^b | + | \bar{\kappa}_a{}^b | ) + O ( z^{-1} ) \, ( | t \mf{H} \phi | + | \psi | ) \text{.}
\end{align}
Furthermore, we can expand the last term in \eqref{eql.einstein_energy_h_1} as
\begin{align*}
\zfac^{-1} \pd{F_{Hz}}{1} \cdot \pd{U_{Hz}}{1} &= \mf{e}_H^2 \, U_H^T ( \pd{Q_H}{1} )^T \cdot \pd{Q_H}{1} \big( \zfac^{-1} F_H \big) \\
&= O ( ( t \mf{H} )^{-1} ) \, U_H^T \big[ I_{ 2 d^2 + 2 } + O ( z^{ - ( \ell_\ast + 1 ) } ) \big] \zfac^{-1} F_H \\
&= O ( ( t \mf{H} )^{-1} ) \, U_H^T \big( \zfac^{-1} F_H \big) + O ( ( t \mf{H} )^{-1} ) \, U_H^T \, O ( z^{ - ( \ell_\ast + 1 ) } ) \, \zfac^{-1} F_H \text{.}
\end{align*}
Controlling the last term in the right-hand side using \eqref{eq.einstein_h_eh}--\eqref{eq.einstein_h_Uhz} and \eqref{eql.einstein_energy_h_10}, we then obtain
\begin{align}
\label{eql.einstein_energy_h_11} \zfac^{-1} \pd{F_{Hz}}{1} \cdot \pd{U_{Hz}}{1} &= O ( ( t \mf{H} )^{-1} ) \, U_H^T \big( \zfac^{-1} F_H \big) + O ( ( t \mf{H} )^{-1} ) \, O ( z^{ -1 - ( \ell_\ast + 1 ) } ) \, | U_H |^2 \\
\notag &= O ( ( t \mf{H} )^{-1} ) \, U_H^T \big( \zfac^{-1} F_H \big) + O ( z^{ -1 - ( \ell_\ast + 1 ) } ) \, | \pd{U_{Hz}}{1} |^2 \text{.}
\end{align}

It remains to control the first term on the right-hand side of \eqref{eql.einstein_energy_h_11}, for which we can expand
\begin{align}
\label{eql.einstein_energy_h_20} U_H^T \big( \zfac^{-1} F_H \big) &= I_1 + I_2 \text{,} \\
\notag I_1 :\!\!&= z^{-1} \sum_{ i, j = 1 }^d \imath t \mf{H} \bar{\eta}_i{}^j \big[ \tfrac{1}{2} t \mf{H} ( t^{ \ell_i - \ell_j } \xi_i \chi^j + t^{ \ell_j - \ell_i } \xi_j \chi^i ) \big] \text{,} \\
\notag I_2 :\!\!&= z^{-1} \sum_{ i, j = 1 }^d \bar{\kappa}_i{}^j ( t^{ 2 + \ell_i + \ell_j } \xi_i \xi_j \nu ) \text{.}
\end{align}
For $I_1$, we first use the symmetry \eqref{eq.einstein_symmetry} of $\bar{\eta}$, along with \eqref{eq.gauge_harmonic}:
\begin{align*}
I_1 &= \imath z^{-1} ( t \mf{H} )^2 \sum_{ i, j = 1 }^d \xi_j t^{ \ell_j - \ell_i } \bar{\eta}_i{}^j \chi^i \\
&= \tfrac{1}{2} \imath z^{-1} ( t \mf{H} )^2 \sum_{ j = 1 }^d \bar{\eta}_j{}^j \sum_{ i = 1 }^d t^{ \ell_j - \ell_i } \xi_i \chi^i \text{.}
\end{align*}
Applying the first part of \eqref{eq.einstein_highc} and the first part of \eqref{eq.einstein_system_h0b} to the above yields
\begin{align}
\label{eql.einstein_energy_h_21} I_1 &= O ( z^{-1} ) \bigg[ \sum_{ j = 1 }^d O (1) \, \bar{\kappa}_j{}^j + O (1) \, \psi \bigg] \bigg[ \sum_{ a, b = 1 }^d O (1) \, \bar{\eta}_a{}^b + O (1) \, \phi \bigg] \\
\notag &= O ( z^{ -1 - ( \ell_\ast + 1 ) } ) \, | U_H |^2 \text{.}
\end{align}
Similarly, for $I_2$, we apply the second part of \eqref{eq.einstein_highc} and the second part of \eqref{eq.einstein_system_h0a}:
\begin{align}
\label{eql.einstein_energy_h_22} I_2 &= z^{-1} \sum_{ i = 1 }^d \bigg( \sum_{ j = 1 }^d t^{ \ell_j - \ell_i } \xi_j \bar{\kappa}_i{}^j \bigg) ( t^{ 2 + 2 \ell_i } \xi_i \nu ) \\
\notag &= O ( z^{-1} ) \, \sum_{ i = 1 }^d \bigg[ \sum_{ a = 1 }^d O (1) \, \bar{\eta}_a{}^a + O (1) \, \phi \bigg] \xi_i \, ( t^{ 2 + 2 \ell_i } \xi_i \nu ) \\
\notag &= O ( z^{-1} ) \, \bigg[ \sum_{ a = 1 }^d O (1) \, \bar{\eta}_a{}^a + O (1) \, \phi \bigg] \bigg[ \sum_{ b = 1 }^d O (1) \, \bar{\kappa}_b{}^b + O (1) \, \psi \bigg] \\
\notag &= O ( z^{ -1 - ( \ell_\ast + 1 ) } ) \, | U_H |^2 \text{.}
\end{align}

Finally, combining \eqref{eq.einstein_h_eh}--\eqref{eq.einstein_h_Uhz} and \eqref{eql.einstein_energy_h_11}--\eqref{eql.einstein_energy_h_22}, we obtain
\begin{align*}
\zfac^{-1} \pd{F_{Hz}}{1} \cdot \pd{U_{Hz}}{1} &= O ( ( t \mf{H} )^{-1} ) \, O ( z^{ -1 - ( \ell_\ast + 1 ) } ) \, | U_H |^2 + O ( z^{ -1 - ( \ell_\ast + 1 ) } ) \, | \pd{U_{Hz}}{1} |^2 \\
&= O ( z^{ -1 - ( \ell_\ast + 1 ) } ) \, | \pd{U_{Hz}}{1} |^2 \text{.}
\end{align*}
Applying \eqref{eql.einstein_energy_h_0}, \eqref{eql.einstein_energy_h_1}, and the above yields
\[
\tfrac{1}{2} \partial_z ( | \pd{U_{Hz}}{1} |^2 ) = O ( z^{ -1 - ( \ell_\ast + 1 ) } ) \, | \pd{U_{Hz}}{1} |^2 \text{,}
\]
which immediately leads to the desired uniform bound \eqref{eq.einstein_energy_h}.
\end{proof}

We can also restate Proposition \ref{thm.einstein_energy_h} in terms of the basic quantities:

\begin{corollary} \label{thm.einstein_energy_hx}
Assume the setting of Proposition \ref{thm.einstein_energy_h}.
Then, for all $( t_0, \xi ), ( t_1, \xi ) \in Z_H$,
\begin{align}
\label{eq.einstein_energy_hx} &\sum_{ i, j = 1 }^d \big[ ( t \mf{H} )^\frac{1}{2} | \bar{\eta}_i{}^j ( t_1, \xi ) | + ( t \mf{H} )^{ - \frac{1}{2} } | \bar{\kappa}_i{}^j ( t_1, \xi ) | \big] + \big[ ( t \mf{H} )^\frac{1}{2} | \phi ( t_1, \xi ) | + ( t \mf{H} )^{ - \frac{1}{2} } | \psi ( t_1, \xi ) | \big] \\
\notag &\quad \lesssim \sum_{ i, j = 1 }^d \big[ ( t \mf{H} )^\frac{1}{2} | \bar{\eta}_i{}^j ( t_0, \xi ) | + ( t \mf{H} )^{ - \frac{1}{2} } | \bar{\kappa}_i{}^j ( t_0, \xi ) | \big] + \big[ ( t \mf{H} )^\frac{1}{2} | \phi ( t_0, \xi ) | + ( t \mf{H} )^{ - \frac{1}{2} } | \psi ( t_0, \xi ) | \big] \text{.}
\end{align}
\end{corollary}

\begin{proof}
This follows immediately from \eqref{eq.einstein_energy_h}, after expanding $\pd{U_{Hz}}{1}$ using \eqref{eq.einstein_h_eh} and \eqref{eq.einstein_h_Uhz}.
\end{proof}

\subsection{Conclusions} \label{sec.einstein_est}

Propositions \ref{thm.einstein_energy_p} and \ref{thm.einstein_energy_h} yield the key energy estimates for our unknowns $\eta$, $\kappa$, $\phi$, $\psi$ (or rather, their renormalizations) on both $Z_P$ and $Z_H$.
Here, we briefly discuss how this would connect to a more extensive result for the linearized Einstein-scalar system.

\subsubsection{Global Estimates}

First, one can derive a corresponding estimate on $Z_I$ in the same manner as in Section \ref{sec.wave_basic}.
Fixing any (single) choice of gauge, defining $U |_{ Z_I }$ by interpolating between the formulas in \eqref{eq.einstein_U}, and using that $z \simeq 1$ on $Z_I$, one then obtains
\begin{equation}
\label{eq.einstein_energy_i} | U ( t_1, \xi ) | \lesssim | U ( t_0, \xi ) | \text{,} \qquad ( t_0, \xi ), ( t_1, \xi ) \in Z_I \text{.}
\end{equation}
Combining Proposition \ref{thm.einstein_energy_p}, Proposition \ref{thm.einstein_energy_h}, and \eqref{eq.einstein_energy_i} results in a global estimate:

\begin{proposition} \label{thm.einstein_energy}
Suppose Assumptions \ref{ass.einstein} and \ref{ass.einstein_gauge} hold, and assume the gauges \eqref{eq.gauge_shift} on $Z_P$ and \eqref{eq.gauge_harmonic} on $Z_H$.
In addition, let $U$ represent a solution to the linearized Einstein-scalar system from Assumption \ref{ass.einstein}, and define the renormalized unknown $\smash{ U_A: ( 0, T ]_t \times \R^d_\xi \rightarrow \C^{ 2 d^2 + 2 } }$ by
\begin{equation}
\label{eq.einstein_UA} U_A |_{ Z_P \setminus Z_I } := \pd{U_{Pz}}{m_P} \text{,} \qquad U_A |_{ Z_I } := U \text{,} \qquad U_A |_{ Z_H \setminus Z_I } := \pd{U_{Hz}}{1} \text{,}
\end{equation}
with $\pd{U_{Pz}}{m_P}$, $\pd{U_{Hz}}{1}$ as in Propositions \ref{thm.einstein_energy_p} and \ref{thm.einstein_energy_h}, respectively.
Then, the following estimate holds:
\begin{equation}
\label{eq.einstein_energy} | U_A ( t_1, \xi ) | \lesssim | U_A ( t_0, \xi ) | \text{,} \qquad 0 \leq t_0, t_1 < T \text{,} \quad \xi \in \R^d \text{.}
\end{equation}
\end{proposition}

\begin{remark}
$U_A$ assumes an inherently microlocal gauge that is simultaneously zero shift on $Z_P$ and spatial harmonic on $Z_H$.
However, from \eqref{eq.einstein_energy}, one could obtain estimates for $\eta$, $\kappa$, $\phi$, $\psi$ in other gauges by deriving estimates for the appropriate gauge transformations $\beta$; see \cite{li_2024}.
\end{remark}

\begin{remark}
When the subcriticality condition \eqref{eq.kasner_subcritical} holds (so that one can take $m_P = 0$), one can show that \eqref{eq.einstein_energy} expands into the energy estimates derived in \cite{li_2024}.
\end{remark}

\subsubsection{Asymptotics}

Recall that Proposition \ref{thm.einstein_energy_p} also yields quantities with asymptotic limits when solving the system backwards in time.
This leads to the following extension of Proposition \ref{thm.einstein_energy}:

\begin{proposition} \label{thm.einstein_asymp}
Assume the setting of Proposition \ref{thm.einstein_energy}.
Then, the following limits are finite:
\begin{equation}
\label{eq.einstein_asymp} u_A ( \xi ) := \lim_{ \tau \searrow 0 } U_A ( \tau, \xi ) \text{,} \qquad \xi \in \R^d \text{.}
\end{equation}
\end{proposition}

In particular, \eqref{eq.einstein_asymp} allows one to read off the leading asymptotic behaviors of $\eta$, $\kappa$, $\phi$, $\psi$ as $t \searrow 0$.
For Kasner backgrounds satisfying the subcriticality condition \eqref{eq.kasner_subcritical}, this is given by the conclusions of Corollary \ref{thm.einstein_energy_subc}, which match the results from \cite{li_2024}.

More generally, from \eqref{eq.einstein_p_Upz}, \eqref{eq.einstein_UA}, and \eqref{eq.einstein_asymp}, we have on $Z_P$ that
\begin{align*}
U_P ( t, \xi ) &\sim ( \pd{Q_P}{m_P} )^{-1} ( t, \xi ) \mc{E}_P^{-1} ( t, \xi ) u_A ( \xi ) \text{,}
\end{align*}
which, when combined with \eqref{eq.einstein_symmetry}, again unwinds to provide leading behaviors for $\eta$, $\kappa$, $\phi$, $\psi$.
A closer analysis would then yield more detailed asymptotics for these quantities.

\begin{remark}
Again, to derive corresponding asymptotics in other gauges, one would also need to derive asymptotic properties for the associated gauge transformations.
\end{remark}

\subsubsection{Scattering}

Proposition \ref{thm.einstein_energy} also provides the key evolutionary estimates needed for a scattering result from $t \searrow 0$.
However, for a full treatment of the scattering theory, one also needs to address the remaining components of the linearized Einstein-scalar system.

For instance, one would need to investigate asymptotic characterizations of the constraint equations, as well as gauge choices and transformations, in terms of higher-order renormalized quantities.
Another closely related issue would be the propagation of the appropriate asymptotic linearized constraint equations from $t \searrow 0$.
In particular, one should determine if any additional conditions need to be imposed in order for these characterizations to be valid.

These issues were addressed in \cite{li_2024} for Kasner backgrounds satisfying the subcriticality condition \eqref{eq.kasner_subcritical} (where no higher-order renormalizations were necessary).
We defer a more detailed study of Kasner backgrounds violating \eqref{eq.kasner_subcritical} to a future paper.

\end{document}